\documentclass[11pt]{article}
\usepackage{amsfonts}
\usepackage{epsfig}

\def\be{\begin{equation}}
\def\ee{\end{equation}}
\def\bea{\begin{eqnarray}}
\def\eea{\end{eqnarray}}
\def\bes{\begin{eqnarray*}}
\def\ees{\end{eqnarray*}}
\def\pmatrix{\left(\begin{array}{cc}}
\def\endpmatrix{\end{array}\right)}

\def\sgn{{\rm sgn}}

\def\Gr{{\rm Gr}}
\def\Lag {{\rm Lag}}
\def\Gr{{\rm Gr}}
\def\sign{{\rm sign}}
\def\det{{\rm det}}

\def\exp {{\rm exp}}

\def\Sp{{\rm Sp}}
\def\sf{{\rm sf}}

\def\nn{\nonumber}
\def\<{\langle}
\def\>{\rangle}
\def\lb{\label}

\def\sp{Sp}

\def\R{{\bf R}}
\def\C{{\bf C}}
\def\Z{{\bf Z}}
\def\N{{\bf N}}
\def\U{{\bf U}}

\def\T{{\bf T}}

\def\aa{{\alpha}}

\def\ga{{\gamma}}
\def\Ga{{\Gamma}}

\def\om{{\omega}}
\def\Om{{\Omega}}

\def\var{{\varepsilon}}
\def\lm{{\lambda}}
\def\Lm{{\Lambda}}

\def\sg{{\sigma}}
\def\Sg{{\Sigma}}

\def\J{{\cal J}}

\def\I{{\cal I}}

\def\td#1{\tilde{#1}}
\def\hb{\vrule height0.18cm width0.14cm $\,$}

\def\td#1{\tilde{#1}}
\def\diag{{\rm diag}}

\title{ Minimal period problems for brake orbits of nonlinear
autonomous reversible semipositive Hamiltonian systems}
\author{ Duanzhi Zhang
\thanks{Partially supported by National Science Foundation of China
(10801078, 11171314), LPMC of Nankai University.
E-mail: zhangdz@nankai.edu.cn}\\ \\
School of Mathematical Sciences and LPMC, Nankai University\\
Tianjin 300071, People's Republic of China}
\date{}
\begin{document}
\maketitle  \begin{abstract} In this paper,  for any positive
integer $n$, we study the Maslov-type index theory of $i_{L_0}$,
$i_{L_1}$ and $i_{\sqrt{-1}}^{L_0}$ with $L_0=\{0\}\times
\R^n\subset \R^{2n}$ and $L_1=\R^n\times \{0\} \subset \R^{2n}$. As
applications we study the minimal period problems for brake orbits
of nonlinear autonomous reversible Hamiltonian systems. For first
order nonlinear autonomous reversible Hamiltonian systems in
$\R^{2n}$, which are semipositive, and superquadratic at zero and
infinity£¬ we prove that for any $T>0$, the considered Hamiltonian
systems possesses a nonconstant $T$ periodic brake orbit $X_T$ with
minimal period no less than $\frac{T}{2n+2}$. Furthermore if
$\int_0^T H''_{22}(x_T(t))dt$ is positive definite, then the minimal
period of $x_T$ belongs to $\{T,\;\frac{T}{2}\}$. Moreover, if the
Hamiltonian system is even, we prove that for any $T>0$, the
considered even semipositive Hamiltonian systems possesses a
nonconstant symmetric brake orbit with minimal period
 belonging to $\{T,\;\frac{T}{3}\}$.

\end {abstract}
\noindent {\bf Key
words:} {symmetric, brake orbit, semipositive and reversible,
Maslov-type index, minimal period, Hamiltonian systems.}

\renewcommand{\theequation}{\thesection.\arabic{equation}}

\setcounter{equation}{0}
\section{Introduction and main results }

 In this paper, let
 $J=\left(\begin{array}{cc}0&-I_n\\I_n&0\end{array}\right)$ and $N=\left(\begin{array}{cc}-I_n&0\\0&I_n\end{array}\right)$,
 where $I_n$ is the identity in $\R^n$ and $n\in \N$. We suppose the following condition

 (H1) $H\in
 C^2(\R^{2n},\R)$ and satisfies the following  reversible condition
 \bea H(Nx)=H(x),\qquad \forall x\in \R^{2n}.\nn\eea

We consider the following  problem:
      \bea &&\dot{x}=JH'(x),\qquad x\in \R^{2n},\lb{1.1}\\
            &&x(-t)=Nx(t),\;\;x(T+t)=x(t),\qquad\forall t\in\R.\lb{1.2} \eea

A solution $(T,x)$ of (\ref{1.1})-(\ref{1.2}) is a special periodic
solution of the Hamiltonian system (\ref{1.1}). We call it a {\it
brake orbit} and $T$ the period of $x$. Moreover, if $x(\R)=-x(\R)$,
we call it a {\it symmetric brake orbit}. It is easy to check that
if $\tau$ is the minimal period of $x$, there must holds
$x(t+\frac{\tau}{2})=-x(t)$ for all $t\in\R$.

Since 1948, when H. Seifert in \cite{Se1} proposed his famous
conjecture of the existence of $n$ geometrically different brake
orbits in the potential well in $\R^n$ under certain conditions,
many people began to study this conjecture and related problems. Let
$^{\#}\td{\mathcal{O}}(\Om)$ and $^{\#}\td{\mathcal{J}}_b(\Sg)$ the
number of geometrically distinct brake obits in $\Om$ for the second
order case and on $\Sg$ for the first order case respectively.
 S.
Bolotin proved first in \cite{Bol}(also see \cite{BolZ}) of 1978 the
existence of brake orbits in general setting.
 K. Hayashi in \cite {Ha1}, H. Gluck and W. Ziller
 in \cite{GZ1}, and V. Benci in \cite {Be1} in 1983-1984
  proved $^{\#}\td{\mathcal{O}}(\Om)\geq 1$ if $V$ is
$C^1$, $\bar{\Om}=\{V\leq h\}$ is compact, and $V'(q)\neq 0$ for all
$q\in \partial{\Om}$. In 1987, P. Rabinowitz in \cite{Ra3} proved
that if $H$ is $C^1$ and satisfies the reversible conditon,
$\Sg\equiv H^{-1}(h)$ is star-shaped, and $x\cdot H'(x)\neq 0$ for
all $x\in \Sg$, then $^{\#}\td{\mathcal{J}}_b(\Sg)\geq 1$. In 1987,
V. Benci and F. Giannoni gave a different proof of the existence of
one brake orbit in \cite{BG}.

In 1989, A. Szulkin in \cite{Sz} proved that
$^{\#}\td{\J_b}(H^{-1}(h))\geq n$, if $H$ satisfies conditions in
\cite{Ra1} of Rabinowitz and the energy hypersurface $H^{-1}(h)$ is
$\sqrt{2}$-pinched. E. van Groesen in \cite{Gro} of 1985 and A.
Ambrosetti, V. Benci, Y. Long in \cite{ABL1} of 1993 also proved
$^{\#}\td{\mathcal{O}}(\Om)\geq n$ under different pinching
conditions.

 In \cite{LZZ} of 2006, Long , Zhu and the author of
this paper proved that there exist at least $2$ geometrically
distinct brake orbits on any central symmetric strictly convex
hypersuface $\Sg$ in $\R^{2n}$ for $n\ge 2$. Recently, in \cite{LZ},
Liu and the author of this paper proved that there exist at least
$[n/2]+1$ geometrically distinct brake orbits on any central
symmetric strictly convex hypersuface  $\Sg$ in $\R^{2n}$ for $n\ge
2$,
 if all brake orbits on $\Sg$ are nondegenerate then there are at
 least $n$ geometrically distinct brake orbits on $\Sg$. For more
 details one can refer to \cite{LZZ}, \cite{LZ} and the reference
 there in.

In his pioneering paper \cite{Ra1} of 1978, P. Rabinowitz proved the
following famous result via the variational method. Suppose $H$
satisfies the following conditions:

($\rm {H}1'$) $H\in C^1(\R^{2n},\R)$.

(H2) There exist constants $\mu >2$ and $r_0>0$ such that
       $$0<\mu H(x)\le H'(x)\cdot x,\quad \forall |x|\le r_0.$$

(H3) $H(x)=o(|x|^2)$ at $x=0$.

(H4) $H(x)\ge 0$ for all $x\in \R^{2n}$.

{\it Then for any $T>0$, the system (\ref{1.1}) possesses a
non-constant $T$-periodic solution.}  Because a $T/k$ periodic
function is also a $T$-periodic function, in \cite{Ra1} Rabinowitz
proposed a conjecture that {\it under conditions (H1$'$) and
(H2)-(H4), there is a non-constant solution possessing any
prescribed minimal period.} Since 1978, this conjecture has been
deeply studied by many mathematicians. A significant progress was
made by Ekeland and Hofer in their celebrated paper \cite{EH1} of
1985, where they proved Rabinowitz's conjecture for the strictly
convex Hamiltonian system. For Hamiltonian systems with convex or
weak convex assumptions, we refer to \cite{AmC1}-\cite{AmM1},
\cite{CZ1}-\cite{DL1}, \cite{Ek1}-\cite{FQ1}, \cite{Long5},
\cite{GirMat0}-\cite{GirMat3}, and references therein for more
details. For the case without convex condition we refer to
\cite{Long1}-\cite{Long3} and Chapter 13 of \cite{Long5} and
references therein. A interesting result is for the semipositive
first order Hamiltonian system, in \cite{FKW1} G. Fei, S.-T. Kim,
and T. Wang proved the existence of a T periodic solution of system
(\ref{1.1}) with minimal period no less than $T/2n$ for any given
$T>0$.

Note that in the second order Hamiltonian systems there are many
results on the minimal problem of brake orbits such  us \cite
{Long1}-\cite{Long3} and \cite{Xiao}. For the even first order
Hamiltonian system, in \cite{Z1}, the author of this paper studied
the minimal period problem of semipositive even Hamiltonian system
and gave a positive answer to Rabinowitz's conjecture in that case.
In \cite{FKW2}, G. Fei, S.-T. Kim, and T. Wang proved the same
result for second order Hamiltonian systems.

So it is natural to consider the minimal period problem of brake
orbits in reversible first order nonlinear Hamiltonian systems. In
\cite{Liu3}, Liu have considered the strictly convex reversible
Hamiltonian systems case and proved the existence of nonconstant
brake orbit of (\ref{1.1}) with minimal period belonging to
$\{T,T/2\}$ for any given $T>0$.

Since \cite{Z1}, we also hope to obtain some interesting results in
the even Hamiltonian system for the minimal period problem of brake
orbits.

It can be found in many papers mentioned above that the Maslov-type
index theory and its iteration theory play a important role in the
study of minimal period problems in Hamiltonian systems. In this
paper we study some monotonicity properties of Maslov-type index and
apply it to prove our main results.

In this paper we denote by $\mathcal{L}(\R^{2n})$ and
$\mathcal{L}_s(\R^{2n})$ the set of all real $2n\times 2n$ matrices
and symmetric matrices respectively. And we denote by $y_1\cdot y_2$
the usual inner product for all $y_1,\;y_2\in\R^k$ with $k$ being
any positive integer. Also we denote by $\N$ and $\Z$ the set of
positive integers and integers respectively.

Let $\Sp(2n)=\{M\in\mathcal{L}(\R^{2n})|M^TJM=J\}$ be the $2n\times
2n$ real symplectic group. For any $\tau>0$, Set
$\mathcal{P}_\tau=\{\ga\in C([0,\tau],\Sp(2n))|\ga(0)=I_{2n}\}$ and
$S_\tau=\R/(\tau\Z)$.

For any $\ga\in \mathcal{P}_\tau$ and $\om\in \U$, where $\U$ is the
unit circle of the complex plane $\C$, the Maslov-type index
$(i_\om(\ga),\nu_\om(\ga))\in \Z\times\{0,1,...2n\}$ was defined by
Long in \cite{Long4}. We have a brief review in Appendix of Section
6.

For convenience to introduce our results, we define the following
({\bf B1}) condition, since the Hamiltonian systems considered here
are reversible, this condition is natural.

 {\bf (B1) Condition}. For any $\tau>0$ and
$B\in C([0,\tau],\mathcal{L}_s(\R^{2n})$ with the $n\times n$ matrix
square block form  $B(t)=\left(
\begin{array}{cc}B_{11}(t)&B_{12}(t)\\B_{21}(t)&B_{22}(t)\end{array}\right)$ satisfying
$B_{12}(0)=B_{21}(0)=0=B_{12}(\tau)=B_{21}(\tau)$, We will call $B$
satisfies the condition {\bf (B1)}.

Throughout this paper, we denote by
 \bea L_0=\{0\}\times\R^n\subset
\R^{2n},\quad L_1=\R^n\times\{0\}\subset\R^{2n}.\lb{s1}\eea
 The definitions of Maslov-type indices $(i_{\sqrt{-1}}^{L_0}(\ga),\nu_{\sqrt{-1}}^{L_0}(\ga))$
 and  $(i_{L_j}(\ga),\nu_{L_j}(\ga))\in \Z\times
 \{0,1,...,n\}$ for $j=0,1$ and $\ga\in \mathcal{P}_\tau(2n)$ with
 $\tau>0$ can be found in \cite{LZZ} and Section 2 below. Also for
$B\in C([0,\tau],\mathcal{L}_s(\R^{2n})$ satisfies condition (B1),
the definitions of
$(i_{\sqrt{-1}}^{L_0}(B),\nu_{\sqrt{-1}}^{L_0}(B))$
 and  $(i_{L_j}(B),\nu_{L_j}(B))\in \Z\times
 \{0,1,...,n\}$ for $j=0,1$ and $\ga\in \mathcal{P}_\tau(2n)$
 can be found in Section 2 and references therein.

 For any $B\in C([0,\tau],\mathcal{L}_s(\R^{2n}))$, denote by
 $\ga_B$ the fundamental solution of the following problem:
        \bea \dot{\ga}_B(t)&=&JB(t)\ga_B(t),\\
             \ga_B(0)&=&I_{2n}.\lb{1.3}\eea
 Then $\ga_B\in \mathcal{P}_\tau$. We call $\ga_B$ the {\it symplectic
 path associated} to $B$.

\noindent{\bf Definition 1.1.} If $H\in C^2(\R^{2n},\R)$ is a
reversible function, for any $x_\tau$ be a $\tau$-periodic brake
orbit solution of (\ref{1.1}), let $B(t)=H''(x(t))$, we define
$\ga_{x_\tau}=\ga_B|_{[0,\frac{\tau}{2}]}$ and call it the
symplectic path associated to $x_\tau$. We define  \be
i_{L_0}(x_\tau)=i_{L_0}(\ga_{x_\tau}),\qquad
\nu_{L_0}(x_\tau)=i_{L_0}(\ga_{x_\tau}).\ee Moreover, if $H$ is even
and $x_\tau$ is a $\tau$-periodic  symmetric brake orbit solution of
(\ref{1.1}), let $B(t)=H''(x(t))$, we define
$\ga_{x_\tau}=\ga_B|_{[0,\frac{\tau}{4}]}$ and call it the
symplectic path associated to $x_\tau$. We define \be
i_{\sqrt{-1}}^{L_0}(x_\tau)=i_{\sqrt{-1}}^{L_0}(\ga_{x_\tau}),\qquad
\nu_{\sqrt{-1}}^{L_0}(x_\tau)=i_{\sqrt{-1}}^{L_0}(\ga_{x_\tau}).\ee

\noindent{\bf Definition 1.2.} For any $\tau$-period and $k\in
\N\equiv\{1,2,...\}$, we define the $k$ times iteration $x^k$ of $x$
by
       \be x^k(t)=x(t-j\tau),\quad j\tau\le t\le (j+1)\tau,\quad
       0\le j\le k.\ee
As in \cite{LZ}, for any $\ga\in \mathcal{P}_\tau$ and $k\in
\N\equiv\{1,2,...\}$, in this paper the $k$-time iteration $\ga^k$
of $\ga\in \mathcal{P}_\tau(2n)$ in brake orbit boundary sense is
defined by $\td{\ga}|_{[0,k\tau]}$ with \bea
   \td{\ga}(t)=\left\{\begin{array}{l}\ga(t-2j\tau)(N\ga(\tau)^{-1}N\ga(\tau))^j,
   \; t\in[2j\tau,(2j+1)\tau], j=0,1,2,...\\
      N\ga(2j\tau+2\tau-t)N(N\ga(\tau)^{-1}N\ga(\tau))^{j+1}\; t\in[(2j+1)\tau,(2j+2)\tau],
      j=0,1,2,...\end{array}\right.\nn\eea

The followings are our main results of this paper.

\noindent{\bf Theorem 1.1.} {\it Suppose that $H$ satisfies
conditions (H1)-(H4) and

(H5) $H''(x)$ is semipositive definite for all $x\in\R^{2n}$.

Then for any $T>0$, the system (1.1)-(1.2) possesses a nonconstant
$T$ periodic brake orbit solution $x_T$ with minimal period no less
that $\frac{T}{2n+2}$. Moreover, for $x=(x_1,x_2)$ with
$x_1,x_2\in\R^n$, denote by $H''_{22}(x)$ the second order
differential of $H$ with respect to $x_2$, if
      \be \int_0^{\frac{T}{2}}H''_{22}(x_T(t))\,dt>0,\lb{zhi6}\ee
then the minimal period of $x_T$ belongs to $\{T,\frac{T}{2}\}$.
       }

\noindent{\bf Remark 1.1.} (Theorem 1.1 of \cite{Liu3}) {\it Suppose
that $H$ satisfies conditions (H1)-(H4) and if $x_T$ satisfies

(H5$'$) $\int_0^{\frac{T}{2}}H''(X_T(t))\,dt>0$.

Then the minimal period of $x_T$ belongs to $\{T,\frac{T}{2}\}$.}

 In the case $n=1$, the result can be better, i.e., the following

\noindent{\bf Theorem 1.2.} {\it For $n=1$, suppose that $H$
satisfies conditions (H1)-(H4).

Then for any $T>0$, the system (1.1)-(1.2) possesses a nonconstant
$T$ periodic brake orbit solution with minimal period belong to
$\{T,\frac{T}{2}\}$.}

Consider the minimal period problem for $H(x)=\frac{1}{2}B_0x\cdot
x+\hat{H}(x)$, where $B_0\in \mathcal{L}_s(\R^{2n})$. This is
motivated by \cite{FKW1}, \cite{GirMat1}, and \cite{Ra1}, where in
\cite{FKW1} $B_0$ was considered to be semipositive, in
\cite{GirMat1} and \cite{Ra1} $B_0$ was considered to be positive.

We have the following general result.

\noindent{\bf Theorem 1.3.} {\it Let $2n\times 2n$ be real
semipositive matrix $B_0=\diag(B_{11},B_{22})$ with $B_{11}$ and
$B_{22}$ being $n\times n$ matrix. Assume $H(x)=\frac{1}{2}B_0x\cdot
x+\hat{H}(x)$ for all $x\in \R^{2n}$, and $\hat{H}$ satisfies
conditions (H1)-(H5).

Then for any $T>0$, (1.1) possesses a nonconstant $T$-periodic brake
orbit $x_T$ with minimal period no less than
$\frac{T}{2i_{L_0}(B_0)+2\nu_{L_0}(B_0)+2n+2}$, where we see $B_0$
as an element in $C([0,T/2], \mathcal{L}_s(\R^{2n}))$ satisfies
condition (B1). }

\noindent{\bf Remark 1.2.} In section 3, we will show
$i_{L_0}(B_0)+\nu_{L_0}(B_0)\ge 0$.

As a direct consequence of Theorem 1.3, we have the following
Corollary 1.1.

\noindent {\bf Corollary 1.1.} {\it For $T>0$ such that
$i_{L_0}(B_0)+\nu_{L_0}(B_0)=0$, where we see $B_0$ as an element in
$C([0,T/2], \mathcal{L}_s(\R^{2n})$ satisfies condition (B1), under
the same assumptions of Theorem 1.2, the system (1.1) possesses a
nonconstant $T$-periodic brake orbit with minimal period no less
that $\frac{T}{2n+2}$.}

We can also prove the following Corollary 1.2 of Theorem 1.3.

\noindent {\bf Corollary 1.2.} {\it If $B_0\neq 0$, then for
$0<T<\frac{\pi}{||B_0||}$ with  $||B_0||$ being the operator norm of
$B_0$, under the same condition of Theorem 1.2,  possesses a
nonconstant $T$-periodic brake orbit $x_T$ with minimal period no
less than $\frac{T}{2n+2}$. Moreover , if
      \bea \int_0^{\frac{T}{2}}H''_{22}(x_T(t))\,dt>0, \nn\eea
then the minimal period of $x_T$ belongs to $\{T,\frac{T}{2}\}$.
       }\bigskip

\noindent{\bf Theorem 1.4.} {\it Suppose that $H$ satisfies
conditions (H1)-(H5) and

(H6) $H(-x)=H(x)$ for all $x\in\R^{2n}$.

Then for any $T>0$, the system (1.1)-(1.2) possesses a nonconstant
symmetric brake orbit with minimal period belonging to $\{T,T/3\}$.}

\noindent{\bf Theorem 1.5.} {\it Let $2n\times 2n$ be real
semipositive matrix $B_0=\diag(B_{11},B_{22})$ with $B_{11}$ and
$B_{22}$ being $n\times n$ matrix, assume $H(x)=\frac{1}{2}B_0x\cdot
x+\hat{H}(x)$ for all $x\in \R^{2n}$, and $\hat{H}$ satisfies
conditions (H1)-(H6). Then for any $T>0$, the system (1.1)-(1.2)
possesses a nonconstant symmetric brake orbit $x_T$ with minimal
period no less than
$\frac{T}{4(i_{\sqrt{-1}}^{L_0}(B_0)+\nu_{\sqrt{-1}}^{L_0}(B_0))+7}$.
Moreover, if $i_{\sqrt{-1}}^{L_0}(B_0)+\nu_{\sqrt{-1}}^{L_0}(B_0)$
is even, then the minimal period of $x_T$ is no less than
$\frac{T}{4(i_{\sqrt{-1}}^{L_0}(B_0)+\nu_{\sqrt{-1}}^{L_0}(B_0))+3}$£¬where
we see $B_0$ as an element in $C([0,T/4], \mathcal{L}_s(\R^{2n})$
satisfies condition (B1).}

\noindent{\bf Remark 1.3.} In section 3, we will show that
$i_{\sqrt{-1}}^{L_0}(B_0)\ge 0$, hence
$i_{\sqrt{-1}}^{L_0}(B_0)+\nu_{\sqrt{-1}}^{L_0}(B_0)\ge 0$.

As a direct consequence of Theorem 1.5, we have the following
Corollary 1.3.

 \noindent {\bf Corollary 1.3.} {\it For $T>0$ such
that $i_{\sqrt{-1}}^{L_0}(B_0)+\nu_{\sqrt{-1}}^{L_0}(B_0)=0$, under
the same assumptions of Theorem 1.4, the system (1.1) possesses a
nonconstant symmetric brake orbit with minimal period belonging to
$\{T,T/3\}$.}

We can also prove the following Corollary 1.4 of Theorem 1.5.

\noindent {\bf Corollary 1.4.} {\it If $B_0\neq 0$, then for
$0<T<\frac{\pi}{||B_0||}$ with $||B_0||$ being the operator norm of
$B_0$, under the same condition of Theorem 1.5, the system (1.1)
possesses a nonconstant symmetric brake orbit with minimal period
belonging to $\{T,T/3\}$.}

This paper is organized as follows. In section 2, we study the
Maslov-type index theory of $i_{L_0}$, $i_{L_1}$ and
$i_{\sqrt{-1}}^{L_0}$. We compute the difference between
$i_{L_0}(\ga)$ and $i_{L_1}(\ga)$. In Section 3, we study the
relation between the Maslov-type index
$(i_{\sqrt{-1}}^{L_0}(B),\nu_{\sqrt{-1}}^{L_0}(B))$ for $B\in
C([0,\tau],\mathcal{L}_s(\R^{2n})$ satisfies condition (B1) and the
Morse indices of the corresponding Galerkin approximation. As
applications we get some monotonicity properties of $i_{L_0}(B)$,
$i_{L_1}(B)$ and $i_{\sqrt{-1}}^{L_0}(B)$ and we prove Theorem 3.2
which is very important in the proof of Theorems 1.4-1.5. In Section
4, based on the preparations in Sections 2 and 3 we prove Theorems
1.1-1.3 and Corollary 1.2. In Section 5, we prove Theorems 1.4-1.5
and corollary 1.4. In Section 6, we give a briefly review of
$(i_\om,\nu_\om)$ index theory  with $\om\in \U$ for symplectic
paths starting with identity as appendix.

\setcounter{equation}{0}
\section{Maslov-type index theory associated with Lagrangian subspaces}
\subsection {A brief review of index function $(i_{L_j},\nu_{L_j})$ with $j=0,1$ and
$(i_{\sqrt{-1}}^{L_0},\nu_{\sqrt{-1}}^{L_0})$} 

Let \be F=\R^{2n}\oplus \R^{2n}\lb{zhang0}\ee possess the standard
inner product. We define the symplectic structure of $F$ by
      \be \{v,w\}=(\mathcal{J}v,w),\;\forall v,w\in F,\;
 {\rm where}\; \mathcal{J}=(-J)\oplus J=\left(\begin{array}{cc} -J &0\\0&J\end{array}\right).
 \;\lb{zhang2}\ee
  We denote by $\Lag(F)$ the set of Lagrangian subspaces of $F$,
  and equip it with the topology as a subspace of the Grassmannian of all
  $ 2n$-dimensional subspaces of $F$.

 It is easy to check that, for any $M\in \Sp(2n)$ its
  graph
      $$\Gr(M)\equiv\left\{\left(\begin{array}{c}x\\Mx\end{array}\right)|x\in
      \R^{2n}\right\}$$
is a Lagrangian subspace of $F$.

Let \bea V_1=\{0\}\times \R^n\times \{0\}\times \R^n\subset \R^{4n},
\quad V_2=\R^n\times \{0\}\times \R^n\times\{0\}\subset
\R^{4n}.\lb{s2}\eea

By Proposition 6.1 of \cite{LZ} and Lemma 2.8 and Definition 2.5 of
\cite{LZZ}, we give the following definition.

 \noindent{\bf Definition 2.1.} For any continuous path $\ga\in\mathcal{P}_\tau(2n)$, we
 define the following Maslov-type indices:

 \bea && i_{L_0}(\ga)=\mu^{CLM}_{F}(V_1, \Gr(\ga),[0,\tau])-n,\\
     && i_{L_1}(\ga)=\mu^{CLM}_{F}(V_2, \Gr(\ga),[0,\tau])-n,\\
      && \nu_{L_j}(\ga)=\dim (\ga(\tau)L_j\cap L_j),\qquad j=0,1,\eea
where we denote by $i^{CLM}_F(V,W,[a,b])$ the Maslov index for
Lagrangian subspace path pair $(V,W)$ in $F$ on $[a,b]$ defined by
Cappell, Lee, and Miller in \cite{CLM}.

For $\omega=e^{\sqrt{-1}\theta}$ with $\theta\in\R$, we define a
Hilbert space $E^{\omega}=E^{\omega}_{L_0}$ consisting of those
$x(t)$ in $L^2([0,\tau], \C^{2n})$ such that $e^{-\theta t J}x(t)$
has Fourier expending
$$e^{-\frac{\theta t}{\tau} J}x(t)=\sum_{j\in \Z}e^{\frac{j\pi t}{\tau}J}\pmatrix 0\\a_j\endpmatrix,\;a_j\in \C^n  $$
with
$$\|x\|^2:=\sum_{j\in\Z}\tau(1+|j|)|a_j|^2<\infty. $$
 For
$\omega=e^{\sqrt{-1}\theta}$, $\theta\in (0,\pi)$, we define two
self-adjoint operators $A^{\omega}, B^{\omega}\in \mathcal
{L}(E^{\omega})$ by \bea  (A^{\omega}x,y)=\int^1_0\<-J\dot
x(t),y(t)\>dt,\;\;(B^{\omega}x,y)=\int^1_0\<B(t)x(t),y(t)\>dt
\nn\eea
 on $E^{\omega}$. Then $B^{\omega}$ is also compact.

\noindent{\bf Definition 2.2.} We define the index function \bea
&&i_{\omega}^{L_0}(B)=I(A^{\omega},
\;\;A^{\omega}-B^{\omega})\equiv-{\rm sf}\{A^{\omega}
 -sB^{\omega},0\le s\le 1\},\nn\\
 &&\nu_{\omega}^{L_0}(B)=m^0(A^{\omega}-B^{\omega}),\;
  \forall\,\omega=e^{\sqrt{-1}\theta},\;\;\theta\in (0,\pi),\nn\eea
where the definition of $\sf$ of spectral flow for the path of
bounded self-adjoint linear operators  one can refer to \cite{ZL}
and references their in.

By (3.21) of \cite{LZ}, we have
 \be i_{L_0}(B)\le i^{L_0}_{\om}(B)\le
i_{L_0}(B)+n.\ee

\noindent{\bf Lemma 2.1.} {\it For $\om=e^{\sqrt{-1}\theta}$ with
$\theta\in(0,\pi)$, let $V_\om=L_0\times (e^{\theta J}L_0)\subset
\R^{4n}\equiv F$. There holds
    \be
    i_{\omega}^{L_0}(B)=\mu_F^{CLM}(V_\om,\Gr(\ga_B),[0,\tau]).\ee}
\quad{\bf Proof.} Without loss of generality we can suppose the
$C^1$ path $\Gr(\ga_B)$ of Lagrangian subspaces intersects $V_\om$
regularly (otherwise we can perturb it slightly with fixed endow
points such that they intersects regularly and the index dose not
change by the homotopy invariant property $\mu_F^{CLM}$ ), where the
definition of intersection form can be found in \cite{RS1}. We
denote by $\mu^{BF}$ the maslov index defined by Booss and Furutani
in \cite{BF1}.

 By the spectral flow formula of Theorem 5.1
in \cite{BF1} or Theorem 1.5 of \cite{BZ1} (cf. also proof of
Proposition 2.3 of \cite{Z2}), we have
   \bea &&{\rm sf}\{A^{\omega}
 sB^{\omega},0\le s\le 1\}\nn\\&=&\mu^{BF}(\Gr(\ga_B),V_\om,
 [0,\tau])\nn\\
 &=&\mu^{BF}((I\oplus e^{-\sqrt{-1}\theta J})\Gr(\ga_B),(I\oplus e^{-\sqrt{-1}\theta J})V_\om,
 [0,\tau])\nn\\
 &=&\mu^{BF}((I\oplus e^{-\sqrt{-1}\theta J})\Gr(\ga_B),V_1,
 [0,\tau])\nn\\
&=& -m^-(-\Ga((I\oplus
 e^{-\sqrt{-1}\theta J})\Gr(\ga_B),V_1,0))+
 \sum_{0<t<\tau}\sign(-\Ga((I\oplus
 e^{-\sqrt{-1}\theta J})\Gr(\ga_B),V_1,t))\nn\\&&+
 m^+(-\Ga((I\oplus
 e^{-\sqrt{-1}\theta J})\Gr(\ga_B),V_1,\tau))\nn\\
&=& -m^+(\Ga((I\oplus
 e^{-\sqrt{-1}\theta J})\Gr(\ga_B),V_1,0))-
 \sum_{0<t<\tau}\sign(\Ga((I\oplus
 e^{-\sqrt{-1}\theta J})\Gr(\ga_B),V_1,t))\nn\\&&+
 m^-(\Ga((I\oplus
 e^{-\sqrt{-1}\theta J})\Gr(\ga_B),V_1,\tau))\nn\\
&=& -\mu_F^{CLM}(V_1, (I\oplus
 e^{-\sqrt{-1}\theta J})\Gr(\ga_B), [0,\tau])\nn\\
 &=& -\mu_F^{CLM}((I\oplus
 e^{\sqrt{-1}\theta J})V_1, \Gr(\ga_B), [0,\tau])\nn\\
 &=& -\mu_F^{CLM}(V_\om,\Gr(\ga_B),[0,\tau]),
 \eea
where in the fourth equality we have used Theorem 2.1 in \cite{BF1}
and the property of index $\mu^{RS}$ for symplectic paths defined in
\cite{RS1}(cf also (2.6)-(2.8) of \cite{Z2}), in the sixth equality
we have used Lemma 2.6 of \cite{LZZ}, in the second and seventh
equalities we used the symplectic invariance property of index
$\mu^{BF}$ and $\mu_F^{CLM}$ respectively. \hfill\hb

\noindent{\bf Definition 2.3.} Let $B\in
C([0,\tau],\mathcal{L}_s(\R^{2n})$ and $\ga_B$ be the symplectic
path associated to $B$. We define
 \bea i_{\om}^{L_0}(\ga_B)=i_{\om}^{L_0}(B),\\
      \nu_{\om}^{L_0}(\ga_B)=\nu_{\om}^{L_0}(B). \eea

By Lemma 2.1, in general we give the following definition.

\noindent{\bf Definition 2.4.} For any $\ga\in \mathcal{P}_\tau(2n)$
and $\om=e^{\sqrt{-1}\theta}$ with $\theta\in (0,\pi)$, we define
 \bea
 i_\om ^{L_0}(\ga)=\mu_F^{CLM}(V_\om,\Gr(\ga_B),[0,\tau]),\nn\\
   \nu_\om^{L_0}(\ga)=\dim \left(\ga(\tau)L_0\cap e^{\sqrt{-1}\theta J}L_0\right).\eea

 For any $\ga\in\mathcal{P}_\tau(2n)$, we define a new
symplectic path $\td{\ga}\in \mathcal{P}_{\tau}(2n)$ by
   \be \td{\ga}(t)=\left\{\begin{array}{l}I_{2n},\quad t\in
   [0,\frac{\tau}{3}],\\ \ga(3t-\tau),\quad t\in[\frac{\tau}{3},\frac{2\tau}{3}],\\
   \ga(\tau),\quad t\in[\frac{2\tau}{3},\tau].\end{array}\right.\ee

So we can perturb $\td{\ga}$ slightly to a $C^1$ path $\hat{\ga}$
such that $\hat{\ga}$ is homotopic to $\td{\ga}$ with fixed end
points and $\hat{\ga}(t)=I_{2n}$ for $t\in [0,\frac{\tau}{6}]$ and
$\hat{\ga}(t)=\ga(\tau)$ for $t\in [\frac{5\tau}{6},\tau]$. Set
$\hat{B}(t)=-J\dot{\hat{\ga}}(t)(\hat{\ga}(t))^{-1}$. So we have \be
\hat{B}(0)=\hat{B}(\tau)=0.\ee Then this $\hat{B}\in
C([0,\tau],\mathcal{L}_s(\R^{2n})$ and satisfies condition (B1).
Also we have $\hat{\ga}$ is is homotopic to $\ga$ with fixed end
points. So we have
    \bea
   && i_1(\hat{\ga}^k)=i_1(\ga^k)=i_1(\ga_{\hat{B}}^k),\qquad \forall
   k\in\N,\\
   && \nu_1(\hat{\ga}^k)=\nu_1(\ga^k)=\nu_1(\ga_{\hat{B}}^k),\qquad \forall k\in\N\eea
and \bea
   && i_{L_0}(\hat{\ga}^k)=i_{L_0}(\ga^k)=i_{L_0}(\ga_{\hat{B}}^k),\qquad \forall
   k\in\N,\\
   && \nu_{L_0}(\hat{\ga}^k)=\nu_{L_0}(\ga^k)=\nu_{L_0}(\ga_{\hat{B}}^k),\qquad \forall k\in\N.\eea
Also by the property of index $\mu_F^{CLM}$ and Definition 2.4 have
 \bea
 &&i_{\sqrt{-1}}^{L_0}(\ga^k)=i_{\sqrt{-1}}^{L_0}(\hat{\ga}^k)=i_{\sqrt{-1}}^{L_0}(\ga_{\hat{B}}^k),
 \qquad \forall k\in\N,\nn\\
 &&\nu_{\sqrt{-1}}^{L_0}(\ga^k)=\nu_{\sqrt{-1}}^{L_0}(\hat{\ga}^k)=\nu_{\sqrt{-1}}^{L_0}(\ga_{\hat{B}}^k),
 \qquad \forall k\in\N.\nn\eea

 Hence, in \cite{LZ} the authors essentially proved the following Bott-type iteration
formula.

\noindent{\bf Theorem 2.1.} (Theorem 4.1 of \cite{LZ}) {\it Let
$\ga\in \mathcal{P}_\tau(2n)$ and $\omega_k=e^{\pi \sqrt{-1}/k}$.
For odd $k$ we have
\bea i_{L_0}(\gamma^k)=i_{L_0}(\gamma^1)+\sum_{i=1}^{(k-1)/2}i_{\omega_k^{2i}}(\gamma^2),\nn\\
\nu_{L_0}(\gamma^k)=\nu_{L_0}(\gamma^1)+\sum_{i=1}^{(k-1)/2}\nu_{\omega_k^{2i}}(\gamma^2),\nn\eea
and for even $k$, we have \bea &&
i_{L_0}(\gamma^k)=i_{L_0}(\gamma^1)+i^{L_0}_{\sqrt{-1}}(\gamma^1)+
\sum_{i=1}^{k/2-1}i_{\omega_k^{2i}}(\gamma^2),\;
\nn\\
&&
\nu_{L_0}(\gamma^k)=\nu_{L_0}(\gamma^1)+\nu^{L_0}_{\sqrt{-1}}(\gamma^1)+
\sum_{i=1}^{k/2-1}\nu_{\omega_k^{2i}}(\gamma^2). \nn\eea}

Obviously we also have
 \be i_{L_0}(\ga)\le i^{L_0}_{\sqrt{-1}}(\ga)\le
i_{L_0}(\ga)+n.\ee

\subsection{The Bott-type iteration formula for
$(i_{\sqrt{-1}}^{L_0},\nu_{\sqrt{-1}}^{L_0})$} 

In order to study the minimal period problem for Even reversible
Hamiltonian systems, we need the iteration formula of the
Maslov-type index of $(i_{\sqrt{-1}}^{L_0},\nu_{\sqrt{-1}}^{L_0})$
for symplectic paths starting with identity. We use Theorem 2.1 to
obtain it.

Precisely we have the following Theorem.

\noindent{\bf Theorem 2.2.} {\it Let $\ga\in \mathcal{P}_\tau(2n)$
and $\omega_k=e^{\pi \sqrt{-1}/k}$. For odd $k$ we have
\bea i_{\sqrt{-1}}^{L_0}(\gamma^k)=i_{\sqrt{-1}}^{L_0}(\gamma^1)+\sum_{i=1}^{(k-1)/2}i_{\omega_k^{2i-1}}(\gamma^2),\lb{c7}\\
\nu_{\sqrt{-1}}^{L_0}(\gamma^k)=\nu_{\sqrt{-1}}^{L_0}(\gamma^1)+\sum_{i=1}^{(k-1)/2}\nu_{\omega_k^{2i-1}}(\gamma^2),\lb{c8}\eea
and for even $k$, we have \bea &&
 i_{\sqrt{-1}}^{L_0}(\gamma^k)=
\sum_{i=1}^{k/2}i_{\omega_k^{2i-1}}(\gamma^2),\;
\lb{c9}\\
&& \nu_{\sqrt{-1}}^{L_0}(\gamma^k)=
\sum_{i=1}^{k/2}\nu_{\omega_k^{2i-1}}(\gamma^2). \lb{c10}\eea}

\noindent{\bf Proof.} For odd $k$, since $\ga^{2k}=(\ga^k)^2$, by
Theorem 2.1 we have
 \bea i_{L_0}(\ga^{2k})=i_{L_0}(\ga^k)+i_{\sqrt{-1}}^{L_0}(\ga^k),\lb{c1}\\
      \nu_{L_0}(\ga^{2k})=\nu_{L_0}(\ga^k)+\nu_{\sqrt{-1}}^{L_0}(\ga^k).\lb{c2}\eea

Also by Theorem 2.1 we have
 \bea i_{L_0}(\gamma^k)=i_{L_0}(\gamma^1)+\sum_{i=1}^{(k-1)/2}i_{\omega_k^{2i}}(\gamma^2),\lb{c3}\\
\nu_{L_0}(\gamma^k)=\nu_{L_0}(\gamma^1)+\sum_{i=1}^{(k-1)/2}\nu_{\omega_k^{2i}}(\gamma^2),\lb{c4}\\
i_{L_0}(\gamma^{2k})=i_{L_0}(\gamma^1)+i_{\sqrt{-1}}^{L_0}(\ga)+\sum_{i=1}^{k-1}i_{\omega_{2k}^{2i}}(\gamma^2),\lb{c5}\\
\nu_{L_0}(\gamma^{2k})=\nu_{L_0}(\gamma^1)+\nu_{\sqrt{-1}}^{L_0}(\ga)+\sum_{i=1}^{k-1}\nu_{\omega_{2k}^{2i}}(\gamma^2).\lb{c6}\eea

Since $\om_k=\om_{2k}^2$, by (\ref{c1}), (\ref{c5}) minus (\ref{c3})
yields (\ref{c7}). By (\ref{c2}), (\ref{c6}) minus (\ref{c4}) yields
(\ref{c8}).

For even k, by similar argument we obtain (\ref{c9}) and
(\ref{c10}). The proof of Theorem 2.2 is complete.\hfill\hb

\subsection {The difference of $i_{L_0}(\ga)$ and $i_{L_1}(\ga)$.}

The precise difference of $i_{L_0}(\ga)$ and $i_{L_1}(\ga)$ for
$\ga\in\mathcal{P}_\tau$ with $\tau>0$ is very important in the
proof of the main results of this paper. In this subsection we use
the H$\ddot{{\rm o}}$rmander index (cf. \cite{Dui}) to compute it.
Note that in \cite{LZZ}, in fact we have already proved that
$|i_{L_0}(\ga)-i_{L_1}(\ga)|\le n$.

For any $P\in\Sp(2n)$ and $\varepsilon\in\R$, we set
 \bea M_\varepsilon(P)=P^T\left(\begin{array}{cc}\sin{2\var}I_n&-\cos{2\var I_n}\\-\cos{2\var}I_n&-\sin
 2\var I_n
 \end{array}\right)P+\left(\begin{array}{cc}\sin{2\var}I_n&\cos{2\var}I_n\\\cos{2\var}I_n&-\sin2\var
 I_n
 \end{array}\right).\eea
Then we have the following theorem.

\noindent{\bf Theorem 2.3.} {\it For $\ga\in\mathcal{P}_\tau$ with
$\tau>0$, we have
   \be i_{L_0}(\ga)-i_{L_1}(\ga)=\frac{1}{2}\sgn
M_\var(\ga(\tau)),\lb{p20}\ee where $\sgn M_\var(\ga(\tau))$ is the
signature of the symmetric matrix $M_\var(\ga(\tau))$ and $\var>0$
is sufficiently small.

we also have,
 \be
 (i_{L_0}(\ga)+\nu_{L_0}(\ga))-(i_{L_1}(\ga)+\nu_{L_1}(\ga))=\frac{1}{2}\sign
M_\var(\ga(\tau)),\lb{p21}\ee where $\var<0$ and $|\var|$ is
sufficiently small.}

{\bf Proof.} By the first geometrical definition of the Maslov-type
index in Section 4 of \cite{CLM}, there exists an $\var>0$ small
enough such that
  \be V_1\cap e^{-\var\mathcal{J}}\Gr(\ga(0))=\{0\},\qquad  V_2\cap
  e^{-\var\mathcal{J}}\Gr(\ga(\tau))=\{0\}.\ee
By definition 2.1, we have
 \bea && i_{L_0}(\ga)=\mu^{CLM}_{F}(V_1, e^{-\var\mathcal{J}}\Gr(\ga),[0,\tau])-n,\lb{p1}\\
     && i_{L_1}(\ga)=\mu^{CLM}_{F}(V_2, e^{-\var\mathcal{J}}\Gr(\ga),[0,\tau])-n.\lb{p2}\eea
 Define $\ga_1(t)=e^{-\var\mathcal{J}}\Gr(\ga(t))$ and $\ga_2(t)=e^{-\var\mathcal{J}}\Gr(\ga(\tau-t))$
 for $t\in [0,\tau]$. Then $\ga_1$ and $\ga_2$ are two paths of Lagrangian subspaces
 of  the symplectic space $(F,\mathcal{J})$ defined in (\ref{zhang0})
 and (\ref{zhang2}). $\ga_1$ connects $e^{-\var\mathcal{J}}\Gr(\ga(0))$
 and $e^{-\var\mathcal{J}}\Gr(\ga(\tau))$ and is transversal to
 $V_1$ and $V_2$. $\ga_2$ connects $e^{-\var\mathcal{J}}\Gr(\ga(\tau))$
 and $e^{-\var\mathcal{J}}\Gr(\ga(0))$ and is transversal to
 $V_1$ and $V_2$. Denote by $\ga$ the catenation of the paths
 $\ga_1$ and $\ga_2$.
       By Definition 3.4.2 of the
$H\ddot{o}rmande\; index$ $s(M_1,M_2;L_1,L_2)$ on p. 66 of
\cite{Dui} and (\ref{p1})-(\ref{p2}), we have
      \bea
      &&s(V_1,V_2;e^{-\var\mathcal{J}}\Gr(\ga(0)),e^{-\var\mathcal{J}}\Gr(\ga(\tau)))\nn\\&=&\<\ga,
      \aa\>\nn\\&=&\mu^{CLM}_F(V_1,\ga_1)+\mu^{CLM}_F(V_2,\ga_2)\lb{p3}\\
      &=& \mu^{CLM}_F(V_1,e^{-\var\mathcal{J}}\Gr(\ga))-\mu^{CLM}_F(V_2,e^{-\var\mathcal{J}}\Gr(\ga))\lb{p4}\\
      &=& i_{L_0}(\ga)-i_{L_1}(\ga),\lb{p14}
      \eea
where $\aa$ is the Maslov-Arnold index defined in Theorem 3.4.9 on
p. 64 of \cite{Dui}. Since $\ga_1$ and $\ga_2$ are transversal to
$V_1$ and $V_2$ (\ref{p3}) holds, (\ref{p4}) holds from the
definition of $\ga_1$ and $\ga_2$.

In the proof of Theorem 3.3 of \cite{LZZ}, we have proved that for
$\var>0$ small enough, there holds
 \be \sgn(V_1,e^{-\var\mathcal{J}}\Gr(I_{2n});V_2)=0,\lb{p5}\ee
 where $\sgn(W_1,W_3;W_2)$ for 3 Lagrangian spaces with $W_3$ transverses to $W_1$ and $W_2$
  is introduced in Definition 3.2.3 on p. 67 of \cite{Dui}. Note
  that by Claim 1 below, we can prove (\ref{p5}) at once.

{\bf Claim 1.} For $\var>0$, small enough, there holds
            \be \sign
            (V_1,e^{-\var\mathcal{J}}
\Gr(\ga(\tau));V_2)=\sgn(M_\var(\ga(\tau))).\lb{p12}\ee

{\bf Proof of Claim 1.} In fact,
     \be e^{-\mathcal{J}}\Gr
     (\ga(\tau))=\left\{\left(\begin{array}{cc} e^{\var J}&0\\
     0& e^{-\var J}\ga(\tau)\end{array}\right)\left(\begin{array}{c}
     p\\ q\\ p\\ q\end{array}\right)=\left(\begin{array}{c}
     cp-sq\\sp+cq\\ (c,s)\ga(\tau)(p,q)^T\\(-s,c)\ga(\tau)(p,q)^T\end{array}\right);\quad
     p,q\in\R^{n}\right\},\ee
 where we denote by $c=\cos\var I_n$ and $s=\sin\var I_n$.
 Hence the transformation $A: V_1\mapsto e^{-\mathcal{J}}\Gr
     (I_{2n})$ satisfies
      \bea&& A(0,-sp-cq, 0,-(-s,c)\ga(\tau)(p,q)^T)\nn\\
      &&=(cp-sq,sp+cq,(c,s)\ga(\tau)(p,q)^T,(-s,c)\ga(\tau)(p,q)^T),\quad
      \forall p,q\in \R^n,\lb{p8}\eea
where $A$ is introduced in Definition 3.4.3 of $\sign (M_1,M_2;L)$
on p. 67 of \cite{Dui}. For the convenience of our computation, we
rewrite (\ref{p8}) as follows.
 \bea && A\left(-\left(\begin{array}{cc}0&0\\s&c\end{array}\right)
 \left(\begin{array}{c}p\\q\end{array}\right),- \left(\begin{array}{cc}0&0\\-s&c\end{array}\right)
  \ga(\tau)\left(\begin{array}{c}p\\q\end{array}\right)\right)\nn\\&&=\left(\left(\begin{array}{cc}c&-s\\s&c
  \end{array}\right)\left(\begin{array}{c}p\\q\end{array}\right),\left(\begin{array}{cc}c&s\\-s&c
  \end{array}\right)\ga(\tau)\left(\begin{array}{c}p\\q\end{array}\right)\right).\eea

Then for $p_1,p_2,q_1,q_2\in\R^n$, the symmetric bilinear form
$Q(V_2):(x,y)\mapsto \mathcal{J}(Ax,y)$ on $V_1$ defined in
Definition 3.4.3 on p. 67 of \cite{Dui} satisfies:
 \bea && Q(V_2)\left(\left(-\left(\begin{array}{cc}0&0\\s&c\end{array}\right)
 \left(\begin{array}{c}p\\q\end{array}\right),- \left(\begin{array}{cc}0&0\\-s&c\end{array}\right)
  \ga(\tau)\left(\begin{array}{c}p\\q\end{array}\right)\right)\right)\nn\\
  &=&\left\<((-J)\oplus J)\left[\left(\begin{array}{cc}c&-s\\s&c
  \end{array}\right)\left(\begin{array}{c}p\\q\end{array}\right),\left(\begin{array}{cc}c&s\\-s&c
  \end{array}\right)\ga(\tau)\left(\begin{array}{c}p\\q\end{array}\right)\right],\right.\;\nn\\
&&\left[-\left(\begin{array}{cc}0&0\\s&c\end{array}\right)
 \left.\left(\begin{array}{c}p\\q\end{array}\right),- \left(\begin{array}{cc}0&0\\-s&c\end{array}\right)
  \ga(\tau)\left(\begin{array}{c}p\\q\end{array}\right)\right]\right\>.\nn\\
  &=&\left\< \left[\left(\begin{array}{cc}0&s\\0&c
  \end{array}\right)J\left(\begin{array}{cc}c&-s\\s&c
  \end{array}\right)-\ga(\tau)^T\left(\begin{array}{cc}0&-s\\0&c
  \end{array}\right)J\left(\begin{array}{cc}c&s\\-s&c
  \end{array}\right)\ga(\tau)\right]\left(\begin{array}{c}p\\q\end{array}\right),\;
  \left(\begin{array}{c}p\\q\end{array}\right)\right\>. \nn \\
  &=& \left\< \left[\left(\begin{array}{cc}sc&-s^2\\c^2&-sc
  \end{array}\right)+\ga(\tau)^T\left(\begin{array}{cc}sc&s^2\\-c^2&-sc
  \end{array}\right)\ga(\tau)\right]\left(\begin{array}{c}p\\q\end{array}\right),\;
  \left(\begin{array}{c}p\\q\end{array}\right)\right\>.\eea

Let $\td{M}_\var(\ga(\tau))=\left(\begin{array}{cc}sc&-s^2\\c^2&-sc
  \end{array}\right)+\ga(\tau)^T\left(\begin{array}{cc}sc&s^2\\-c^2&-sc
  \end{array}\right)\ga(\tau)$. Then by definition of the symmetric
  bilinear form $Q(V_2)$, $\td{M}_\var(\ga(\tau)$ is an invertible symmetric
  $2n\times 2n$ matrix. We define
   \be M_\var(\ga(\tau))=2\td{M}_\var(\ga(\tau))=\td{M}_\var(\ga(\tau))+\td{M}_\var^T(\ga(\tau)).\ee
 Then we have
 \bea M_\var(\ga(\tau))=\ga(\tau)^T\left(\begin{array}{cc}\sin{2\var}I_n&-\cos{2\var I_n}\\-\cos{2\var}I_n&-\sin
 2\var I_n
 \end{array}\right)\ga(\tau)+\left(\begin{array}{cc}\sin{2\var}I_n&\cos{2\var}I_n\\\cos{2\var}I_n&-\sin2\var
 I_n
 \end{array}\right).\eea
It is clear that
 \be \sgn Q(V_2)=\sgn \td{M}_\var(\ga(\tau))=\sgn M_\var(\ga(\tau)).\lb{p10}\ee

 By the definition of $\sgn(V_1,e^{-\var\mathcal{J}}
\Gr(\ga(\tau));V_2)$, we have
 \be \sgn(V_1,e^{-\var\mathcal{J}}
\Gr(\ga(\tau));V_2)=\sgn Q(V_2).\lb{p11}\ee

Then (\ref{p12}) holds from (\ref{p10}) and (\ref{p11}), and the
proof of Claim 1 is complete.

 Thus by (\ref{p14}), (\ref{p5}) and Claim 1, we have
\bea  && i_{L_0}(\ga)-i_{L_1}(\ga)\nn\\
&=&
s(V_1,V_2;e^{-\var\mathcal{J}}\Gr(\ga(0)),e^{-\var\mathcal{J}}\Gr(\ga
                       (\tau)))\nn\\
&=&\frac{1}{2}\sgn(V_1,e^{-\var\mathcal{J}}\Gr(\ga(\tau));V_2)-\frac{1}{2}\sgn(V_1,e^{-\var\mathcal{J}}\Gr(\ga(0));V_2)
                       \nn\\
&=&\frac{1}{2}\sgn(V_1,e^{-\var\mathcal{J}}\Gr(\ga(\tau));V_2)-\frac{1}{2}\sgn(V_1,e^{-\var\mathcal{J}}\Gr(I_{2n});V_2)\nn\\
 &=&\frac{1}{2}\sgn(V_1,e^{-\var\mathcal{J}}\Gr(\ga(\tau));V_2)\nn\\
 &=& \frac{1}{2}\sgn M_\var(\ga(\tau)).\nn\eea
 Here in the second equality, we have used Theorem 3.4.12 of on p. 68
of \cite{Dui}. Thus (\ref{p20}) holds.

Choose $\var<0$ such that $|\var|$ is sufficiently small, by the
discussion of $\mu^{CLM}_{F}$ index we have
\bea && i_{L_0}(\ga)=\mu^{CLM}_{F}(V_1, e^{-\var\mathcal{J}}\Gr(\ga),[0,\tau])-\nu_{L_0}(\ga),\\
     && i_{L_1}(\ga)=\mu^{CLM}_{F}(V_2, e^{-\var\mathcal{J}}\Gr(\ga),[0,\tau])-\nu_{L_1}(\ga).\eea
Then by the same proof as above, we have
 \be
 i_{L_0}(\ga)+\nu_{L_0}(\ga)-i_{L_1}(\ga)-\nu_{L_1}(\ga)=\frac{1}{2}\sgn
 M_\var(\ga(\tau)),\ee
where $\var<0$ is small enough. Hence (\ref{p21}) holds. The proof
of Theorem 2.3 is complete. \hfill\hb

We have the following consequence.

 \noindent{\bf Corollary 2.1.} (Theorem 2.3 of \cite{LZ}) {\it
 For $\ga\in \mathcal{P}_\tau(2n)$ with $\tau>0$, there hold
\bea |i_{L_0}(\ga))-i_{L_1}(\ga))|\le n,\quad
|i_{L_0}(\ga)+\nu_{L_0}(\ga)-i_{L_1}(\ga)-\nu_{L_1}(\ga)|\le
n.\lb{p24}\eea Moreover if $\ga(1)$ is a orthogonal matrix then
there holds
 \be i_{L_0}(\ga)=i_{L_1}(\ga).\lb{p23}\ee  }

{\bf Proof.} (\ref{p24}) holds directly from Theorem 2.3, so we only
need to prove (\ref{p23}). Since $\ga(\tau)$ is an orthogonal and
symplectic matrix, we have
    \be \ga^T(\tau)J\ga(\tau)=J,\quad
    \ga^T(\tau)\ga(\tau)=I_{2n}.\ee
So we have
  \be \ga(\tau)J=J\ga(\tau),\quad
  \ga(\tau)^TJ=J\ga(\tau)^T.\lb{p25}\ee

It is easy to check that for any $\var\in\R$, there holds
 \be J\left(\begin{array}{cc}\sin{2\var}I_n&\pm\cos{2\var I_n}\\\pm \cos{2\var}I_n&-\sin
 2\var I_n
 \end{array}\right)J=\left(\begin{array}{cc}\sin{2\var}I_n&\pm\cos{2\var I_n}\\\pm\cos{2\var}I_n&-\sin
 2\var I_n
 \end{array}\right).\lb{p26}\ee
Hence by (\ref{p25}) and (\ref{p26}), we have
 \bea JM_\var(\ga(\tau)) J&=&J\left[\ga(\tau)^T\left(\begin{array}{cc}\sin{2\var}I_n&-\cos{2\var I_n}\\-\cos{2\var}I_n&-\sin
 2\var I_n
 \end{array}\right)\ga(\tau)+\left(\begin{array}{cc}\sin{2\var}I_n&\cos{2\var}I_n\\\cos{2\var}I_n&-\sin2\var
 I_n
 \end{array}\right)\right]J\nn\\
 &=& J\ga(\tau)^T\left(\begin{array}{cc}\sin{2\var}I_n&-\cos{2\var I_n}\\-\cos{2\var}I_n&-\sin
 2\var I_n
 \end{array}\right)\ga(\tau)J+J\left(\begin{array}{cc}\sin{2\var}I_n&\cos{2\var}I_n\\\cos{2\var}I_n&-\sin2\var
 I_n
 \end{array}\right)J\nn\\
 &=& \ga(\tau)^TJ\left(\begin{array}{cc}\sin{2\var}I_n&-\cos{2\var I_n}\\-\cos{2\var}I_n&-\sin
 2\var I_n
 \end{array}\right)J\ga(\tau)+J\left(\begin{array}{cc}\sin{2\var}I_n&\cos{2\var}I_n\\\cos{2\var}I_n&-\sin2\var
 I_n
 \end{array}\right)J\nn\\
&=& \ga(\tau)^T\left(\begin{array}{cc}\sin{2\var}I_n&-\cos{2\var
I_n}\\-\cos{2\var}I_n&-\sin
 2\var I_n
 \end{array}\right)\ga(\tau)+\left(\begin{array}{cc}\sin{2\var}I_n&\cos{2\var}I_n\\\cos{2\var}I_n&-\sin2\var
 I_n
 \end{array}\right)\nn\\
 &=& M_\var(\ga(\tau)).\eea
So we have
 \be M_\var(\ga(\tau)) J=-JM_\var(\ga(\tau)).\lb{p30}\ee
Thus for any $x\in\R^{2n}$ and $\lm\in \R$ satisfying
 \be M_\var(\ga(\tau)) x=\lm x.\ee
By (\ref{p30}) we have \be M_\var(\ga(\tau)) (J
x)=-JM_\var(\ga(\tau)) x=-\lm (Jx).\lb{p31}\ee

 Since for $\var>0$ small enough $M_\var(\ga(\tau))$ is an invertible symmetric
matrix, by (\ref{p31}) we have \be
m^+(M_\var(\ga(\tau)))=m^-(M_\var(\ga(\tau)))=n\ee which yields \be
\sgn
M_\var(\ga(\tau))=m^+(M_\var(\ga(\tau)))-m^-(M_\var(\ga(\tau)))=0.\ee
Then (\ref{p23}) holds from Theorem 2.3. \hfill\hb

\noindent{\bf Lemma 2.2.} {\it For a symplectic path $P: [0,
\tau]\to \Sp(2n)$ with $\tau>0$, if for $j=0,1$ there holds
$\nu_{L_j}(P(t))=constant$ for all $t\in[0,\tau]$, then for $\var>0$
small enough we have
       \be \sgn M_\var(P(0))=\sgn M_\var(P(\tau)).\lb{zhi0}\ee}
\quad{\bf Proof.} Since $\Sp(2n)$ is path connected, we can choose a
path $\ga\in\mathcal{P}_\tau$ with $\ga(\tau)=P(0)$. By Proposition
2.11 of \cite{LZZ} and the definition of $\mu_j$ for $j=1,2$ in
\cite{LZZ}, we have
      \be \mu_F^{CLM}(V_j,\Gr(P), [0,\tau])=0,\qquad j=0,1.\ee
 So by the Path Additivity and Reparametrization Invariance properties of $\mu_F^{CLM}$ in \cite{CLM}, we
 have
     \bea i_{L_j}(P*\ga)&=& \mu_F^{CLM}(V_j,\Gr(P*\ga),[0,\tau])-n\nn\\
                        &=&\mu_F^{CLM}(V_j,\Gr(\ga),[0,\tau])+\mu_F^{CLM}(V_j,\Gr(P),[0,\tau])-n\nn\\
                        &=&\mu_F^{CLM}(V_j,\Gr(\ga),[0,\tau])-n\nn\\
                        &=& i_{L_j}(\ga),\lb{zhi1}\eea
 where the definition of joint path $\eta*\xi$ is given by (\ref{zzz}) in Section 6 below.
Then by Theorem 2.3 we have
 \bea && i_{L_0}(\ga)-i_{L_1}(\ga)=\frac{1}{2}\sgn (M_\var(P(0))),\lb{zhi2}\\
      && i_{L_0}(P*\ga)-i_{L_1}(P*\ga)=\frac{1}{2}\sgn (M_\var(P(\tau))).
 \lb{zhi3}\eea
Then (\ref{zhi0}) holds from (\ref{zhi1})-(\ref{zhi3}). The proof of
Lemma 2.2 is complete.\hfill\hb

\noindent{\bf Remark 2.1.} {\it It is easy to check that for
 $n_j \times n_j$ symplectic matrix $P_j$ with $j=1,2$ and
 $n_j\in\N$, we have
    \bea &&M_\var(P_1\diamond P_2)=M_\var(P_1)\diamond
    M_\var(P_2),\nn\\
     &&\sgn M_\var(P_1\diamond P_2)=\sgn M_\var(P_1)+
    \sgn M_\var(P_2).\nn\eea}

 By direct computation according to
Theorem 2.3 and Corollary 2.1, for $\ga\in\mathcal{P}_\tau(2)$,
$b>0$, and $\var>0$ small enough we have
 \bea
&& \sgn M_\var(R(\theta))=0,\quad
 {\rm for}\;\theta\in \R,\lb{guo1}\\
&&  \sgn M_\var(P)=0,\quad {\rm if}\; P=\pm
 \left(\begin{array}{cc}1&b\\0&1\end{array}\right)\;{\rm or}\;
 \pm
 \left(\begin{array}{cc}1&0\\-b&1\end{array}\right),\lb{guo2}\\
&&\sgn M_\var(P)=2,\quad {\rm if}\; P=\pm
 \left(\begin{array}{cc}1&-b\\0&1\end{array}\right),\lb{guo3.5}\\
 &&\sgn M_\var(P)=-2,\quad {\rm if}\; P=
 \pm
 \left(\begin{array}{cc}1&0\\b&1\end{array}\right).\lb{guo3}\eea
Also we give a example as follows to finish this section
 \be \sgn M_\var(P)=2,\quad {\rm if}\; P=\pm
 \left(\begin{array}{cc}2&-1\\-1&1\end{array}\right).\lb{guo4}\ee

\setcounter{equation}{0}
\section{Relation between $i_{L_0}$, $i_{L_1}$, $i_{\sqrt{-1}}^{L_0}$ and the corresponding Morse indices,
and their monotonicity properties. }

 In \cite{Liu2}, Liu studied the relation between the $L$-index of solutions of Hamiltonian systems with
 $L$-boundary conditions and the Morse index of the corresponding
 functional defined via the Galerkin approximation method on the
 finite dimensional truncated space at its corresponding critical
 points. In order to prove the main results of this paper, in this section
 we use the results of \cite{Liu2} to study
 some monotonicity properties of $i_{L_0}$ and $i_{L_1}$. We also
 study the index $i_{\sqrt{-1}}^{L_0}(B)$ with $B$ being a continuous
 symmetric matrices path satisfying condition (B1) defined in Section 1 and the Morse index of the corresponding
 functional defined via the Galerkin approximation method. Then as
 applications we study some monotonicity properties of
 $i_{\sqrt{-1}}^{L_0}(B)$ which will be important in the proof of
 Theorems 1.4-1.5 in Section 5 below.

For any $\tau>0$ and $B\in C([0,\tau/4],\mathcal{L}_s(\R^{2n}))$ (in
order to apply the results in this section conveniently Section 5,
we always assume $B\in C([0,\tau/4],\mathcal{L}_s(\R^{2n})$)
satisfying condition (B1). We extend $B$ to $[0,\frac{\tau}{2}]$ by
 \be B(\frac{\tau}{4}+t)=NB(\frac{\tau}{4}-t)N, \;\forall
 t\in[0,\frac{\tau}{4}].\ee
 Then since $B(\frac{\tau}{2})=B(0)$, we can extend it
 $\frac{\tau}{2}$-periodically to $\R$, so we can see $B$ as an element in $
 C(S_{\tau/2},\mathcal{L}_s(\R^{2n}))$.

   Let
$E_{\tau}=\{x\in
W^{1/2,2}(S_{\tau},\R^{2n})|\,x(-t)=Nx(t)\;a.e.\;t\in \R\}$ with the
usual norm and inner product denoted by $||\cdot||$ and
$\langle\cdot\rangle$ respectively.

By the Sobolev embedding theorem, for any $s\in[1,+\infty)$, there
is a constant $C_s>0$ such that
     \be ||z||_{L^s}\le C_s||z||,\quad \forall z\in
     E_{2\tau}.\lb{feng}\ee

Note that $B$ can also be seen as an element in
$C(S_{\tau},\mathcal{L}_s(\R^{2n}))$. We define two selfadjoint
operators $A_{\tau}$ and $B_{\tau}$ on $E_{\tau}$ by the following
bilinear forms

\bea \langle A_{\tau}x,y\rangle=\int_0^{\tau}-J\dot{x}\cdot
y\,dt,\qquad \langle B_{\tau}x,y\rangle=\int_0^{\tau}B(t)x\cdot y
\,dt.\lb{he8}\eea

Then $A_{\tau}$ is a bounded operator on $E_{\tau}$ and dim $\ker
A_{\tau}=n$, the Fredholm index of $A_{\tau}$ is zero, and
$B_{\tau}$ is a compact operator on $E_{\tau}$.

Set
       $$E_{\tau}(j)=\left\{z\in E_\tau\left|z(t)=\exp(\frac{2j\pi
       t}{\tau}J)a+\exp(-\frac{2j\pi t}{\tau}J)b,\; \forall t\in \R;\; \forall a,\,b\in L_0\right.\right\}.$$ and
       $$E_{\tau, m}=E_{\tau}(0)+E_{\tau}(1)+\cdots +E_{\tau}(m).$$
Let $\Ga_{\tau}=\{P_{\tau,m}: m=0,1,2,...\}$ be the usual Galerkin
approximation scheme w.r.t. $A_{\tau}$, just as in \cite{Liu2},
i.e., $\Ga_{\tau}$ is a sequence of orthogonal projections
satisfies:

  (1) $E_{\tau,0}=P_{\tau,0}E_{\tau}=\ker A_{\tau},\; E_{\tau,m}=P_{\tau,m}E_{\tau}$ is finite dimension for $m\ge 0$;

  (2) $P_{\tau,m}\to x$ as $m\to \infty$ for any $x\in E_\tau$;

  (3) $P_{\tau,m}A_{\tau}=A_{\tau} P_{\tau,m}$, $\forall m\ge 0$.

For $d>0$, we denote by $M^+_d(\cdot)$, $M^-_d(\cdot)$ and
$M^0_d(\cdot)$ the eigenspace corresponding to the eigenvalue $\lm$
belong to $[d,+\infty)$, $(-\infty, -d]$ and $(-d,d)$ respectively,
 and $M^+(\cdot)$,  $M^-(\cdot)$ and $M^0(\cdot)$ the positive,
 negative and null subspace of of the selfadjoint operator defining
 it respectively. For any bounded selfadjoint linear operator on
 $E$, We denote $L^\#=(L|_{Im L})^{-1}$, and we also denote by
 $P_{\tau,m}LP_{\tau,m}=(P_{\tau,m}LP_{\tau,m})|_{E_{\tau,m}}: E_{\tau,m}\to E_{\tau,m}$.

Similarly we define two subspaces of $E_\tau$ by $\hat{E}=\{x\in
E|x(t+\frac{\tau}{2})=-x(t), a.e.\, t\in \R\}$ and $\td{E}=\{x\in
E|x(t+\frac{\tau}{2})=x(t), a.e.\, t\in \R\}$ be the symmetric ones
and $\frac{\tau}{2}$-periodic ones of $E_\tau$ respectively.

We define two selfadjoint operators $\hat{A}$ and $\hat{B}$ on
$\hat{E}$ by the following bilinear forms

\bea \langle \hat{A}x,y\rangle=\int_0^{\tau}-J\dot{x}\cdot
y\,dt,\qquad \langle \hat{B}x,y\rangle=\int_0^{\tau}B(t)x(t)\cdot
y(t)\,dt.\eea

Then $\hat{A}$ is a bounded Fredholm operator on $\hat{E}$ and dim
$\ker \hat{A}=0$, the Fredholm index of $\hat{A}$ is zero. $\hat{B}$
is a compact operator on $\hat{E}$.

For any positive integer $m$,  we define
$$ \hat{E}_m=\Sg_{j=1}^{m}E_{\tau}(2j-1).$$
For $m\ge 1$, let $\hat{P}_m$ be the orthogonal projection from
$\hat{E}$ to $\hat{E}_m$. Then $\{\hat{P}_m\}$ is a Galerkin
approximation scheme w.r.t. $\hat{A}$.

 \noindent{\bf Theorem 3.1.} {\it For any $B(t)\in C([0,\frac{\tau}{4}], \mathcal{L}_s(\R^{2n}))$
 satisfying condition (B1)
 and $0<d\le \frac{1}{4}||(A_{\tau}-B_{\tau})^\#||^{-1}$, there
 exists $m^*>0$ such that for $m\ge m^*$ there hold
  \bea
  \dim M_d^+(\hat{P}_m(\hat{A}-\hat{B})\hat{P}_m)&=&mn-i_{\sqrt{-1}}^{L_0}(B)-\nu_{\sqrt{-1}}^{L_0}(B),\lb{d1}\\
  \dim M_d^-(\hat{P}_m(\hat{A}-\hat{B})\hat{P}_m)&=&mn+i_{\sqrt{-1}}^{L_0}(B),\lb{d2}\\
  \dim
  M_d^0(\hat{P}_m(\hat{A}-\hat{B})\hat{P}_m)&=&\nu_{\sqrt{-1}}^{L_0}(B).\lb{d3}\eea}
{\bf Proof.} The method of the proof here is similar as that of
Theorem 2.1 in \cite{Z1}.

For any positive integer $m$,  we define
$$ \td{E}_m=\sum_{j=0}^{m}E_{\tau}(2j).$$
For $m\ge 1$, let $\td{P}_m$ be the orthogonal projection from
$\td{E}$ to $\td{E}_m$. Then $\{\td{P}_m\}$ is a Galerkin
approximation scheme w.r.t. $\td{A}$.

For any $y\in \hat{E}_m$ and $z\in \td{E}_m$, it is easy to check
that
 \bea \langle (P_{\tau,m}(A_{\tau}-B_{\tau})P_{\tau,m}y,z)\rangle
 =0.\eea

So we have the following $P_{\tau,m}(A_{\tau}-B_{\tau})P_{\tau,m}$
orthogonal decomposition
   \be E_{\tau,2m}=\hat{E}_m \oplus \td{E}_m.\ee
Similarly, we have the following $A_{\tau}-B_{\tau}$ orthogonal
decomposition
    \be E_{\tau}=\hat{E} \oplus \td{E}.\ee
Hence, under above decomposition we have
 \be
(A_{\tau}-B_{\tau})=(\hat{A}-\hat{B})\oplus (\td{A}-\td{B}).\ee

Thus
 \bea ||(A_{\tau}-B_{\tau})^\#||^{-1}\le ||(\hat
{A}-\hat{B})^\#||^{-1}\lb{z1}\\
||(A_{\tau}-B_{\tau})^\#||^{-1}\le ||(\td
{A}-\td{B})^\#||^{-1}\lb{z2} \eea

By the definitions of $M_d^*(\cdot)$ for
$P_{\tau,2m}(A_{\tau}-B_{\tau})P_{\tau,2m}$,
$\hat{P}_m(\hat{A}-\hat{B})\hat{P}_m$, and
$\td{P}_m(\td{A}-\td{B})\td{P}_m$ with $*=+,-,0$. So for
$*\in\{+,-,0\}$ we have \be
        \dim M_d^*(P_{\tau,2m}(A_{\tau}-B_{\tau})P_{\tau,2m})=\dim M_d^*(\hat{P}_m(\hat{A}-\hat{B})\hat{P}_m)
+\dim M_d^*(\td{P}_m(\td{A}-\td{B})\td{P}_m).\lb{2.1}\ee

Note that, the space $E_\tau$ and the operators $A_\tau$, $B_{\tau}$
and $P_{\tau,m}$ are also defined in the same way. So by the
definition we see that $\td{E}$ is the $\tau$-periodic extending of
$E_\tau$ from $S_\tau$ to $S_{2\tau}$, and $\td{E}_m$ is the
$\tau$-periodic extending of $E_{\tau,2m}$ from $S_\tau$ to
$S_{2\tau}$ too.

Thus we have \be ||(A_{\tau}-B_{\tau})^\#||^{-1}= ||(\td
{A}-\td{B})^\#||^{-1}.\lb{z3}\ee

\noindent By (\ref{z2}) and (\ref{z3}) we have \be
||(A_{2\tau}-B_{2\tau})^\#||^{-1}\le
||(A_{\tau}-B_{\tau})^\#||^{-1}. \lb{z5}\ee For $*\in\{+,-,0\}$ we
have \be \dim M_d^*(P_{\tau,m}(A_\tau-B_\tau)
P_{\tau,m})=M_d^*(\td{P}_m(\td{A}-\td{B})\td{P}_m).\lb{z4}\ee

Then for $0<d\le \frac{1}{4}||(A_{\tau}-B_{\tau})^\#||^{-1}$, by
Theorem 2.1 in \cite{Liu2} there exists $m_1>0$ such that for $m\ge
m_1$ we have
 \bea &&\dim
M_d^+(P_{\tau,2m}(A_{\tau}-B_{\tau})P_{\tau,2m})=2mn-i_{L_0}(\ga_B^2)-\nu_{L_0}(\ga_B^2),\lb{z6}\\
&& \dim
M_d^-(P_{\tau,2m}(A_{\tau}-B_{\tau})P_{\tau,2m})=2mn+n+i_{L_0}(\ga_B^2),\lb{z7}\\
&&\dim
M_d^0(P_{\tau,2m}(A_{\tau}-B_{\tau})P_{\tau,2m})=\nu_{L_0}(\ga_B^2).\lb{z8}
\eea

By (\ref{z5}), we have $0<d\le
\frac{1}{4}||(A_{\tau}-B_{\tau})^\#||^{-1}$. By Theorem 2.1 in
\cite{Liu2} again there exists $m_2>0$, such that for $m\ge m_2$ we
have
  \bea &&\dim
M_d^+(P_{\tau,m}(A_{\tau}-B_{\tau})P_{\tau,m})=mn-i_{L_0}(\ga_B)-\nu_{L_0}(\ga_B)),\lb{z9}\\
&& \dim
M_d^-(P_{\tau,m}(A_{\tau}-B_{\tau})P_{\tau,m})=mn+n+i_{L_0}(\ga_B)),\lb{z10}\\
&&\dim
M_d^0(P_{\tau,m}(A_{\tau}-B_{\tau})P_{\tau,m})=\nu_{L_0}(\ga_B)).\lb{z11}
\eea

Let $m^*=\max\{m_1,m_2\}$. Then for $m\ge m^*$, all of
(\ref{z6})-(\ref{z11}) hold.

So by (\ref{2.1}), (\ref{z4}), and (\ref{z6})-(\ref{z11}) we have
  \bea
  \dim M_d^+(\hat{P}_m(\hat{A}-\hat{B})\hat{P}_m)&=&mn-(i_{L_0}(\ga_B^2)-i_{L_0}(\ga_B))-
  (\nu_{L_0}(\ga_B^2)-\nu_{L_0}(\ga_B)),\lb{z14}\\
  \dim M_d^-(\hat{P}_m(\hat{A}-\hat{B})\hat{P}_m)&=&mn+i_{L_0}(\ga_B^2)-i_{L_0}(\ga_B),\lb{z15}\\
  \dim M_d^0(\hat{P}_m(\hat{A}-\hat{B})\hat{P}_m)&=&\nu_{L_0}(\ga_B^2)-\nu_{L_0}(\ga_B).\lb{z16}\eea

Thus (\ref{d1})-(\ref{d3}) hold from (\ref{z14})-(\ref{z16}),
Definition 2.3, and Theorem 2.2. The proof of Theorem 3.1 is
complete. \hfill\hb

\noindent{\bf Remark 3.1.} {\it Let any $B\in C([0,\frac{\tau}{4}],
\mathcal{L}_s(\R^{2n}))$ be a constant
 matrix path satisfying condition (B1).  By Theorem 5.1
of \cite{LZZ}, for $d=0$ the same conclusions of Theorem 2.1 of
\cite{Liu2} still holds . Hence for $d=0$  the same conclusions of
Theorem 3.1 still hold, i.e., there
 exists $m^*>0$ such that for $m\ge m^*$ there hold
   \bea
  \dim M^+(\hat{P}_m(\hat{A}-\hat{B})\hat{P}_m)&=&mn-i_{\sqrt{-1}}^{L_0}(B)-\nu_{\sqrt{-1}}^{L_0}(B),\nn\\
  \dim M^-(\hat{P}_m(\hat{A}-\hat{B})\hat{P}_m)&=&mn+i_{\sqrt{-1}}^{L_0}(B),\nn\\
  \dim
  M^0(\hat{P}_m(\hat{A}-\hat{B})\hat{P}_m)&=&\nu_{\sqrt{-1}}^{L_0}(B).\nn\eea}

In the following, we study some monotonicity of the the Maslov-type
$i_{\sqrt{-1}}^{L_0}$
  index. In this paper, for any two symmetric matrices $B_1$ and $B_2$, we say  $B_1>B_2$ if $B_1-B_2$ is
  positive definite and we say $B_1\ge B_2$ if $B_1-B_2$ is semipositive.
  Similarly for two symmetric matrix paths $B_1$, $B_2\in C([0,\tau],\mathcal{L}_s(R^{2n}))$,
   we say  $B_1>B_2$ if  $B_1(t)-B_2(t)$ is
  positive definite for all $t\in [0,\tau]$ and we say $B_1\ge B_2$ if $B_1(t)-B_2(t)$ is
  semipositive definite for all $t\in [0,\tau]$.

\noindent{\bf Lemma 3.1.} {\it For any $\tau>0$ and $B_1,\; B_2\in
C([0,\frac{\tau}{4}],\mathcal{L}_s(\R^{2n}))$ satisfying condition
(B1).
       If $B_1\ge B_2$,  then there hold
\be i_{\sqrt{-1}}^{L_0}(B_1)\ge i_{\sqrt{-1}}^{L_0}(B_2)\lb{w4}\ee
and
      \bea i_{\sqrt{-1}}^{L_0}(B_1)+  \nu_{\sqrt{-1}}^{L_0}(B_1)\ge i_{\sqrt{-1}}^{L_0}(B_2)+\nu_{\sqrt{-1}}^{L_0}(B_2).\lb{www}\eea
   Moreover, if \be \int_0^{\frac{\tau}{4}}(B_1(t)-B_2(t))dt>0,\ee then there
          holds
\bea i_{\sqrt{-1}}^{L_0}(B_1)\ge i_{\sqrt{-1}}^{L_0}(B_2)+
\nu_{\sqrt{-1}}^{L_0}(B_2).\lb{w6}\eea}

{\bf Proof.} Let the space $\hat{E}$ and the orthogonal projection
operator $\hat{P}_m$ be the ones defined in Section 2.
Correspondingly we define the compact operators $\hat{B}_1$ and
$\hat{B}_2$.  By Theorem 3.1, for $d>0$ small enough, there exists
$m^*>0$ such that
        \bea
  \dim M_d^+(\hat{P}_m(\hat{A}-\hat{B}_1)\hat{P}_m)&=&mn-i_{\sqrt{-1}}^{L_0}(B_1)-\nu_{\sqrt{-1}}^{L_0}(B_1),\lb{w1}\\
  \dim M_d^-(\hat{P}_m(\hat{A}-\hat{B}_1)\hat{P}_m)&=& mn+i_{\sqrt{-1}}^{L_0}(B_1),\lb{w2}\\
  \dim
  M_d^0(\hat{P}_m(\hat{A}-\hat{B}_1)\hat{P}_m)&=&\nu_{\sqrt{-1}}^{L_0}(B_1).\lb{w3}\eea
and
   \bea
  \dim M_d^+(\hat{P}_m(\hat{A}-\hat{B_2})\hat{P}_m)&=&mn-i_{\sqrt{-1}}^{L_0}(B_2)-\nu_{\sqrt{-1}}^{L_0}(B_2),\lb{r1}\\
  \dim M_d^-(\hat{P}_m(\hat{A}-\hat{B}_2)\hat{P}_m)&=&mn+i_{\sqrt{-1}}^{L_0}(B_2),\lb{r2}\\
  \dim
  M_d^0(\hat{P}_m(\hat{A}-\hat{B}_2)\hat{P}_m)&=&\nu_{\sqrt{-1}}^{L_0}(B_2).\lb{r3}\eea
If $B_1\ge B_2$, we have
$\hat{P}_m(\hat{A}-\hat{B}_1)\hat{P}_m\le\hat{P}_m(\hat{A}-\hat{B}_2)\hat{P}_m$,
So \be \dim M_d^-(\hat{P}_m(\hat{A}-\hat{B}_1)\hat{P}_m)\ge \dim
M_d^-(\hat{P}_m(\hat{A}-\hat{B}_2)\hat{P}_m).\ee
 Then by (\ref{w2}) and (\ref{r2}), (\ref{w4}) holds. Also we have
\be \dim M_d^+(\hat{P}_m(\hat{A}-\hat{B}_1)\hat{P}_m)\le \dim
M_d^+(\hat{P}_m(\hat{A}-\hat{B}_2)\hat{P}_m).\ee Then by (\ref{w1})
and (\ref{r1}), (\ref{www}) holds.

If $\int_0^{\frac{\tau}{4}}(B_1(t)-B_2(t))dt>0$, then
 \be\hat{P}_m(\hat{A}-\hat{B}_1)\hat{P}_m<\hat{P}_m(\hat{A}-\hat{B}_2)\hat{P}_m.\ee
So we have
  \be \dim M_d^-(\hat{P}_m(\hat{A}-\hat{B}_1)\hat{P}_m)\ge \dim
M_d^-(\hat{P}_m(\hat{A}-\hat{B}_2)\hat{P}_m)+M_d^0(\hat{P}_m(\hat{A}-\hat{B}_2)\hat{P}_m).\ee
Then by (\ref{w2}), (\ref{r2}) and (\ref{r3}), (\ref{w6}) holds and
the proof of Lemma 3.1 is complete.\hfill\hb

\noindent{\bf Corollary 3.1.} {\it For any $\tau>0$ and $B\in
C([0,\frac{\tau}{4}],\mathcal{L}_s(\R^{2n}))$ satisfying condition
(B1) and $B\ge 0$, there holds
           \be i_{\sqrt{-1}}^{L_0}(B)\ge 0.\ee}
{\bf proof.} By Lemma 3.1, we have \be i_{\sqrt{-1}}^{L_0}(B)\ge
i_{\sqrt{-1}}^{L_0}(0).\ee Then the conclusion holds from the fact
that
  \be i_{\sqrt{-1}}^{L_0}(0)=i_{\sqrt{-1}}^{L_0}(\ga_0)=0,\ee
  Where $\ga_0$ is the identity symplectic path.\hfill\hb

By Theorem 2.1 of \cite{Liu2} and the Remark below Theorem 2.1 in
\cite{Liu2} and the similar proof of Lemma 3.1 we have the following
lemma.

\noindent {\bf Lemma 3.2.} {\it If $\tau>0$ and $B_1,\; B_2\in
C([0,\frac{\tau}{4}],\mathcal{L}_s(\R^{2n}))$ satisfying condition
(B1) and $B_1\ge B_2$, then for $j=0,1$ there hold \be
i_{L_j}(B_1)\ge i_{L_j}(B_2)\ee and
      \bea i_{L_j}(B_1)+  \nu_{L_j}(B_1)\ge i_{L_j}(B_2)+\nu_{L_j}(B_2).\lb{w5}\eea
   Moreover, if $\int_0^{\frac{\tau}{4}} (B_1(t)-B_2(t))dt>0$, then there
          holds
\bea i_{L_j}(B_1)\ge i_{L_j}(B_2)+ \nu_{L_j}(B_2).\eea}

 Since $i_{L_j}(0)=-n$ and $\nu_{L_j}(0)=n$ for $j=0,1$, a direct consequence of Lemma 3.2 is the following

\noindent{\bf Corollary 3.2.} {\it If $\tau>0$ and $B\in
C([0,\frac{\tau}{2}],\mathcal{L}_s(\R^{2n}))$ satisfying condition
(B1) and $B\ge 0$, then for $j=0,1$ there hold \be
i_{L_j}(B)+\nu_{L_j}(B)\ge 0,\qquad i_{L_j}(B)\ge -n.\ee Moreover if
$\int_0^\frac{\tau}{2} B(t)dt>0$, there holds \be i_{L_j}(B)\ge
0.\ee} \quad Moreover we can give a stronger version of Corollary
3.2, i.e., the following Lemma 3.3.

\noindent{\bf Lemma 3.3.} {\it Let $\tau>0$ and $B\in
C([0,\frac{\tau}{2}],\mathcal{L}_s(\R^{2n}))$ with the $n\times n$
matrix square block form  $B(t)=\left(
\begin{array}{cc}B_{11}(t)&B_{12}(t)\\B_{21}(t)&B_{22}(t)\end{array}\right)$ satisfying condition
(B1) and $B\ge 0$.

If $\int_0^{\frac{\tau}{2}} B_{22}(t)dt>0$, there holds \be
i_{L_0}(B)\ge 0.\ee

If $\int_0^{\frac{\tau}{2}} B_{11}(t)dt>0$, there holds \be
i_{L_1}(B)\ge 0.\ee}

{\bf Proof.} Without loss of generality, assume $\lm>0$ such that
    \be \int_0^{\frac{\tau}{2}}B_{22}(t)\ge \lm I_n.\ee
Also we can extend $B$ to $[0,\tau]$ by
 \be B(\frac{\tau}{2}+t)=NB(\frac{\tau}{2}-t)N, \;\forall
 t\in[0,\frac{\tau}{2}].\ee
 Then since $B(\tau)=B(0)$, we can extend it
 $\tau$-periodically to $\R$, so we can see $B$ as an element in $
 C(S_{\tau},\mathcal{L}_s(\R^{2n}))$.
 Then we have
 \be \int_0^\tau B_{22}(t)\ge 2\lm I_n.\ee

For any $m\in \N$, we define two subspaces of $E$ as follows
$$E^-_{\tau,m}=\left\{z\in
E_\tau\left|z(t)=\sum_{j=1}^m\exp(-\frac{2j\pi t}{\tau}J)b_j,\;
\forall t\in \R;\; \forall b_j\in L_0\right.\right\},$$
$$E_{\tau}(0)=\left\{z\in E_\tau\left|z(t)\equiv b,\; b\in
L_0\right.\right\}.$$

Then for any $z=\alpha x+\beta y \in E_{\tau}(0)\oplus E^-_{\tau,m}$
with $\alpha^2+\beta^2=1$ and $||x||=||y||=1$, we have
 \bea \langle(A_\tau-B_\tau)z,z\rangle &=& \langle(A_\tau-B_\tau)(\alpha x+\beta y),\alpha x+\beta
 y\rangle\nn\\
 &=& -\beta^2\langle A_\tau y,y\rangle-\langle B_\tau(\alpha x+\beta y), \alpha x+\beta
 y\rangle\nn\\
 &\le & -||A_\tau^\#||^{-1}\beta^2-\langle B_\tau(\alpha x+\beta y), \alpha x+\beta
 y\rangle.\lb{liang1}\eea
Since $B\ge 0$, note that $x(t)\equiv b=(0,b_1)\in L_0$ for all
$t\in S_\tau$ with $\tau|b_1|^2=1$, we have \bea &&\langle
B_\tau(\alpha x+\beta y),\alpha
  x+\beta y\rangle\nn\\
  &=& \int_0^\tau (\alpha^2 Bx\cdot x+\beta^2 By\cdot y+2\alpha\beta Bx\cdot
  y)\,dt\nn\\
  &\ge& \alpha^2\int_0^\tau  Bx\cdot x\,dt+\beta^2\int_0^\tau  By\cdot y\,dt-2|\alpha||\beta|(\int_0^\tau  Bx\cdot x\,dt)^{1/2}
   (\int_0^\tau  By\cdot y\,dt)^{1/2}\nn\\
   &\ge&\alpha^2\int_0^\tau  Bx\cdot x\,dt+\beta^2\int_0^\tau  By\cdot y\,dt-
   \frac{1}{1+\var}\alpha^2\int_0^\tau  Bx\cdot x\,dt-(1+\var)\beta^2\int_0^\tau  By\cdot y\,dt\nn\\
   &=&\frac{\var\alpha^2}{1+\var}\int_0^\tau  Bx\cdot x\,dt-\var\beta^2\int_0^\tau  By\cdot y\,dt\nn\\
    &=&\frac{\var\alpha^2}{1+\var}\left(\int_0^\tau  B(t)dt\right)b\cdot
    b-\var\beta^2\int_0^\tau  By\cdot y\,dt\nn\\
    &=&\frac{\var\alpha^2}{1+\var}\left(\int_0^\tau  B_{22}(t)dt\right)b_1\cdot
    b_1-\var\beta^2\int_0^\tau  By\cdot y\,dt\nn\\
   &\ge&\frac{\var\alpha^2}{1+\var}2\lm |b_1|^2-
     \var\beta^2||B_\tau||\,||y||^2\nn\\
   &=& \frac{2\var\lm\alpha^2}{(1+\var)\tau}-\var\beta^2||B_\tau||\lb{liang2}    \eea
for any $\var>0$.

Let $\var=\min\{1,\frac{||A_\tau^\#||^{-1}||B_\tau||^{-1}}{2}\}$. By
(\ref{liang1}) and (\ref{liang2}), we have
 \bea \langle(A_\tau-B_\tau)z,z\rangle &\le
 &-||A_\tau^\#||^{-1}\beta^2-\frac{2\var\lm\alpha^2}{(1+\var)\tau}+\var\beta^2||B_\tau||\nn\\
 &\le &
 -\frac{||A_\tau^\#||^{-1}\beta^2}{2}-\frac{\var\lm\alpha^2}{\tau}\nn\\
 &\le& -d_0(\alpha^2+\beta^2)\nn\\
& =&-d_0,\lb{liang3}\eea
 where $d_0=\min\{\frac{||A_\tau^\#||^{-1}}{2},\,\frac{\var\lm}{\tau}\}
 =\min\{\frac{||A_\tau^\#||^{-1}}{2},\,\frac{\lm}{\tau},\,\frac{\lm||A_\tau^\#||^{-1}||B_\tau||^{-1}}{2\tau}\}$.
Note that $d_0$ is independent of $m$, so for
$0<d\le\min\{d_0,\frac{||(A_\tau-B_\tau)^\#||^{-1}}{4}\}$, by
Theorem 2.1 of \cite{Liu2} there exists $m^*>0$ such that, for $m\ge
m^*$, we have
 \be \dim
 M^-_d(P_{\tau,m}(A_\tau-B_\tau)P_{\tau,m})=mn+n+i_{L_0}(B).\lb{liang4}\ee
By (\ref{liang3}) we have
  \bea \dim
 M^-_d(P_{\tau,m}(A_\tau-B_\tau)P_{\tau,m})\ge \dim (E_{\tau}(0)\oplus
 E^-_{\tau,m})=mn+n.\lb{liang5}\eea
Then by (\ref{liang4}) and (\ref{liang5}) we have $i_{L_0}(B)\ge 0$.

For $\int_0^{\frac{\tau}{2}} B_{11}(t)dt>0$,  by similar proof we
have  $i_{L_1}(B)\ge 0$. The proof of Lemma 3.3 is complete.
\hfill\hb

 Now we give the following Theorem 3.2 which will play
a important role in the proof of our main results in Section 5. This
results implies that the corresponding Maslov-type index of a
periodic symmetric solution of a first order even semipositive
Hamilton increases with the increasing of the iteration time of the
solution.

\noindent {\bf Theorem 3.2.} {\it If $\tau>0$ and $B\in
C([0,\frac{\tau}{4}],\mathcal{L}_s(\R^{2n}))$ satisfying condition
(B1) and $B\ge 0$, then for any two positive integers $p> q$ there
holds \be i_{\sqrt{-1}}^{L_0}(\ga_B^{p})\ge
i_{\sqrt{-1}}^{L_0}(\ga_B^{q}).\lb{zh0}\ee}

\noindent{\bf Proof.} Extend $\ga_B(t)$ to $[0,\frac{p\tau}{4}]$ as
$\ga_B^p$, we still denote it by $\ga_B$. By definition of
$i_{\sqrt{-1}}^{L_o}$ and the Path additivity and Symplectic
invariance property of $\mu_F^{CLM}$ in \cite{CLM}, we have
 \bea &&i_{\sqrt{-1}}^{L_0}(\ga_B^{p})-i_{\sqrt{-1}}^{L_0}(\ga_B^{q})\nn\\
         &=& \mu_F^{CLM}(L_0\times J L_0,\Gr(\ga_B), [0,\frac{p\tau}{4}])-
         \mu_F^{CLM}(L_0\times J L_0,\Gr(\ga_B), [0,\frac{q\tau}{4}])\nn\\
         &=&\mu_F^{CLM}(L_0\times J L_0,\Gr(\ga_B),
         [\frac{q\tau}{4},\frac{p\tau}{4}])\nn\\
         &=& \mu_F^{CLM}(L_0\times L_0, \Gr(-J\ga_B),
         [\frac{q\tau}{4},\frac{p\tau}{4}]).\lb{liu1}\eea

By the first geometrical definition of the index $\mu_F^{CLM}$ in
section 4 of \cite{CLM}, there is a $\var>0$ small enough such
  that
  \bea (e^{-\var \mathcal{J}}\Gr (-J\ga_B(\frac{p\tau}{4}))\cap (L_0\times L_0)=\{0\}=
  (e^{-\var J}\Gr (\ga_B(\frac{q\tau}{4}))\cap
  (L_0\times L_0)\lb{wang2}\eea
and
   \bea &&\mu_F^{CLM}(L_0\times L_0, \Gr(-J\ga_B),
         [\frac{q\tau}{4},\frac{p\tau}{4}])\nn\\
         &=& \mu_F^{CLM}(L_0\times L_0, e^{-\var\mathcal{J}}\Gr(-J\ga_B),
         [\frac{q\tau}{4},\frac{p\tau}{4}])\nn\\
         &=&  \mu_F^{CLM}(L_0\times L_0, \Gr(-e^{-\var J}J\ga_Be^{-\var J}),
         [\frac{q\tau}{4},\frac{p\tau}{4}]),\lb{wang1}\eea
where in the second equality we have used Symplectic invariance
property of $\mu_F^{CLM}$ index in \cite{CLM}. Choose a $C^1$ path
$\ga\in\mathcal{P}_{\frac{p\tau}{4}}$ such that
 $\ga(t)=-e^{-\var J}J\ga_Be^{-\var J}$ for all $t\in[\frac{q\tau},\frac{p\tau}{4}]$.
Denote by $D(t)=-J\dot{\ga}(t)\ga(t)^{-1}$ for
$t\in[0,\frac{p\tau}{4}]$.
 For $t\in[\frac{q\tau},\frac{p\tau}{4}]$, by direct computation we have
   \bea D(t)=-J\frac{d}{dt}(-e^{-\var J}J\ga e^{-\var J})(-e^{-\var J}J\ga e^{-\var
   J})^{-1}=-Je^{-\var J}B(t)e^{\var J}J.\lb{wang9}\eea
 Since $B\ge 0$ we have $D(t)\ge 0$ for $t\in[q\tau,p\tau]$ and $D\in
C([0,\frac{p\tau}{4}],\mathcal{L}_s(\R^{2n}))$. For $s\ge 0$, we
define $D_s(t)=D(t)+sI_{2n}$ and symplectic path $\ga_s(t)$ by
 \bea &&\frac{d}{dt}\ga_s(t)=JD_s(t)\ga_s(t),\quad
 t\in[0,\frac{p\tau}{4}]\nn\\
  &&\ga_s(0)=I_{2n}.\nn\eea
It is clear that
 \be \ga_0=\ga.\lb{wang7}\ee

 By the same argument of step2 of the proof of Theorem 5.1 in
\cite{LZZ}, we have
 \bea &&-J\frac{d}{ds}\ga_s(t)(\ga_s(t))^{-1}>0,\quad {\rm for}\, t=\frac{p\tau}{4},\frac{q\tau}{4}.\lb{wang3}\eea
By (\ref{wang2}) and definition of $\ga_s$ we have
 \be \nu_{L_0}(\ga_0(\frac{p\tau}{4}))=0=\nu_{L_0}(\ga_0(\frac{q\tau}{4})).\ee
So by (\ref{wang3}), there is a $\sg>0$ small enough such that
  \be \nu_{L_0}(\ga_s(\frac{p\tau}{4}))=0=\nu_{L_0}(\ga_s(\frac{q\tau}{4})),\quad \forall s\in[0,\sg].\lb{wang5}\ee
So we have
  \bea \mu_F^{CLM}(L_0\times L_0, \Gr(\ga_s(\frac{p\tau}{4})), s\in [0,\sg])=0,\nn\\
        \mu_F^{CLM}(L_0\times L_0, \Gr(\ga_s(\frac{q\tau}{4})), s\in
        [0,\sg])=0.\lb{wang6}\eea
By the Homotopy invariance with respect to end points and Path
additivity properties of $\mu_F^{CLM}$ index in \cite{CLM}, we have
 \bea &&\mu_F^{CLM}(L_0\times L_0, \Gr(\ga_s(\frac{p\tau}{4})), s\in [0,\sg])+\mu_F^{CLM}(L_0\times L_0a,
 \Gr(\ga_\sg(t)), t\in
 [\frac{q\tau}{4},\frac{p\tau}{4}])\nn\\
 &=&\mu_F^{CLM}(L_0\times L_0, \Gr(\ga_0(t)), t\in
 [\frac{q\tau}{4},\frac{p\tau}{4}])+\mu_F^{CLM}(L_0\times L_0, \Gr(\ga_s(\frac{p\tau}{4})), s\in
 [0,\sg]).\lb{wang10}\eea
So by (\ref{liu1}), (\ref{wang1}), (\ref{wang7}),(\ref{wang6}) and
(\ref{wang10}), we have
 \be i_{\sqrt{-1}}^{L_0}(\ga_B^{p})-i_{\sqrt{-1}}^{L_0}(\ga_B^{q})=\mu_F^{CLM}(L_0\times L_0, \Gr(\ga_\sg(t)), t\in
 [\frac{q\tau}{4},\frac{p\tau}{4}]).\lb{wang11}\ee
Since $D(t)\ge 0$ for $t\in [\frac{q\tau}{4},\frac{p\tau}{4}]$, we
have

\be D_\sg(t)>0,\quad \forall
t\in[\frac{q\tau}{4},\frac{p\tau}{4}].\lb{wang12}\ee

So by the proof of Lemma 3.1 of \cite{LZZ} and Lemma 2.6 of
\cite{LZZ}, we have
 \be \mu_F^{CLM}(L_0\times L_0,
\Gr(\ga_\sg(t)), t\in
 [\frac{q\tau}{4},\frac{p\tau}{4}])=\sum_{t\in [\frac{q\tau}{4},\frac{p\tau}{4})}\nu_{L_0}(\ga_\sg(t))\ge
 0.\lb{wang13}\ee
Thus by (\ref{wang11}) and (\ref{wang13}), (\ref{zh0}) holds.
 The proof of Theorem 3.1 is complete.\hfill\hb

By similar proof of Theorem 3.2 we have the following Theorem 3.3.

\noindent{\bf Theorem 3.3.} {If $\tau>0$ and $B\in
C([0,\frac{\tau}{4}],\mathcal{L}_s(\R^{2n}))$ satisfying condition
(B1) and $B\ge 0$, then for $j=0,1$ and any two positive integers
$p\ge q$ there holds \be i_{L_j}(\ga_B^{p})\ge
i_{L_j}(\ga_B^{q}).\ee}

\setcounter{equation}{0}
\section{Proof of Theorems 1.1-1.3 and Corollary 1.2}

In this section we study the minimal period problem for brake orbits
of the reversible Hamiltonian system (\ref{1.1}) and complete the
proof of Theorems 1.1-1.3 and Corollary 1.2.

For $T>0$, we set $E=W^{1/2,2}(S_T,\R^{2n})$ with the usual norm and
inner product denoted by $||\cdot||$ and $\langle\cdot\rangle$
respectively, and two subspaces of $E$ by $E_{T}=\{x\in
W^{1/2,2}(S_{\tau},\R^{2n})|\,x(-t)=Nx(t)\;a.e.\;t\in \R\}$ and
$\check{E}_{T}=\{x\in
W^{1/2,2}(S_{\tau},\R^{2n})|\,x(-t)=-Nx(t)\;a.e.\;t\in \R\}$. Then
we have
 \be E=E_T\oplus\check{E}_T.\ee

As in Section 3, we define two selfadjoint operators $A_{T}$ on
$E_{T}$ by the same way as (\ref{he8}). We also define two
selfadjoint operators $\check{A}_{T}$  on $\check{E}_{T}$ by the
following bilinear form: \bea \langle
\check{A}_{\T}x,y\rangle=\int_0^{T}-J\dot{x}\cdot y\,dt.\eea Then
$A_{T}$ is a bounded operator on $E_{T}$ and dim $\ker A_{T}=n$, the
Fredholm index of $A_{T}$ is zero, and $\check{A}_{T}$ is a bounded
operator on $\check{E}_{T}$ and dim $\ker \check{A}_{T}=n$, the
Fredholm index of $\check{A}_{T}$ is zero.

Set
       $$E_{T}(j)=\left\{z\in E_T\left|z(t)=\exp(\frac{2j\pi
       t}{T}J)a+\exp(-\frac{2j\pi t}{T}J)b,\; \forall t\in \R;\; \forall a,\,b\in L_0\right.\right\},$$
       $$E_{T, m}=E_{T}(0)+E_{T}(1)+\cdots +E_{T}(m)$$ and
 $$\check{E}_{T}(j)=\left\{z\in E_T\left|z(t)=\exp(\frac{2j\pi
       t}{T}J)a+\exp(-\frac{2j\pi t}{T}J)b,\; \forall t\in \R;\; \forall a,\,b\in L_1\right.\right\},$$
       $$\check{E}_{T, m}=\check{E}_{T}(0)+\check{E}_{T}(1)+\cdots +\check{E}_{T}(m).$$

Let $P_{T,m}$ be the orthogonal projection from $E_T$ to $E_{T, m}$
and $\check{P}_{T,m}$ be the orthogonal projection from
$\check{E}_T$ to $\check{E}_{T, m}$ for $m=0,1,2,...$, then
$\Ga_{T}=\{P_{T,m}: m=0,1,2,...\}$  and
$\check{\Ga}_{T}=\{\check{P}_{T,m}: m=0,1,2,...\}$ are the usual
Galerkin approximation schemes w.r.t. $A_{T}$ and $\check{A}_T$
respectively.

For $z\in {E_T}$, we define
  \be f(z)=\frac{1}{2}\langle A_Tz,z\rangle-\int_0^TH(z)dt.\ee

It is well known that $f\in C^2(E_T,\R)$ whenever,
    \be H\in C^2(\R^{2n})\qquad {\rm and \qquad} |H''(x)|\le
    a_1|x|^s+a_2\lb{n1}\ee
for some $s\in (1,+\infty)$ and all $x\in\R^{2n}$.

By similar argument of Lemma 4.1 of \cite{Z1}, looking for
$T$-periodic brake orbit solutions of (\ref{1.1}) is equivalent to
look for critical points of $f$.

In order to get the information about the Maslov-type indices, we
need the following theorem which was proved in \cite{Gh1,LaSo1,So1}.

\noindent{\bf Theorem 4.1.} {\it Let $W$ be a real Hilbert space
with orthogonal decomposition $E=X\oplus Y$, where $\dim X<+\infty$.
Suppose $f\in C^2(W,\R)$ satisfies (PS) condition and the following
conditions:

(i) There exist $\rho,\;\delta>0$ such that $f(w)\ge\delta$ for any
$w\in W$;

(ii) There exist $e\in \partial B_1(0)\cap Y$ and $r_0>\rho>0$ such
that for any $w\in\partial Q$, $f(w)<\delta$ where
$Q=(B_{r_0}(0)\cap X)\oplus\{re:0\le r\le r_0\}$, $B_r(0)=\{w\in W:
||w||\le r\}$.

Then (1) $f$ possesses a critical value $c\ge \delta$, which is
given by
     $$ c=\inf_{h\in \Ga}\max_{w\in Q}f(h(w)),$$
where $\Ga=\{h\in C(Q,E):h=id\; {\rm on}\; \partial Q\}$;

(2) There exists $w_0\in \mathcal{K}_c\equiv\{w\in E:\,f'(w)=0,\,
f(w)=c\}$ such that the Morse index $m^-(w_0)$ of $f$ at $w_0$
satisfies
      $$ m^-(w_0)\le \dim X+1.$$}

{\bf Proof of Theorem 1.3.} For any given $T>0$, we prove the
existence of $T$-periodic brake solution of (\ref{1.1}) whose
minimal period satisfies the inequalities in the conclusion of
Theorem 1.2. We divide the proof into five steps.

 {\bf Step 1.} We truncate the function $\hat{H}$ suitably and evenly such that it satisfies the growth condition (\ref{n1}).
 Hence corresponding new reversible function $H$ satisfies condition (\ref{n1}).

We follow the method in Rabinowitz's pioneering work \cite{Ra1} (cf.
also \cite{FKW1}, \cite{Ra2} and  \cite{Z1}). Let $K>0$ and $\chi\in
C^\infty(\R,\R)$ such that $\chi\equiv 1$ if $y\le K$, $\chi\equiv
0$ if $y\ge K$ and $\chi'(y)<0$ if $y\in(K,K+1)$, Where $K$ will be
determined later. Set
 \be \hat{H}_K(z)=\chi(|z|)\hat{H}(z)+(1-\chi(|z|))R_K|z|^4\ee

 and
 \be H_K(z)=\frac{1}{2}B_0x\cdot x+ \hat{H}_K(z),\ee
where the constant $R_K$ satisfies
        \be R_K\ge \max_{K\le |z|\le K+1} \frac{H(z)}{|z|^4}.\ee
 Then $H_K\in C^2(\R^{2n},\R)$. Since $\hat{H}$ satisfies (H3), $\forall \varepsilon>0$, there is a $\delta_1>0$
 such that $\hat{H}_K(z)\le\varepsilon|z|^2$ for $|z|\le \delta_1$. It is
 easy to see that $H_K(z)|z|^4$ is uniformly bounded as $|z|\to
 +\infty$, there is an $M_1=M_1(\varepsilon,K)$ such that $\hat{H}_K(z)\le M_1|z|^4$
 for $|z|\ge \delta_1$. So
  \be \hat{H}_K(z)\le \varepsilon |z|^2+M_1|z|^4,\quad \forall
  z\in\R^{2n}.\lb{he20}\ee

Set \bea f_K(z)=\frac{1}{2}\langle {A_T} z,z\rangle-\int_0^T
H_K(z)dt,\qquad \forall z\in\hat{E}.\nn\eea

Then $f_K\in C^2(E_T,\R)$ and \bea f_K(z)=\frac{1}{2}\langle
({A_T}-{B_0}_T) z,z\rangle-\int_0^T \hat{H}_K(z)dt,\qquad \forall
z\in\hat{E},\nn\eea where ${B_0}_T$ is the selfadjoint linear
compact operator on ${E_T}$ defined by \bea \langle
{B_0}_Tz,z\rangle=\int_0^TB_0z(t)\cdot z(t)\,dt.\nn\eea

  {\bf Step 2.} For $m>0$, let $f_{Km}=f|E_{T,m}$. We
  show $f_{Km}$ satisfies the hypotheses of Theorem 4.1.

 We set
 \bea &&X_m=M^-(P_{T,m}({A_T}-{B_0}_T)P_{T,m})\oplus M^0(P_{T,m}({A_T}-{B_0}_T)P_{T,m}),\nn\\
        &&  Y_m=M^+(P_{T,m}({A_T}-{B_0}_T)P_{T,m}).\nn\eea

For $z\in Y_m$, by (\ref{he20}), (\ref{feng}),  and the fact that
$P_{T,j}{B_0}_T=P_{T,j}{B_0}_T$ for $j>0$, we have
   \bea f_{Km}(z)&=&\frac{1}{2}\langle ({A_T}-{B_0}_T)
z,z\rangle-\int_0^T \hat{H}_K(z)dt\nn\\
             &\ge& \frac{1}{2}||({A_T}-{B_0}_T)^\#||^{-1}||z||^2-(\varepsilon
               ||z||_{L^2}^2+M_1||z||_{L^4}^4)\nn\\
               &\ge& \frac{1}{2}||({A_T}-{B_0}_T)^\#||^{-1}||z||^2-(\varepsilon
              C_2^2+M_1C_4^4||z||^2)||z||^2,\eea
where $C_2$ and $C_4$ are constants for $s=2,\,4$ for the Sobolev
embedding of inequality (\ref{feng}), and they are independent of
$m$ and $K$.

\noindent So if choose $\varepsilon>0$ small enough such that
$\varepsilon C_2^2< \frac{1}{4}||(A_T-{B_0}_T)^\#||^{-1}$, then
there exists $\rho=\rho(K)>0$ small enough and $\delta=\delta(K)>0$,
which are independent of $m$, such that
      \be f_m(z)\ge \delta,\qquad \forall z\in \partial
      B_\rho(0)\cap Y_m.\ee

Let $e\in B_1(0)\cap Y_m$ and set
   \bea Q_m=\{re:0\le r\le r_1\}\oplus (B_{r_1}(0)\cap X_m),\nn\eea
where $r_1$ will be determined later. Let $z=z_{-}+z_0\in
B_{r_1}(0)\cap X_m$, we have

\bea f_{Km}(z+re)&=&\frac{1}{2}\langle
({A_T}-{B_0}_T)z_,z_\rangle+\frac{1}{2}r^2\langle({A_T}-{B_0}_T)e,e\rangle-\int_0^T\hat{H}_K(z+re)dt\nn\\
   &\le& \frac{1}{2}||{A_T}-{B_0}_T||r^2-\frac{1}{2}||(A_T-{B_0}_T)^\#||^{-1}||z_{-}||^2-\int_0^T\hat{H}_K(z+re)dt.\lb{f1}\eea

Since $\hat{H}$ satisfies (H2) we have
    \bea \hat{H}_K(x)\ge a_1|x|^\alpha-a_2,\qquad\forall
    x\in\R^{2n},\nn\eea
where $\alpha=\min\{\mu,4\}$, $a_1>0$, $a_2$ are two constants
independent of $K$ and $m$. Then there holds

\be \int_0^T \hat{H}_K(z+re)dt\ge a_1\int_0^T|z+re|^\alpha-Ta_2\ge
a_3(||z_0||_{L^\alpha}^\alpha+r^\alpha)-a_4,\lb{f2}\ee where $a_3$
and $a_4$ are constants independent of $K$ and $m$. By (\ref{f1})
and (\ref{f2}) we have \bea f_{Km}(z+re)\le
\frac{1}{2}||{A_T}-\hat{B_0}||r^2-\frac{1}{2}||(A-B_0)^\#||^{-1}||z_{-}||^2-a_3(||z_0||_{L^\alpha}^\alpha+r^\alpha)+a_4.\nn\eea
Since $\alpha>2$ there exists a constant $r_1>\rho>0$, which are
independent of $K$ and $m$, such that \be f_{Km}\le0,\qquad\forall
z\in \partial Q_{m}.\ee
 Then by Theorem 4.1, $f_{Km}$ has a critical value $c_{Km}$,
 which is given by
 \be c_{Km}=\inf_{g\in\Ga_m} \max_{z\in Q_m} f_{Km}(g(z)),\ee
 where $\Ga_m=\{g\in C(Q_m,\hat{E}_m|g=id; {\rm on}\; \partial
 Q_m\}$. Moreover there is a critical point $x_{Km}$ of
 $f_{Km}$ which satisfies
   \bea m^-(x_{Km})\le \dim X_m+1.\lb{y1}\eea

{\bf Step 3.} We prove that there exists a $T$-periodic brake orbit
solution $x_T$ of (\ref{1.1}) which satisfies $i_{L_0}(x_T)\le
i_{L_0}(B_0)+\nu_{L_0}(B_0)+1$.

 Note that $id\in\Ga_m$, by (\ref{f1}) and condition
(H4), we have

\bea c_{Km}\le \sup_{z\in
Q_m}f_{Km}(z)\le\frac{1}{2}||{A_T}-{B_0}_T||r_1^2.\nn\eea

Then $\{c_{Km}\}$ possesses a convergent subsequence, we still
denote it by $\{c_{Km}\}$ for convenience. So there is a
$c_K\in[\delta,] $ such that $c_{Km}\to c_K$.

By the same arguments as in section 6 of \cite{Ra2} we have $f_K$
satisfies $(PS)_c^*$ condition for $c\in\R$, i.e., any sequence
${z_m}$ such that $z_m\in E_{T,m}$, $f_{Km}'(z_m)\to 0$ and
$f_{Km}(z_m)\to c$ possesses a convergent subsequence in $E_T$.
Hence in the sense of subsequence we have
 \bea x_{Km}\to x_K,\qquad f_K(x_K)=c_K,\qquad
 f'_K(x_K)=0.\lb{y2}\eea
By similar argument in \cite{Ra2}, $x_K$ is a classical nonconstant
symmetric $T$-periodic solution of

 \bea \dot{x}=JH_K'(x), \quad x\in\R^{2n}.\nn\eea

 Set $B_K(t)=H''_K(x_K(t))$, Then $B_K\in
C([0,T/2],\mathcal{L}_s(\R^{2n}))$ and satisfies condition (B1). Let
${B_K}_T$ be the operator defined by the same way of the definition
of ${B_0}_T$. It is easy to show that \bea
||f''(z)-({A_T}-{B_K}_T)||\to 0\qquad {\rm as}\;\; ||z-x_K||\to
0.\nn\eea

So for $0<d\le\frac{1}{4}||(A_T-B_{K_T})^\#||^{-1}$, there exists
$r_2>0$ such that
      \bea ||f_{Km}''(z)-P_{T,m}({A_T}-{B_K}_T)P_{T,m}||\le ||f''(z)-({A_T}-{B_K}_T)||\le\frac{1}{2}d,
      \;
      \forall z\in\{z\in E_T:||z-x_K||\le r_2\}.\nn\eea

Then for $z\in \{z\in E_T: ||z-x_K||\le r_2\}\cap E_{T,m}$, $\forall
u\in M^-_d(P_{T,m}({A_T}-{B_K}_T)P_{T,m})\setminus\{0\}$, we have
   \bea \langle f_{Km}''(z)u,u\rangle &\le& \langle P_{T,m}({A_T}-{B_K}_T)P_{T,m}u,u\rangle
   +\|f_{Km}''(z)-P_{T,m}({A_T}-{B_K}_T)P_{T,m}\|\|u\|^2\nn\\
       &\le& -\frac{1}{2}d\|u\|^2.\nn\eea
So we have
   \be m^-(f_{Km}''(z))\ge \dim
   M^-_d(P_{T,m}({A_T}-{B_K}_T)P_{T,m}).\lb{y3} \ee

By Theorem 2.1 of \cite{Liu2} and Remark 3.1,  there is $m^*>0$ such
that for $m\ge m^*$ we have
      \bea &&\dim
      X_m=mn+n+i_{L_0}(B_0)+\nu_{L_0}(B_0),\lb{y5}\\
      &&\dim
      M^-_d(P_{T,m}({A_T}-{B_K}_T)P_{T,m})=mn+n+i_{L_0}(B_K).\lb{y6}\eea
Then by (\ref{y1}), (\ref{y2}), and (\ref{y3})-(\ref{y6}), we have
       \bea i_{L_0}(B_K)\le
       i_{L_0}(B_0)+\nu_{L_0}(B_0)+1.\nn\eea

By the similar argument as in the section 6 of \cite{Ra2}, there is
a constant $M_2$ independent of $K$ such that $||x_K||_\infty\le
M_2$. Choose $K>M_2$. Then $x_K$ is a non-constant symmetric
$T$-periodic solution of the problem (\ref{1.1}). From now on in the
proof of Theorem 1.3, we write $B=B_K$ and $x_T=x_K$. Then $x_T$ is
a non-constant symmetric $T$-periodic solution of the problem
(\ref{1.1}), and $B$ satisfies
    \bea i_{L_0}(x_T)=i_{L_0}(B)\le
       i_{L_0}(B_0)+\nu_{L_0}(B_0)+1.\lb{y7}\eea

      Since $x_T$ obtained in Step 3 is a nonconstant and symmetric
      $T$-period solution, its minimal period $\tau=\frac{T}{k}$
      for some $k\in\N$.

  We denote by $x_\tau=x_T|_{[0,\tau]}$, then it is a brake orbit solution of (\ref{1.1}) with the minimal $\tau$ and
  $X_T=x_\tau^{k}$ being the $k$ times iteration of $x_\tau$.
  As in Section 1, let $\ga_{x_T}$ and $\ga_{x_\tau}$ be the symplectic path
  associated to $(\tau, x)$ and $(T,x_T)$ respectively.
  Then $\ga_{x_\tau}\in C([0,\frac{\tau}{2}],\Sp(2n))$ and
 $\ga_{x_T}\in C([0,\frac{T}{2}],\Sp(2n))$. Also we have
 $\ga_{x_T}=\ga_{x_\tau}^{k}$.

{\bf Step 4.} We prove that
 \bea i_{L_1}(\ga_{x_\tau})+\nu_{L_1}(\ga_{x_\tau})\ge 1.\nn\eea

We follow the way of the proof of Theorem 1.2 of \cite{FKW1}. By the
same way as $\check{E}_T$ and $\check{A}_T$ we can define the space
$\check{E}_\tau$ and the operator $\check{A}_\tau$ on it. Also we
can define the orthogonal projection $\check{P}_\tau,m$ and the
subspaces $\check{E}_{\tau,m}$ for $m=0,1,2,...$.  Let
$\check{B}_\tau$ be the selfadjoint linear compact operator on
${\check{E}_T}$ defined by:
 \bea \langle \check{B}_\tau
z,z\rangle=\int_0^\tau B(t)z(t)\cdot z(t)\,dt,\quad \forall
z\in\check{E}_\tau. \nn\eea

For $z\in \check{E}_\tau$, set
 \bea f_\tau(z)=\frac{1}{2}\langle(\check
{A}_\tau-\check{B}_\tau)z,z\rangle=\frac{1}{2}\langle\check {A}_\tau
z,z\rangle-\frac{1}{2}\int_0^\tau H''(x_\tau(t))z\cdot z\,dt\nn\eea
 and
  \bea f_{\tau m}(w)=f_\tau(w),\qquad \forall w\in
  \check{E}_{\tau,m}.\nn\eea

Let
$$X=\{z\in L_1|\, B_0z=0 \;{\rm and}\;
\hat{H}''(x_\tau(t))z=0,\;\forall t\in\R\}$$
 and $Y$ be the orthogonal complement of $X$ in $L_1$, i.e.,
 $L_1=X\oplus Y$. Since $H''(x_\tau(t))=B_0+\hat{H}''(x_\tau(t))$,
 by (H4) it is easy to see that there exists $\lm_0>0$ such that
  \bea \int_0^\tau H''(x_\tau(t))z_0\cdot z_0\,dt\ge
  \lm_0||z_0||,\qquad \forall z_0\in Y.\nn\eea

Thus for any $z=z_-+z_0\in
\check{P}_{\tau,m}M^-(\check{A}_\tau)\oplus Y$ with $||z||=1$, we
have
 \bea f_{\tau m}(z)&=&\frac{1}{2}\langle(\check
{A}_\tau-\check{B}_\tau)z,z\rangle=\frac{1}{2}\langle\check {A}_\tau
z_-,z_-\rangle-\frac{1}{2}\int_0^\tau H''(x_\tau(t))z\cdot z\, dt\lb{q1}\\
&\le &
-\frac{1}{2}||\check{A}_\tau^\#||^{-1}||z_-||^2-\frac{1}{2}\int_0^\tau
H''(x_\tau(t))z_0\cdot z_0\,dt-\int_0^\tau H''(x_\tau(t))z_-\cdot
z_0\,dt\nn\\
&\le &
-\frac{1}{2}||\check{A}_\tau^\#||^{-1}||z_-||^2-\frac{\lm_0}{2}||z_0||^2
+\max_{t\in[0,\tau]}||H''(x_\tau(t))||\,||z_-||\,||z_0||.\eea
 Since
 \bea ||z_-||\,||z_0||\le \frac{\var}{4}||z_-||^2+\frac{1}{\var} ||z_0||^2, \quad\forall \var>0
 .\nn\eea
By choosing $\var$ suitably one can see that there exists $0<c_0<1$
with $|1-c_0|$ small enough such that if $||z_0||\le c_0$,
 \be f_{\tau m}(z)\le -\frac{\lm_0}{4}c_0^2.\ee

When $||z_0||\le c_0$, we have $||z_-||^2\ge 1-c_0^2$. By (\ref{q1})
and (H4)
 \bea f_{\tau m}(z)\le
 -\frac{1}{2}||\check{A}_\tau^\#||^{-1}||z_-||^2\le
 -\frac{1}{2}||\check{A}_\tau^\#||^{-1}(1-c_0^2).\nn\eea
Hence we always have
    \be f_{\tau m}(z)\le -c||z||^2,\qquad \forall z\in
\check{P}_{\tau,m}M^-(\check{A}_\tau)\oplus Y,\lb{q2}\ee where
$c=\max\{\frac{\lm_0}{4}c_0^2,
\frac{1}{2}||\check{A}_\tau^\#||^{-1}(1-c_0^2)\}$ is independent of
$m$. Let $$d=\min\{
\frac{1}{4}||(\check{A}_\tau-\check{B}_\tau)^\#||^{-1},\frac{c}{2}\}.$$
By (\ref{q2}) and Theorem 2.1 of \cite{Liu2} and Remark 3.1 and the
definition of $i_{L_1}(\ga(x_\tau))$, for $m$ large enough,  we have
 \bea mn+n+i_{L_1}(\ga(x_\tau))&=&\dim M^-_d(\check{P}_{\tau,m}
 (\check{A}_\tau-\check{B}_{\tau})\check{P}_{\tau,m})\nn\\
 &\ge& \dim (\check{P}_{\tau,m}M^-(\check{A}_\tau)\oplus Y)\nn\\
 &=& mn + n-\dim X,\eea
which implies that
 \be i_{L^1}(\ga(x_\tau))\ge -\dim X.\lb{q3}\ee
Since $x_\tau$ is a nonconstant brake solution of (\ref{1.1}), by
the definition of $X$ we have
 \be \nu_{L^1}(\ga(x_\tau))\ge \dim X +1.\lb{q4}\ee
Hence by (\ref{q3}) and (\ref{q4}) we have
 \be i_{L^1}(\ga(x_\tau))+\nu_{L^1}(\ga(x_\tau))\ge 1.\lb{q5}\ee

{\bf Step 5.} Finish the proof of Theorem 1.3.

By Theorem 2.1 and Theorem 6.2 below (also Theorem 2.6 of
\cite{Liu3}) we have
 \bea &&i_{L_0}(\ga_{x_\tau}^k)\ge i_{L_0}(\ga_{{x_\tau}})+\frac{k-1}{2}
 (i_1(\ga^2)+ \nu_1(\ga^2)-n),\quad {\rm if}\, k\in 2\N-1,\lb{q8}\\
&&i_{L_0}(\ga_{x_\tau}^k)\ge
i_{L_0}(\ga_{{x_\tau}})+i_{\sqrt{-1}}^{L_0}(\ga_{{x_\tau}})+(\frac{k}{2}-1)
 (i_1(\ga^2)+ \nu_1(\ga^2)-n),\quad {\rm if}\, k\in 2\N.\lb{q9}\eea

Since $B_0$ is semipositive and  $\hat{H}$ satisfies (H4), by
Corollary 3.2, we have
   \be i_{L_0}(\ga_{x_\tau})+\nu_{L_0}(\ga_{x_\tau})\ge 0.\lb{q6}\ee
 By Proposition C of \cite{LZZ} and the definitions of $i_{L_0}$ and $i_{L_1}$
 we have
 \bea i_1(\ga^2)=i_{L_0}(\ga)+i_{L_1}(\ga)+n,\nn\\
     \nu_1(\ga^2)=\nu_{L_0}(\ga)+\nu_{L_1}(\ga).\nn\eea

 So by (\ref{q5}) and (\ref{q6}) we have
  \be i_1(\ga^2)+ \nu_1(\ga^2)-n\ge 1.\lb{q7}\ee

So by (\ref{q8}), (\ref{q9}) and (\ref{q7}) we have

\bea &&i_{L_0}(\ga_{x_\tau}^k)\ge
i_{L_0}(\ga_{{x_\tau}})+\frac{k-1}{2},\quad {\rm if}\, k\in 2\N-1,\lb{q10}\\
&&i_{L_0}(\ga_{x_\tau}^k)\ge
i_{L_0}(\ga_{{x_\tau}})+\frac{k-1}{2},\quad {\rm if}\, k\in
2\N.\lb{q11}\eea

By (\ref{y7}) and the definition of $\ga_{x_\tau}$ we have
 \be i_{L_0}(\ga_{{x_\tau}})^k)\le i_{L_0}(B_0)+\nu_{L_0}(B_0)+1.\lb{q15}\ee

By Corollary 3.2, we have
 \be i_{L_0}(\ga_{{x_\tau}})\ge -n.\lb{q14}\ee

So by (\ref{q10})-(\ref{q14}) we have

\be k\le 2(i_{L_0}(B_0)+\nu_{L_0}(B_0))+2n+4.\ee

{\bf Claim 2.} $k$ can not be $2(i_{L_0}(B_0)+\nu_{L_0}(B_0))+2n+3$
and $2(i_{L_0}(B_0)+\nu_{L_0}(B_0))+2n+4$.

Hence by Claim 2, $k\le 2(i_{L_0}(B_0)+\nu_{L_0}(B_0))+2n+2$, and
Theorem 1.3 holds.

{\bf Proof of Claim 2.} We first show that $k$ can not be
$2(i_{L_0}(B_0)+\nu_{L_0}(B_0))+2n+3$. Otherwise, we have
 \be k=2(i_{L_0}(B_0)+\nu_{L_0}(B_0))+2n+3.\lb{q12}\ee
The equality in (\ref{q8}) holds, then by (\ref{q7}), in this case
there must hold that
 \bea i_1(\ga^2)+ \nu_1(\ga^2)-n=1 \lb{q18}\eea
  and
  \bea i_{L_0}(\ga_{{x_\tau}})=-n.\lb{q16}\eea
By Corollary 3.2 again we have that
   \be \nu_{L_0}(\ga_{{x_\tau}})=n.\lb{q17}\ee
Also by (\ref{q18}) we have
 \bea i_{L^1}(\ga(x_\tau)+\nu_{L^1}(\ga(x_\tau))=1.\nn\eea
Denote by $\nu_{L^1}(\ga(x_\tau))=r$. Then we have
 \bea &&i_{L^1}(\ga(x_\tau))=1-r,\lb{q19}\\
     &&\nu_{L^1}(\ga(x_\tau))\ge 1.\lb{q20}\eea

By (\ref{q16}) and (\ref{q19}) we have
 \be i_{L_0}(\ga_{{x_\tau}})-i_{L_1}(\ga_{{x_\tau}})=r-n-1.\lb{q40}\ee

So we can write $\ga_{{x_\tau}}(\frac{\tau}{2})=
\left(\begin{array}{cc}A&0\\C&D\end{array}\right)$ with $A,C,D$ to
be $n\times n$ real matrices. Hence by (4.2) of \cite{LZ} we have
 \bea \ga_{{x_\tau}}^2(\tau)=N\ga_{{x_\tau}}(\frac{\tau}{2})^{-1}N\ga_{{x_\tau}}(\frac{\tau}{2})
 = \left(\begin{array}{cc}D^TA&0\\C^TA&A^TD\end{array}\right).\nn\eea
 Since $\ga_{x_\tau
 }(\frac{\tau}{2})$ is a symplectic matrix we have
  \bea A^TD=D^TA=I_n, \quad C^TA=A^TC.\nn\eea
So we have
 \bea
 \ga_{{x_\tau}}^2(\tau)=\left(\begin{array}{cc}I_n&0\\C^TA&I_n\end{array}\right).\nn\eea
Note that here $C^TA$ is a symmetric matrix and $A$ is invertible.
So by (\ref{q20}) there exists a orthogonal matrix $Q$ such that
 \be
 Q(C^TA)Q^T=\diag(0,0,...,0,\lm_1,\lm_2,...,\lm_p,\lm_{p+1},...,\lm_{n-p-r})\lb{qq}\ee with
$\lm_j>0$ for $j=1,2,...,p$ and $\lm_j<0$, for
$j=p+1,p+2,...,n-p-r$, where $1\le p\le n-r$. Then it is easy to
check that $(I_2)^{\diamond^r}\diamond
N_1(1,-1)^{\diamond^p}\diamond N_1(1,1)^{\diamond^{(n-p-r)}}\in
\Om^0(\ga_{{x_\tau}})$ with $\Om^0(\ga_{x\tau })$ to be defined in
Section 6 below. Then by Theorem 6.2 below or Theorem 2.6 of
\cite{Liu3}, when the equality in (\ref{q8}) holds, there must hold
$p=n-r$. Hence we have \bea
 &&Q(C^TA)Q^T=\diag(0,0,...,0,\lm_1,\lm_2,...,\lm_{n-r})\lb{q24},\\
&&\lm_j>0,\qquad {\rm for}\; j=1,2,...,n-r.\lb{q23}\eea

{\bf Case 1.} If $\det A>0$, then there exists a invertible matrix
path $\rho(s)$ for $s\in [0,\frac{\tau}{2}]$ connecting it and $I_n$
such that $\rho(0)=I_n$ and $\rho(1)=A$.

We define a symplectic path $\phi_1$ by
   \bea \phi_1(s)=\left(\begin{array}{cc}\rho(s)^{-1}&0\\0&\rho(s)^T\end{array}\right)
   \left(\begin{array}{cc}A&0\\C&D\end{array}\right), \quad \forall
   s\in[0,\frac{\tau}{2}].\nn\eea
Then $\nu_{L_j}(\phi_1(s)=constant$ for $j=0,1$ and
$s\in[0,\frac{\tau}{2}]$. So by Definition 2.5 and Lemma 2.8 and
Proposition 2.11 of \cite{LZZ}, for $j=1,2$ we have
 \be \mu_F^{CLM}(V_j,\Gr(\phi_1),[0,\frac{\tau}{2}])=0.\lb{q30}\ee
Also we have
$\phi_1(0)=\left(\begin{array}{cc}A&0\\C&D\end{array}\right)$ and
$\phi_1(0)=\left(\begin{array}{cc}I_n&0\\A^TC&I_n\end{array}\right)$.

Note that we can always choose the orthogonal matrix $Q$ in
(\ref{q24}) such that $\det Q=1$ (otherwise we replace it by
$\diag(-1,1,...,1)Q$). Then there exists a invertible matrix path
$\rho_2(s)$ for $s\in [0,\frac{\tau}{2}]$ connecting it and $I_n$
such that $\rho_2(0)=I_n$ and $\rho_2(\frac{\tau}{2})=Q$. We define
a symplectic path $\phi_2$ by
   \bea \phi_2(s)=
   \left(\begin{array}{cc}I_n&0\\\rho_2(s)A^TC\rho_2(s)^T&I_n\end{array}\right), \quad \forall
   s\in[0,\frac{\tau}{2}].\nn\eea
Then $\nu_{L_j}(\phi_2(s)=constant$ and for $j=0,1$ and
$s\in[0,\frac{\tau}{2}]$. So by Definition 2.5 and Lemma 2.8 and
Proposition 2.11 of \cite{LZZ} again, for $j=1,2$ we have
 \be \mu_F^{CLM}(V_j,\Gr(\phi_2),[0,\frac{\tau}{2}])=0.\lb{q31}\ee
 Also we have
 \bea&&\phi_2(0)=\left(\begin{array}{cc}I_n&0\\A^TC&I_n\end{array}\right)\nn\\
 &&\phi_2(\frac{\tau}{2})=\left(\begin{array}{cc}I_n&0\\QA^TCQ^T&I_n\end{array}\right)=(I_2)^{\diamond^r}\diamond
N_1(1,\lm_1)\diamond\cdots\diamond N_1(1,\lm_{n-r}).\lb{q33}\eea

By the Reparametrization invariance and Path additivity of the
Maslov index $\mu_F^{CLM}$ in \cite{CLM} and (\ref{q30}) and
(\ref{q31}), for $j=1,2$ we have
 \bea \mu_F^{CLM}(V_j,\Gr(\ga_{x_\tau}),[0,\frac{\tau}{2}])=
 \mu_F^{CLM}(V_j,\Gr(\phi_2*(\phi_1*\ga_{{x_\tau}})),[0,\frac{\tau}{2}]),\nn\eea
where the joint path $\phi_2*(\phi_1*\ga_{{x_\tau}})$ is defined by
(\ref{zzz}).
 So by definition for $j=0,1$ we have
  \be i_{L_j}(\ga_{{x_\tau}})=i_{L_j}(\phi_2*(\phi_1*\ga_{{x_\tau}})).\lb{q32}\ee
Then by Theorem 2.3 and (\ref{q33}) we have
  \bea i_{L_0}(\ga_{{x_\tau}})-i_{L_1}(\ga_{{x_\tau}})=\frac{1}{2}\sgn M_\var((I_2)^{\diamond^r}\diamond
N_1(1,\lm_1)^T\diamond\cdots\diamond N_1(1,\lm_{n-r})^T).\eea
 By Remark
2.1 and the computations (\ref{guo1})-(\ref{guo3}) at the end of
Section 2, for $\var>0$ small enough we have
 \be \sgn
M_\var((I_2)^{\diamond^r}\diamond
N_1(1,\lm_1)^T\diamond\cdots\diamond N_1(1,\lm_{n-r})^T)=2(r-n).\ee
 So we have
  \be i_{L_0}(\ga_{{x_\tau}})-i_{L_1}(\ga_{{x_\tau}})=r-n,\ee
 which contradicts to (\ref{q40}).

{\bf Case 2.} If $\det A<0$, then there exists a invertible matrix
path $\rho(s)$ for $s\in [0,\frac{\tau}{2}]$ such that
$\rho(0)=\diag(-1,1,1,...,1)$ and $\rho(1)=A$. by similar arguments
we can show that
     \bea i_{L_0}(\ga_{{x_\tau}})-i_{L_1}(\ga_{{x_\tau}})=\frac{1}{2}\sgn
     M_\var((-I_2)\diamond (I_2)^{\diamond^(r-1)}\diamond
N_1(1,\lm_1)\diamond\cdots\diamond N_1(1,\lm_{n-r}))=r-n,\eea which
still contradicts to (\ref{q40}).

Hence we have proved that $k$ can not be
$2(i_{L_0}(B_0)+\nu_{L_0}(B_0))+2n+3$. By the same argument we can
prove that $k$ can not be $2(i_{L_0}(B_0)+\nu_{L_0}(B_0))+2n+4$.
Thus Claim 2 is proved and the proof of Theorem 1.3 is
complete.\hfill\hb

{\bf Proof of Theorem 1.1.} Note that this is the case $B_0=0$ of
Theorem 1.3. Then by Theorem 1.3 and the fact that $i_{L_0}(0)=-n$
and $\nu_{L_0}(0)=n$, the minimal period of $x_T$ is no less than
$\frac{T}{2n+2}$. In the following we prove that if (\ref{zhi6})
holds then the minimal period of $x_T$ belongs to
$\{T,\frac{T}{2}\}$.

Let $x_T$ is the $k$-time iteration of $x_\tau$ with $\tau$ being
the minimal period of $x_\tau$ and $\tau=\frac{T}{k}$. Then by the
proof of Theorem 1.3 with $B_0=0$ we have (\ref{q5}), (\ref{q8}) and
(\ref{q9}) hold. Since (\ref{zhi6}) holds, by Lemma 3.3 we have
 \be i_{L_0}(\ga_{x\tau})\ge 0.\lb{zhi7}\ee
So by (\ref{q8}) if $k$ is odd, we have
 \be 1\ge 0+\frac{k-1}{2}.\ee
 Hence $k\le 3$. Now we prove that $k$ can not be $3$, other wise
  we have
   \bea &&i_{L_0}(\ga_{x_\tau})=0,\lb{zhi9}\\
     && \nu_{L_0}(\ga_{x_\tau})=0,\lb{zhi10}\\
  &&i_{L_1}(\ga_{x_\tau})+\nu_{L_1}(\ga_{x_\tau})=1.\lb{zhi11}\eea
And by Theorem 2.1 and Theorem 6.2 we have
 \bea 1\ge
 i_{L_0}(\ga_{x_\tau}^3)=i_{L_0}(\ga_{x_\tau})+i_{e^{2\pi/3}}(\ga_{x_\tau}^2)\ge
 (i_1(\ga_{x_\tau}^2)-\nu_1(\ga_{x_\tau}^2)-n)\ge 1.\lb{zhi8}\eea
Then all the equalities of (\ref{zhi8}) hold. By Lemma 6.2 and 2 of
Theorem 6.2 again, there exist $p\ge 0$, $q\ge 0$ with $p+q\le n$
and $0<\theta_1\le \theta_2\le...\le \theta_{n-(p+q)}\le 2\pi/3$
such that
       \be (I_2)^{\diamond^p}\diamond
       N_1(1,-1)^{\diamond^q}\diamond R(\theta_1)\diamond
       R(\theta_2)\diamond...\diamond R(\theta_{n-p-q})\in \Om^0((\ga^2_{x_\tau})(\tau)),\lb{zhi12}\ee
where $\Om^0(M)$ for a symplectic matrix $M$ is defined in Section
6. By (\ref{zhi12}) we have
 \be -1\notin \sg((\ga^2_{x_\tau})(\tau)).\lb{zhi13}\ee

Now we denote by
$\ga_{x_\tau}(\frac{\tau}{2})=\left(\begin{array}{cc}A&B\\C&D\end{array}\right)$
with $A,B,C,D$ are all $n\times n$ matrices.

{\bf Claim 1.} Both $D$ and $A$ are invertible.

     We first prove $D$ is invertible. Otherwise, there exists a
     $n\times n$ invertible matrix $P$ such that
     $P^{-1}DP=\left(\begin{array}{cc}0&0\\0&R\end{array}\right)$
 and $R$ is a $(n-r)\times (n-r)$ matrix with $r\ge 1$.
 So we have
      \bea \left(\begin{array}{cc}P^T&0\\0&P^{-1}\end{array}\right)
          \left(\begin{array}{cc}A&B\\C&D\end{array}\right)
           \left(\begin{array}{cc}(P^{-1})^T&0\\0&P\end{array}\right)
           :=\left(\begin{array}{cc}\td{A}&\td{B}\\\td{C}&\td{D}\end{array}\right)\nn\eea
with $\td{D}=\left(\begin{array}{cc}0&0\\0&R\end{array}\right)$.
Since
$\left(\begin{array}{cc}\td{A}&\td{B}\\\td{C}&\td{D}\end{array}\right)$
is a symplectic matrix, we have
 \be \td{A}^TD-\td{C}^T\td{B}=I_n.\lb{zhi14}\ee
Since  $\td{D}=\left(\begin{array}{cc}0&0\\0&R\end{array}\right)$,
$\td{B}^T\td{D}$ and $\td{A}^T\td{D}$ both have form
 \bea \td{B}^T\td{D}=
 \left(\begin{array}{cc}0&*\\0&*\end{array}\right),\;\;\td{A}^T\td{D}
 =\left(\begin{array}{cc}0&*\\0&*\end{array}\right).\lb{zhi15}\eea
So by (\ref{zhi14}) and (\ref{zhi15}) we have
 \be
 \td{A}^T\td{D}+\td{C}^T\td{B}=2\td{A}^T\td{D}-I_n=\left(\begin{array}{cc}-I_r&*\\0&*\end{array}\right).\lb{zhi16}\ee
 By direct
computation and (\ref{zhi15}) and (\ref{zhi16}) we have
 \bea
 && N\left(\begin{array}{cc}\td{A}&\td{B}\\\td{C}&\td{D}\end{array}\right)^{-1}N
  \left(\begin{array}{cc}\td{A}&\td{B}\\\td{C}&\td{D}\end{array}\right)\nn\\
 &= &\left(\begin{array}{cc}P^T&0\\0&P^{-1}\end{array}\right)N\left(\begin{array}{cc}A&B\\C&D\end{array}\right)^{-1}
   N\left(\begin{array}{cc}A&B\\C&D\end{array}\right)\left(\begin{array}{cc}(P^{-1})^T&0\\0&P\end{array}\right)\lb{zhi17}\\
  &=&
  \left(\begin{array}{cc}\td{D}^T\td{A}+\td{B}^T\td{C}&2\td{B}^T\td{D}\\2\td{A}^T\td{C}&\td{A}^T\td{D}+\td{C}^T\td{B}\end{array}\right)\nn\\
  &=& \left(\begin{array}{cccc}*&*&0&*\\**&*&0&*\\ ** &*&-I_r&*\\ **
  &*&0&*\end{array}\right)\lb{zhi18}
\eea

Since by (4.2) o f{\cite{LZ} we have
 \bea \ga^2_{x_\tau}(\tau)=N\left(\begin{array}{cc}A&B\\C&D\end{array}\right)^{-1}
   N\left(\begin{array}{cc}A&B\\C&D\end{array}\right),\nn\eea
by (\ref{zhi17}) and (\ref{zhi18}) we have
 \bea -1\in \sg(\ga^2_{x_\tau}(\tau)),\nn\eea
 which contradicts to (\ref{zhi13}). Thus we have proved that $D$ is
 invertible. Similarly we can prove $A$ is invertible, and Claim 1
 is proved.

{\bf Claim 2.} There exists a invertible $n\times n$ real matrix $Q$
with $\det Q>0$ such that \be Q^{-1}(B^TC)Q=\diag
({0,0,...,0,\lm_1,\lm_2,...\lm_{n-r}})\ee
 with
$r=\nu_{L_1}(\ga_{x_\tau})$ and $\lm_i\in (-1,0)$ for
$i=1,2,...,n-r$.

In fact \bea
\ga^2_{x_\tau}(\tau)&=&N\left(\begin{array}{cc}A&B\\C&D\end{array}\right)^{-1}
   N\left(\begin{array}{cc}A&B\\C&D\end{array}\right)\nn\\
   &=&\left(\begin{array}{cc}D^T&B^T\\C^T&A^T\end{array}\right)\left(\begin{array}{cc}A&B\\C&D\end{array}\right)\nn\\
   &=&\left(\begin{array}{cc}I+2B^TC&2B^TD\\2A^TC&I+2C^TB\end{array}\right).\lb{zhi19}\eea

Since $B$ and $D$ are both invertible, for any $\om\in\C$, we have
  \bea &&\left(\begin{array}{cc}I_n&0\\-\frac{1}{2}(I_n+2C^TB-\om
  I_n)D^{-1}(B^T)^{-1}&I_n\end{array}\right)\left(\begin{array}{cc}I+2B^TC-\om I_n&2B^TD\\2A^TC&I+2C^TB-\om
       I_n\end{array}\right)\nn\\
       &=&\left(\begin{array}{cc}I+2B^TC-\om I_n&2B^TD\\-\frac{1}{2}(I_n+2C^TB-\om
  I_n)D^{-1}(B^T)^{-1}(I+2B^TC-\om
  I_n)+2A^TC&0\end{array}\right).\nn\eea
So we have
     \bea \det(\ga^2_{x_\tau}(\tau)-\om I_{2n})&=&\det(B^TD)\det((I_n+2C^TB-\om
  I_n)D^{-1}(B^T)^{-1}(I+2B^TC-\om I_n)-4A^TC)\nn\\
  &=&\det (D[(I_n+2C^TB-\om
  I_n)D^{-1}(B^T)^{-1}(I+2B^TC-\om I_n)-4A^TC]B^T)\nn\\
  &=&\det (D[I_n+2C^TB-\om
  I_n)D^{-1}(B^T)^{-1}(I+2B^TC-\om I_n]B^T-4DA^TCB^T)\nn\\
  &=& \det ((I+2B^TC-\om I_n)^2-4(1+CB^T)CB^T)\nn\\
  &=&\det (\om^2 I_n-2\om(I+2CB^T)+I).\lb{zhi21}\eea
By (\ref{zhi12}) we have
 \be \sg(\ga^2_{x_\tau}(\tau))\subset \U.\ee
So for $\om\in \U$ by (\ref{zhi21}) we have
 \be \det(\ga^2_{x_\tau}(\tau)-\om
 I_{2n})=(-4)^n\om^n\det(CB^T-\frac{1}{2}({\rm Re}\, \om -1)).\ee
Hence by (\ref{zhi12}) again we have $\sg(CB^T)\subset (-1,0]$,
moreover there exists a invertible $n\times n$ matrix $S$ such that
 \be S^{-1}CB^TS=\diag (0,0,...,0,\lm_1,\lm_2,...,\lm_{n-r}).\ee
with $r=\nu_{L_1}(\ga_{x_\tau})$ and $\lm_i\in (-1,0)$ for
$i=1,2,...,n-r$. Since $S^{-1}CB^TS=  (B^TS)^{-1}B^TC(B^TS)$, let
$Q=B^TS$, if $\det Q<0$ we replace it by $B^TS\diag (-1,1,1,...,1)$,
Claim 2 is proved.

{\bf Continue the proof of Theorem 1.1.}

If $\det B>0$, there is a continuous symplectic matrix path joint
$\left(\begin{array}{cc}B^{-1}&0\\0&B^T\end{array}\right)$ and
$I_{2n}$. Since
 \bea \left(\begin{array}{cc}B^{-1}&0\\0&B^T\end{array}\right)\left(\begin{array}{cc}A&B\\C&D\end{array}\right)
 =\left(\begin{array}{cc}B^{-1}A&I_n\\B^TC&B^TD\end{array}\right).\nn\eea
 By Lemma 2.2, for $\var>0$ small enough, we have
  \be \sgn M_\var \left(\left(\begin{array}{cc}A&B\\C&D\end{array}\right)\right)=
       \sgn M_\var
       \left(\left(\begin{array}{cc}B^{-1}A&I_n\\B^TC&B^TD\end{array}\right)\right). \lb{zhi30}\ee
 If $\det B<0$, there is a continuous symplectic path joint $\left(\begin{array}{cc}B^{-1}&0\\0&B^T\end{array}\right)$
 and $(-I_2)\diamond I_{2(n-1)}$. By direct computation we have
    \bea \sgn M_\var
    \left(\left(\begin{array}{cc}A&B\\C&D\end{array}\right)\right)=
    \sgn M_\var\left(\left((-I_2)\diamond I_{2(n-1)}\right)
    \left(\begin{array}{cc}A&B\\C&D\end{array}\right)\right).\nn\eea
So by Lemma 2.2 again we have (\ref{zhi30}) holds. So whenever
$\det(B)>0$ or not, (\ref{zhi30}) always holds.

Denote by
$\left(\begin{array}{cc}P^T&0\\0&P^{-1}\end{array}\right)\left(\begin{array}{cc}B^{-1}A&I_n\\B^TC&B^TD\end{array}\right)
    \left(\begin{array}{cc}P&0\\0&(P^{-1})^T\end{array}\right)
    =
\left(\begin{array}{cc}\td{A}&I_n\\\td{C}&\td{D}\end{array}\right)$.
By Claim 2, we have
 \be \td{C}=\diag (0,0,..,0,\lm_1,\lm_2,...,\lm_{n-r}).\ee

Since
$\left(\begin{array}{cc}\td{A}&I_n\\\td{C}&\td{D}\end{array}\right)$
is a symplectic matrix, we have $\td{A}$ and $\td{D}$ are both
symmetric and have the follow forms:
      \bea
      \td{A}=\left(\begin{array}{cc}A_{11}&0\\0&A_{22}\end{array}\right),\quad
       \td{D}=\left(\begin{array}{cc}D_{11}&0\\0&D_{22}\end{array}\right),\nn\eea
where $A_{11}$ and $D_{11}$ are $r\times r$ invertible matrices,
$A_{22}$ and $D_{22}$ are $(n-r)\times (n-r)$ invertible matrices.
So we have
  \be
  \left(\begin{array}{cc}\td{A}&I_n\\\td{C}&\td{D}\end{array}\right)=
  \left(\begin{array}{cc}A_{11}&I_r\\0&D_{11}\end{array}\right)
  \diamond \left(\begin{array}{cc}A_{22}&I_{n-r}\\\Lm
  &D_{22}\end{array}\right),\lb{zhi50}\ee
where $\Lm=\diag(\lm_1,\lm_2,...,\lm_{n-r})$.

Since
$N\left(\begin{array}{cc}A_{11}&I_r\\0&D_{11}\end{array}\right)^{-1}N
\left(\begin{array}{cc}A_{11}&I_r\\0&D_{11}\end{array}\right)
=\left(\begin{array}{cc}I_r&2D_{11}\\0&I_r\end{array}\right)$, by
(\ref{zhi12}) $D_{11}$ is negative definite. So we can joint it to
$-I_r$ by a invertible symmetric matrix path. Then by Lemma 2.2,
Remark 2.1, and computations below Remark 2.1 in Section 2, we have
 \bea \sgn
 M_\var
 \left(\left(\begin{array}{cc}A_{11}&I_r\\0&D_{11}\end{array}\right)\right)&=&
\sgn
 M_\var
 \left(\left(\begin{array}{cc}-I_r&I_r\\0&-I_r\end{array}\right)\right)\nn\\&=&
 r\,\sgn M_\var(N_1(-1,1))\nn\\&=&2r.\lb{zhi51}\eea

Since $M_\var\left(\left(\begin{array}{cc}A_{22}&I_{n-r}\\\Lm
  &D_{22}\end{array}\right)\right)$ is invertible for $\var=0$, for
  $\var>0$ small enough, we have
   \bea &&\sgn M_\var\left(\left(\begin{array}{cc}A_{22}&I_{n-r}\\\Lm
  &D_{22}\end{array}\right)\right)=\sgn M_0\left(\left(\begin{array}{cc}A_{22}&I_{n-r}\\\Lm
  &D_{22}\end{array}\right)\right)\nn\\
  &=&\sgn \left\{\left(\begin{array}{cc}A_{22}&\Lm\\I_{n-r}&D_{22}\end{array}\right)
  \left(\begin{array}{cc}0&-I_{n-r}\\-I_{n-r}
  &0\end{array}\right)\left(\begin{array}{cc}A_{22}&I_{n-r}\\\Lm
  &D_{22}\end{array}\right)+\left(\begin{array}{cc}0&I_{n-r}\\
  I_{n-r}&0\end{array}\right)\right\}\nn\\
  &=&\sgn \left\{2\left(\begin{array}{cc}-A_{22}\Lm&-\Lm\\-\Lm
  &-D_{22}\end{array}\right)\right\}\nn\\
  &=&\sgn \left(\begin{array}{cc}-A_{22}\Lm&-\Lm\\-\Lm
  &-D_{22}\end{array}\right).\lb{zhi41}\eea
Since $\left(\begin{array}{cc}A_{22}&I_{n-r}\\\Lm
  &D_{22}\end{array}\right)$ is a symplectic matrix, we have
  \bea A_{22}D_{22}-\Lm=I_{n-r},\nn\\
       A_{22}\Lm=\Lm A_{22}.\lb{zhi40}
  \eea
  Hence
  \bea
  A_{22}^{-1}\Lm-D_{22}=A_{22}^{-1}(\Lm-A_{22}D_{22})=-A_{22}^{-1}.\nn\eea
So we have
 \bea &&\left(\begin{array}{cc}I_{n-r}&0\\-A_{22}^{-1}
  &I_{n-r}\end{array}\right)\left(\begin{array}{cc}-A_{22}\Lm&-\Lm\\-\Lm
  &-D_{22}\end{array}\right)\left(\begin{array}{cc}I_{n-r}&-A_{22}^{-1}\\0
  &I_{n-r}\end{array}\right)\nn\\&=&\left(\begin{array}{cc}-A_{22}\Lm&0\\0
  &A_{22}^{-1}\Lm-D_{22}\end{array}\right)\nn\\&=&\left(\begin{array}{cc}-A_{22}\Lm&0\\0
  &-A_{22}^{-1}\end{array}\right).\lb{zhi42}\eea
By (\ref{zhi40}), there exist invertible matrix $R$ such that
 \bea &&R^{-1}A_{22}R=\diag(\alpha_1,\alpha_2,...,\alpha_{n-r}), \quad
 \alpha_i\in \R\setminus\{0\}, \; i=1,2,...,n-r,\lb{zhi46}\\
   &&R^{-1}\Lm
   R=\diag(\lm_{i_1},\lm_{i_2},...,\lm_{i_{n-r}}),\;\;\{i_1,i_2,...,i_{n-r}\}=\{1,2,...,n-r\}.\lb{zhi43}\eea
So we have
   \bea R^{-1}(-A_{22}\Lm)R=\diag (-\lm_{i_1}\alpha_1,
   -\lm_{i_2}\alpha_2,...,-\lm_{i_{n-r}}\alpha_{n-r}),\lb{zhi44}\\
    R^{-1}(-A_{22}^{-1})R=\diag
    (-\frac{1}{\alpha_1},-\frac{1}{\alpha_2},...,\frac{1}{\alpha_{n-r}}).\lb{zhi45}\eea
Since $\lm_i\in (-1,0)$ for $i=1,2,...,n-r$, by
(\ref{zhi46})-(\ref{zhi45}) we have
 \be\sgn(-A_{22}\Lm)+\sgn(-A_{22}^{-1})=0.\lb{zhi47}\ee
Hence by (\ref{zhi41}), (\ref{zhi42}) and (\ref{zhi47}) we have
 \bea \sgn M_\var\left(\left(\begin{array}{cc}A_{22}&I_{n-r}\\\Lm
  &D_{22}\end{array}\right)\right)=\sgn(-A_{22}\Lm)+\sgn(-A_{22}^{-1})=0.\lb{zhi48}\eea

Since $\det Q>0$ we can joint it to $I_n$ by a invertible matrix
 path. Hence by Lemma 2.2 and Remark 2.1, (\ref{zhi50}), (\ref{zhi51}) and (\ref{zhi48}), we have
\bea \sgn M_\var
       \left(\left(\begin{array}{cc}B^{-1}A&I_n\\B^TC&B^TD\end{array}\right)\right)&=& \sgn
 M_\var
 \left(\left(\begin{array}{cc}A_{11}&I_r\\0&D_{11}\end{array}\right)\right)+\sgn M_\var\left(\left(\begin{array}{cc}A_{22}&I_{n-r}\\\Lm
  &D_{22}\end{array}\right)\right)\nn\\&=&2r+0\nn\\&=&2r.\lb{zhi52}\eea
 Then by Theorem 2.3, (\ref{zhi30}) and (\ref{zhi52}) we have
 \be i_{L_0}(\ga_{x_\tau})-i_{L_1}(\ga_{x_\tau})=r.\lb{zhi53}\ee
However by (\ref{zhi9}), (\ref{zhi11}) and
$\nu_{L_1}(\ga_{x_\tau})=r$ we have
 \bea i_{L_0}(\ga_{x_\tau})-i_{L_1}(\ga_{x_\tau})=r-1,\nn\eea
which contradicts to (\ref{zhi53}).

Thus we have prove that $k$ can not be $3$. So if $k$ is odd, it
must be $1$. By the same proof we have if $k$ is even, it must be
$2$. Then $\tau\in\{T,\frac{T}{2}\}$. The proof of Theorem 1.1 is
complete. \hfill\hb

 {\bf
Proof of Corollary 1.2.} Since $0<T<\frac{\pi}{||B_0||}$, there is
$\varepsilon>0$ small enough such that
 \bea 0\le B_0\le
||B_0||I_{2n}<(\frac{\pi}{T}-\varepsilon)I_{2n}.\nn\eea

It is easy to see that
 \bea
 \ga_{(\frac{\pi}{T}-\varepsilon)I_{2n}}(t)=\exp(({\frac{\pi}{T}-\varepsilon})tJ)£¬\quad
 \forall t\in[0,\frac{T}{2}].\nn\eea

So we have
 \bea &&\nu_{L_0}( \ga_{(\frac{\pi}{T}-\varepsilon)I_{2n}})=0,\nn\\
   && i_{L_0}((\frac{\pi}{T}-\varepsilon)I_{2n})=0.\nn
 \eea
 Then by (\ref{y20}) and Lemma 3.1 and
Corollary 3.1 we have
 \bea 0\le
i_{-1}(B_0)+\nu_{-1}(B_0)\le
i_{-1}((\frac{\pi}{T}-\varepsilon)I_{2n}) =0.\nn\eea So we have
 \bea i_{-1}(B_0)+\nu_{-1}(B_0)=0.\nn\eea

Hence by the same proof of Theorem 1.1, the conclusions of Corollary
1.2 holds.\hfill\hb

\noindent{\bf Remark 4.1.} Under the same conditions of Theorem 1.3,
if $\int_0^{\frac{T}{2}}H''_{22}(x_T(t))\,dt>0$, by the same proof
of Theorem 1.1, we have \bea \tau\ge
\frac{T}{2(i_{L_0}(B_0)+\nu_{L_0}(B_0))+2}.\nn\eea Moreover, if
$0<T<\frac{\pi}{||B_0||}$ or $i_{L_0}(B_0)+\nu_{L_0}(B_0)=0$, we
have $\tau\in\{T,\frac{T}{2}\}$.

{\bf Proof of Theorem 1.2.} This is the case $n=1$ and $B_0=0$ of
Theorem 1.3, by the proof Theorem 1.3, for any $T>0$ we obtain an
T-periodic brake solution $x_T$ satisfies
 \bea i_{L_0}(\ga_{x_T})\le 1.\eea
If it's minimal period is $\tau=T/k$ for some $k\in \N$, we denote
$x_\tau=x_T|_{[0,\tau]}$. Then by the proof of Theorem 1.3 we have
 \be i_1(\ga_{x\tau
 }^2)+\nu_1(\ga_{{x_\tau}}^2)\ge 2.\lb{q41}\ee

 In the following we prove Theorem 1.2 in 2 steps.

 {\bf Step 1.} For $k=2p+1$ for some $p\ge 0$, we prove that $p=0$.

Firstly by the proof of Theorem 1.3 we have
 \be 1\ge i_{L_0}(\ga_{{x_\tau}}^{2p+1})\ge
 p(i_1(\ga_{{x_\tau}}^2)+\nu_1(\ga_{{x_\tau}}^2)-1)+i_{L_0}(\ga).\lb{q42}\ee
We divide the argument into three cases.

{\bf Case 1.} $i_1(\ga_{x\tau }^2)+\nu_1(\ga_{{x_\tau}}^2)=2$.
      If $\nu_1(\ga_{{x_\tau}}^2)=1$, then $i_1(\ga_{{x_\tau}}^2)=1\in
      2\Z+ 1$. By Lemma 6.3, we have
      $N_1(1,1)\in\Om^0(\ga_{{x_\tau}}^2(\tau))$.
     Since
     \be
     1=i_1(\ga_{{x_\tau}}^2)=i_{L_0}(\ga_{{x_\tau}})+i_{L_1}(\ga_{x\tau
     })+1.\lb{q43}\ee
  By Corollary 2.1 we have
     \be |i_{L_0}(\ga_{{x_\tau}})-i_{L_1}(\ga_{{x_\tau}})|\le 1.\lb{q44}\ee
  Then by (\ref{q43}) and (\ref{q44}) we have
   \be i_{L_0}(\ga_{{x_\tau}})=i_{L_1}(\ga_{{x_\tau}})=0.\ee

So by Theorem 2.1, Lemma 6.2, and (\ref{ddd}), we have \bea
     i_{L_0}(\ga_{{x_\tau}}^3)&=&i_{L_0}(\ga_{{x_\tau}})+i_{e^{2\pi\sqrt{-1}/3}}(\ga_{x\tau
     }^2)\nn\\
                           &=&i_{L_0}(\ga_{{x_\tau}})+i_1(\ga_{{x_\tau}}^2)+S_{N_1(1,1)}(1)\nn\\
                           &=&0+1+1\nn\\
                           &=&2>1\ge i_{L_0}(\ga_{x\tau
                           }^{2p+1}).\eea
Then by Theorem 3.3 we have
 \bea 2p+1< 3.\nn\eea
 Hence $p=0$.

If $\nu_1(\ga_{{x_\tau}}^2)=2$, then $i_1(\ga_{{x_\tau}}^2)=0$. But
now $\ga_{{x_\tau}}^2(\tau)=I_2$, by Lemma 6.3
$i_1(\ga_{{x_\tau}}^2)\in 2\Z+1$, which yields a contradiction. So
this case can not happen. So in Case 1, we have proved $p=0$.

{\bf Case 2.} $i_1(\ga_{x\tau }^2)+\nu_1(\ga_{{x_\tau}}^2)=3$.

  If $\nu_1(\ga_{{x_\tau}}^2)=1$, then
  \bea i_1(\ga_{x\tau
  }^2)=2\in 2\Z.\lb{q50}\eea
   By Lemma
  6.3 we have $N_1(1,-1)\in \Om^0((\ga_{{x_\tau}}^2)(\tau))$. So if $p\ge 1$,
  by Theorem 3.3,
  Theorem 2.1, Lemma 6.2 and (\ref{ddd}), we have
  we have
  \bea 1\ge i_{L_0}(\ga_{{x_\tau}}^{2p+1})
&\ge& i_{L_0}(\ga_{{x_\tau}}^3)\nn\\
  &=&i_{L_0}(\ga_{{x_\tau}})+i_{e^{2\pi\sqrt{-1}/3}}(\ga_{{x_\tau}}^2)\nn\\
                                      &=& i_{L_0}(\ga_{{x_\tau}})+i_1(\ga_{x\tau
                                      }^2)+S_{N_1(1,-1)}(1)\nn\\
                                       &\ge& -1 +2+0\nn\\
                                       &=&1.\eea
  So there must hold
       \bea i_{L_0}(\ga_{x\tau
       })=-1.\nn\eea
    Then by Corollary 2.1 we have
     \bea i_{L_1}(\ga_{{x_\tau}})\le 0.\nn\eea
     So we have
      \bea i_1(\ga_{{x_\tau}}^{2})=i_{L_0}(\ga_{{x_\tau}}+
      i_{L_1}(\ga_{{x_\tau}})+1\le -1+0+1=0,\nn\eea
which contradicts (\ref{q50}). Thus we have $p=0$.

If $\nu_1(\ga_{{x_\tau}}^2)=2$, then \be
i_1(\ga_{{x_\tau}}^2)=1,\quad \ga_{x_\tau}^2(\tau)=I_2.\lb{q51}\ee
 If $p\ge 1$,
  by Theorem 3.3,
  Theorem 2.1, Corollary 3.2, Lemma 6.2 and (\ref{ddd}), we have
  we have
  \bea 1\ge i_{L_0}(\ga_{{x_\tau}}^{2p+1})
&\ge& i_{L_0}(\ga_{{x_\tau}}^{2+1})\nn\\
  &=&i_{L_0}(\ga_{{x_\tau}})+i_{e^{2\pi\sqrt{-1}/3}}(\ga_{{x_\tau}}^2)\nn\\
                                      &=& i_{L_0}(\ga_{{x_\tau}})+i_1(\ga_{{x_\tau}}^2)+S_{I_2}(1)\nn\\
                                       &\ge& -1 +1+1\nn\\
                                       &=&1.\nn\eea
So there must hold
       \bea i_{L_0}(\ga_{x\tau
       })=-1.\nn\eea
Then by Corollary 2.1 we have
     \bea i_{L_1}(\ga_{{x_\tau}})\le 0.\nn\eea
     So we have
      \bea i_1(\ga_{x\tau
      }^{2})=i_{L_0}(\ga_{{x_\tau}}+
      i_{L_1}(\ga_{{x_\tau}})+1\le -1+0+1=0,\nn\eea
which contradicts (\ref{q51}). Thus we have $p=0$.

{\bf Case 3.} $i_1(\ga_{{x_\tau}}^2)+\nu_1(\ga_{x_\tau }^2)\ge 4$.

In this case $i_1(\ga_{{x_\tau}}^2)+\nu_1(\ga_{{x_\tau}}^2)-1\ge 3$.
By Corollary 3.2 we have $i_{L_0}\ge -1$. So by (\ref{q42}) we have
 \be p\le 2/3,\ee
 which yields $p=0$.
So we finish Step 1.

{\bf Step 2.} For $k=2p+2$ for some $p\ge 0$, we prove that $p=0$.

 In fact, apply Bott-type iteration formula of Theorem 2.1 to the the case of the iteration time equals to 4
  and note that by Corollary 3.1
 $i_{\sqrt{-1}}(\ga_{{x_\tau}})\ge 0$. Then by the same argument of
 Step 1, we can prove that $p=0$.

 Thus by Steps 1 and 2, Theorem 1.2 is proved.\hfill\hb

A natural question is that can we prove the minimal period is $T$ in
this way? We have the following remark.

\noindent{\bf Remark 4.2.} Only use the Maslov-type index theory to
estimate the iteration time of the $T$-periodic brake solution $x_T$
obtained by the first 4 steps in the proof of Theorem 1.3 with
$B_0=0$, we can not hope to prove $T$ is the minimal period of
$x_T$. Even $H''(z)>0$ for all $z\in \R^{2n}\setminus\{0\}$. For
$n=1$ and $T=4\pi$, we can not exclude the following case:
   \bea
    &&x_T(t)=\left(\begin{array}{c}\sin t\\ \cos
    t\end{array}\right),\nn\\
     && H'(x_T(t))=x_T(t),\nn\\
    &&H''(x_T(t))\equiv I_{2n}.\nn\eea
It is easy to check that $\ga_{x_T}(t)=R(t)$ for $t\in[0,2\pi]$.
Hence by Lemma 5.1 of \cite{Liu4} or the proof of Lemma 3.1 of
\cite{LZZ} we have
  \bea
  i_{L_0}(\ga_{x_T})=\sum_{0<s<2\pi}\nu_{L_0}(\ga_{x_T})(s)=1.\nn\eea
 In this case the minimal period of $x_T$ is
$\frac{T}{2}$. Similarly for $n>1$ we can construct examples to
support this remark.

\setcounter{equation}{0}
\section{Proof of Theorems 1.4-1.5 and Corollary 1.4}

In this section we study the minimal period problem for symmeytric
brake orbit solutions of the even  reversible Hamiltonian system
(\ref{1.1}) and complete the proof of Theorems 1.4-1.5 and Corollary
1.4.

For $T>0$, let $E_{T}=\{x\in
W^{1/2,2}(S_{\tau},\R^{2n})|\,x(-t)=Nx(t)\;a.e.\;t\in \R\}$ with the
usual $W^{1/2,2}$ norm and inner product. Correspondingly $\hat{E}$
and $\td{E}$ are defined to be the symmetric ones and the
$\frac{T}{2}$-periodic ones in $E_{T}$ respectively. Also
$\{P_{T,m}\}$ and $\{\hat{P}_m\}$ are the Galerkin approximation
scheme w.r.t. $A_{T}$ and $\hat{A}$ respectively, where
$\{P_{T,m}\}$, $\{\hat{P}_m\}$, $A_T$, and $\hat{A}$ are defined
 by the same way as in Section 2, we only need to replace $\tau$ by $T$.

 For $z\in E_{T}$, we define
  \be f(z)=\frac{1}{2}\langle A_Tz,z\rangle-\int_0^TH(z)dt.\ee
For $z\in \hat{E}$, we define
  \be \hat{f}(z)=\frac{1}{2}\langle \hat{A}z,z\rangle-\int_0^TH(z)dt.\ee

We have the following lemma.

\noindent{\bf Lemma 5.1.} {\it Let $z\in \hat{E}$. If
$\hat{f}'(z)=0$, then $f'(z)=0$.}

{\bf Proof.} Let $z\in \hat{E}$ and $\hat{f}'(z)=0$. So for any
            $y\in\hat{E}$ we have
          \be \langle \hat{f}'(z),y\rangle=
          \int_0^TJ\dot{z}(t)\cdot y(t)\,dt-\int_0^TH'(z(t))\cdot y(t)\,dt=0,\quad
          \forall y\in\hat{E}.\ee

          Since $H$ is even and  $z\in \hat{E}$, we have
          \be H'(z(t+\frac{T}{2}))=H'(-z(t))=-H'(z(t)).\lb{n2}\ee
          So $H'(z)\in \hat{E}$ and

        \be \langle f'(z),y\rangle=
          \int_0^TJ\dot{z}(t)\cdot y(t)\,dt-\int_0^TH'(z(t))\cdot y(t)\,dt=0,\quad
          \forall y\in\td{E}.\lb{n3}\ee
By (\ref{n2}) and (\ref{n3}), we have
            \be \langle f'(z),y\rangle=
         \int_0^TJ\dot{z}(t)\cdot y(t)\,dt-\int_0^TH'(z(t))\cdot y(t)\,dt=0,\quad
          \forall y\in E_T.\ee
 Hence $f'(z)=0$  \hfill\hb

By Lemma 5.1 and arguments in the proof of Theorem 1.3 in Section 4,
to look for the $T$-period symmetric solutions of (\ref{1.1}) is
equivalent to look for critical points of $\hat{f}$.

{\bf Proof of Theorem 1.5.} For any given $T>0$, we prove the
existence of $T$-periodic symmetric  brake orbit solution of
(\ref{1.1}) whose minimal period satisfies the inequalities in the
conclusion of Theorem 1.5. Since the proof of existence of
$T$-periodic symmetric  brake orbit solution $x_T$ of (\ref{1.1}) is
similar to that of the proof of Theorem 1.3, we will only give the
sketch. We divide the proof into several steps.

 {\bf Step 1.} Similarly as Step 1 in the proof of Theorem 1.3,  for any $K>0$
 we can truncate the function $\hat{H}$
 suitably and evenly to $\hat{H}_K$ such that it satisfies the growth condition (\ref{n1}).
Correspondingly we obtain a new even and reversible function $H_K$
satisfies condition (\ref{n1}).

Set \be \hat{f}_K(z)=\frac{1}{2}\langle \hat{A} z,z\rangle-\int_0^T
H_K(z)dt,\qquad \forall z\in\hat{E}.\ee
 Then $\hat{f}_K\in
C^2(\hat{E},\R)$ and \be \hat{f}_K(z)=\frac{1}{2}\langle
(\hat{A}-\hat{B}_0) z,z\rangle-\int_0^T \hat{H}_K(z)dt,\qquad
\forall z\in\hat{E},\ee where $\hat{B}_0$ is the selfadjoint linear
compact operator on $\hat{E}$ defined by \be \langle
\hat{B}_0z,z\rangle=\int_0^TB_0z(t)\cdot z(t)\,dt.\ee

{\bf Step 2.} For $m>0$, let $\hat{f}_{Km}=\hat{f}|\hat{E}_m$, where
$\hat{E}_m=\hat{P}_m \hat{E}$. Set
 \bea &&X_m=M^-(\hat{P}_m(\hat{A}-\hat{B}_0)\hat{P}_m)\oplus M^0(\hat{P}_m(\hat{A}-\hat{B}_0)\hat{P}_m),\nn\\
        &&  Y_m=M^+(\hat{P}_m(\hat{A}-\hat{B}_0)\hat{P}_m).\nn\eea

  By the same argument of Step 2 in
the proof of Theorem 1.3, we can
  show that $\hat{f}_Km$ satisfies the hypotheses of Theorem 4.1.
Moreover, we obtain a critical point $x_{Km}$ of
 $\hat{f}_{Km}$ with critical value $C_Km$ which satisfies
   \bea m^-(x_{Km})\le \dim X_m+1.\lb{duan1}\eea
and
 \be \delta\le C_{Km}\le \frac{1}{2}||\hat{A}-\hat{B}_0||r_1^2,\ee
 where $\delta$ is a positive number depending on $K$ and $r_1>0$
is independent of $K$ and $m$.

{\bf Step 3.} We prove that there exists a symmetric $T$-periodic
brake orbit solution $x_T$ of (\ref{1.1}) which satisfies
 \be
i_{\sqrt{-1}}^{L_0}(\ga_{x_T})\le
i_{\sqrt{-1}}^{L_0}(B_0)+\nu_{\sqrt{-1}}^{L_0}(B_0)+1.\ee

From the proof of Theorem 1.3 we have $f_K$ satisfies $(PS)_c^*$
condition for $c\in\R$, by the same proof of Lemma 5.1, we have
$\hat{f}_K$ satisfies $(PS)_c^*$ condition for $c\in\R$, i.e., any
sequence ${z_m}$ such that $z_m\in\hat{E}_m$, $\hat{f}_{Km}'(z_m)\to
0$ and $\hat{f}_{Km}(z_m)\to c$ possesses a convergent subsequence
in $\hat{E}$. Hence in the sense of subsequence we have
 \be x_{Km}\to x_K,\qquad \hat{f}_K(x_K)=c_K,\qquad
 \hat{f}'_K(x_K)=0.\lb{duan2}\ee
By similar argument as in \cite{Ra2}, $x_K$ is a classical
nonconstant symmetric $T$-periodic solution of

 \be \dot{x}=JH_K'(x), \quad x\in\R^{2n}.\ee

 Set $B_K(t)=H''_K(x_K(t))$, Then $B_K\in
C(S_{T/2},\mathcal{L}_s(\R^{2n}))$. Let $\hat{B}_K$ be the operator
defined by the same way of the definition of $\hat{B}_0$. It is easy
to show that \be ||\hat{f}''(z)-(\hat{A}-\hat{B}_K)||\to 0\qquad
{\rm as}\;\; ||z-x_K||\to 0.\ee

So for $0<d<\frac{1}{4}||(A_T-B_{K_T})^\#||^{-1}$, there exists
$r_2>0$ such that
      \be ||\hat{f}_{Km}''(z)-\hat{P}_m(\hat{A}-\hat{B}_K)\hat{P}_m||\le ||\hat{f}''(z)-(\hat{A}-\hat{B}_K)||<\frac{1}{2}d, \qquad
      \forall z\in\{z\in \hat{E}:||z-x_K||\le r_2\}.\ee

Then for $z\in \{z\in \hat{E}: ||z-x_K||\le r_2\}\cap \hat{E}_m$,
$\forall u\in
M^-_d(\hat{P}_m(\hat{A}-\hat{B}_T)\hat{P}_m)\setminus\{0\}$, we have
   \bea \langle \hat{f}_{Km}''(z)u,u\rangle &\le& \langle \hat{P}_m(\hat{A}-\hat{B}_K)\hat{P}_mu,u\rangle
   +\|\hat{f}_{Km}''(z)-\hat{P}_m(\hat{A}-\hat{B}_K)\hat{P}_m\|\|u\|^2\nn\\
       &\le& -\frac{1}{2}d\|u\|^2.\nn\eea
So we have
   \be m^-(\hat{f}_{Km}''(z))\ge \dim
   M^-_d(\hat{P}_m(\hat{A}-\hat{B}_K)\hat{P}_m).\lb{duan3} \ee

By Theorem 3.1, Remark 3.1,  there is $m^*>0$ such that for $m\ge
m^*$ we have
      \bea &&\dim
      X_m=mn+i_{\sqrt{-1}}^{L_0}(B_0)+\nu_{\sqrt{-1}}^{L_0}(B_0),\lb{duan5}\\
      &&\dim
      M^-_d(\hat{P}_m(\hat{A}-\hat{B}_K)\hat{P}_m)=mn+i_{\sqrt{-1}}^{L_0}(B_K).\lb{duan6}\eea
Then by (\ref{duan1}), (\ref{duan2}), and
(\ref{duan3})-(\ref{duan6}), we have
       \be i_{\sqrt{-1}}^{L_0}(B_K)\le
       i_{\sqrt{-1}}^{L_0}(B_0)+\nu_{\sqrt{-1}}^{L_0}(B_0)+1.\ee

By the similar argument as in the section 6 of \cite{Ra2}, there is
a constant $M_3$ independent of $K$ such that $||x_K||_\infty\le
M_3$. Choose $K>M_3$. Then $x_K$ is a non-constant symmetric
$T$-periodic brake orbit solution of the problem (\ref{1.1}). From
now on in the proof of Theorem 1.2, we write $B=B_K$ and $x_T=x_K$.
Then $x_T$ is a non-constant symmetric $T$-periodic solution of the
problem (\ref{1.1}), and $B$ satisfies
    \be i_{\sqrt{-1}}^{L_0}(\ga_{x_T})=_{\sqrt{-1}}^{L_0}(B)\le
       i_{\sqrt{-1}}^{L_0}(B_0)+\nu_{\sqrt{-1}}^{L_0}(B_0)+1.\lb{duan7}\ee

{\bf Step 4.} Finish the proof of Theorem 1.5.

      Since $x_T$ obtained in Step 3 is a nonconstant and symmetric
      $T$-period brake orbit solution, its minimal period $\tau=\frac{T}{4r+s}$
      for some nonnegative integer $r$ and $s=1$ or $s=3$. We now estimate $r$.

  We denote by $x_\tau=x_T|_{[0,\tau]}$, then it is a symmetric
  period solution of (\ref{1.1}) with the minimal $\tau$ and
  $X_T=x_\tau^{4r+s}$ being the $4r+s$ times iteration of $x_\tau$.
  As in Section 1, let $\ga_{x_T}$ and $\ga_{x_\tau}$ the symplectic path
  associated to $(\tau, x)$ and $(T,x_T)$ respectively.
  Then $\ga_{x_\tau}\in C([0,\frac{\tau}{4}],\Sp(2n))$ and
 $\ga_{x_T}\in C([0,\frac{T}{4}],\Sp(2n))$. Also we have
 $\ga_{x_T}=\ga_{x_\tau}^{4r+s}$, which is the $4r+s$ times iteration of
 $\ga_{x_\tau}$.

By (\ref{duan7}) we have
  \be
i_{\sqrt{-1}}^{L_0}(\ga_{x_\tau}^{4r+s})\le
       i_{\sqrt{-1}}^{L_0}(B_0)+\nu_{\sqrt{-1}}^{L_0}(B_0)+1.\lb{duan9}\ee

Since $x_\tau$ is also a nonconstant symmetric periodic solution of
(\ref{1.1}). It is clear that
  \bea \nu_{-1}(x_\tau^2)&\ge& 1.\lb{duan8}\eea
Since $\hat{H}$ satisfies condition (H5) and $B_0$ is semipositive,
by Corollary 3.1 of \cite{Z1} (also by Theorem 6.2) we have
     \be i_{-1}(\ga_{x_\tau}^2)\ge 0.\lb{duan10}\ee
By Corollary 3.2 of \cite{Z1} (cf. aslo \cite{Liu1}), we have
       \be i_1(\ga_{x_\tau}^2)+\nu_1(\ga_{x_\tau}^2)\ge n.\lb{duan11}\ee
It is easy to see that \be
\ga_{x_\tau}^4(\frac{\tau}{2}+t)=\ga_{x_\tau}^2(t)\,\ga_{x_\tau}^2(\frac{\tau}{2}),\qquad
\forall t\in [0,\frac{\tau}{2}].\lb{duan19}\ee So by Theorem 6.1 of
Bott-type iteration formula we have
        \bea
         i_1(\ga_{x_\tau}^4)+\nu_1(\ga_{x_\tau}^4)&=&
         i_1(\ga_{x_\tau}^2)+\nu_1(\ga_{x_\tau}^2)+i_{-1}(\ga_{x_\tau}^2)+\nu_{-1}(\ga_{x_\tau}^2)\nn\\
         &\ge& n+0+1\nn\\
          &=& n+1.\lb{duan12}\eea

If $r\ge 1$, then by Theorems 2.2 and 6.2 and (\ref{duan12}) we have
       \bea
       i_{-1}(\ga_{x_\tau}^{4r}) &=&i_{-1}((\ga_{x_\tau}^{2})^{2p})\nn\\
                                 &=&\sum_{j=1}^ri_{\om_{2r}^{2j-1}}(\ga_{x_\tau}^4)\nn\\
                                 &\ge&
                                 \sum_{j=1}^r(i_1(\ga_{x_\tau}^4)+\nu_1(\ga_{x_\tau}^4)-n)\lb{duan16}\\
                                  &=&
                                  r(i_1(\ga_{x_\tau}^4)+\nu_1(\ga_{x_\tau}^4)-n)\nn\\
                                  &\ge& r,\lb{duan13}\eea
where $\om_{2r}=e^{\pi\sqrt{-1}/(2r)}$ as defined in Theorem 2.2.

By Theorem 3.2, we have
             \be  i_{\sqrt{-1}}^{L_0}(\ga_{x_\tau}^{4r+s})\ge
             i_{\sqrt{-1}}^{L_0}(\ga_{x_\tau}^{4r}).\lb{duan14}\ee
Then (\ref{duan9}), (\ref{duan13}) and (\ref{duan14}) yield
        \be r\le
        i_{\sqrt{-1}}^{L_0}(B_0)+\nu_{\sqrt{-1}}^{L_0}(B_0)+1.\lb{duan15}\ee

Thus for $i_{\sqrt{-1}}^{L_0}(B_0)+\nu_{\sqrt{-1}}^{L_0}(B_0)$ is
odd, by (\ref{duan15}) we have
 \be 4r+s\le 4r+3\le
 4(i_{\sqrt{-1}}^{L_0}(B_0)+\nu_{\sqrt{-1}}^{L_0}(B_0))+7.\lb{duan22}\ee

{\bf Claim 3.} For
$i_{\sqrt{-1}}^{L_0}(B_0)+\nu_{\sqrt{-1}}^{L_0}(B_0)$ is even, the
equality in (\ref{duan15}) can not hold.

Otherwise, $r\ge 1$ and the equality in (\ref{duan16}) holds i.e.,
 \be
 i_{\om_{2r}^{2j-1}}(\ga_{x_\tau}^4)=i_1(\ga_{x_\tau}^4)+\nu_1(\ga_{x_\tau}^4)-n=1,
 \quad j=1,2,...,r.\lb{duan17}\ee

By the definition of $\om_{2r}$, we have $\om_{2r}^{2j-1}\neq -1$
for $j=1,2,...,r$. So by  (\ref{duan17}) and 2 of Theorem 6.2, we
have $I_{2p}\diamond N_1(1,-1)^{\diamond q}\diamond K\in
\Om^0(\ga_{x_\tau}^4(\tau))$ for some non-negative integers $p$ and
$q$ satisfying $0\le p+q\le n$ and $K\in \Sp(2(n-p-q))$ with
$\sg(K)\in \U\setminus\{1\}$ satisfying the condition that all
eigenvalues of $K$ located with the arc between $1$ and $\om_{2r}$
in $\U^+\setminus\{\pm 1\}$ possess total multiplicity $n-p-q$. So
there are no eigenvalues of $K$ on the arc between $\om_{2r}^{2j-1}$
and $-1$ except $\om_{2r}^{2r-1}$ with $r=1$. However, whether
$\om_{2r}^{2r-1}\in\sg(\ga_{x_\tau}^4(\tau))$ or not, we always have
 \bea &&S^+_{\ga_{x_\tau}^4(\tau)}(\om_{2r}^{2r-1})=0,\\
     &&i_{\om_{2r}^{2r-1}}(\ga_{x_\tau}^4)=1.\eea
So (\ref{ddd}) and Lemma 6.2, we have
 \bea i_{-1}(\ga_{x_\tau}^4)&=&
 i_{\om_{2r}^{2r-1}}(\ga_{x_\tau}^4)+S^+_{\ga_{x_\tau}^4(\tau)}(\om_{2r}^{2r-1})\nn\\
     &=& 1+0=1.\lb{duan18}\eea
But by (\ref{duan19}), Lemma 6.1, and Theorem 6.1, we have
       \bea
       i_{-1}(\ga_{x_\tau}^{4r})&=&i_{-1}((\ga_{x_\tau}^{2r})^2)\nn\\
                                &=&i_{\sqrt{-1}}(\ga_{x_\tau}^{2r})+i_{-\sqrt{-1}}(\ga_{x_\tau}^{2r})\nn\\
                                &=&
                                2i_{\sqrt{-1}}(\ga_{x_\tau}^{2r}).\lb{duan20}\nn\eea
 Then  $i_{-1}(\ga_{x_\tau}^{4r})$ is an even integer, which yields a contradiction to (\ref{duan18}).
So Claim 3 holds, and we have
 \be r\le i_{\sqrt{-1}}^{L_0}(B_0)+\nu_{\sqrt{-1}}^{L_0}(B_0).\ee
Hence
\be 4r+s\le 4r+3\le
 4(i_{\sqrt{-1}}^{L_0}(B_0)+\nu_{\sqrt{-1}}^{L_0}(B_0))+3.\lb{duan21}\ee

   Theorem 1.5 holds from (\ref{duan22}) and (\ref{duan21}).\hfill\hb

{\bf Proof of Theorem 1.4.} This is the case $B_0\equiv0$ of Theorem
1.2. From Theorem 3.1 it is easy to see that \be
i_{\sqrt{-1}}^{L_0}(0)=0,\qquad \nu_{\sqrt{-1}}^{L_0}(0)=0.\ee
 Then
$i_{\sqrt{-1}}^{L_0}(0)+\nu_{\sqrt{-1}}^{L_0}(0)=0$ and is also
even. So Theorem 1.4 holds from Theorem 1.5.\hfill\hb

 {\bf
Proof of Corollary 1.2.} Since $0<T<\frac{\pi}{||B_0||}$, there is
$\varepsilon>0$ small enough such that
 \be 0\le B_0\le
||B_0||I_{2n}<(\frac{\pi}{T}-\varepsilon)I_{2n}.\lb{y20}\ee

It is easy to see that
 \be
 \ga_{(\frac{\pi}{T}-\varepsilon)I_{2n}}(t)=\exp(({\frac{\pi}{T}-\varepsilon})tJ)£¬\quad
 \forall t\in[0,\frac{T}{4}].\ee

 Since
\be\nu_{L_0}(\exp(({\frac{\pi}{T}-\varepsilon})tJ))=0,\qquad \forall
t\in[0,\frac{T}{2}].\ee We have
 \be i_{L_0}(\ga_{(\frac{\pi}{T}-\varepsilon)I_{2n}}^2)=0,\;
 i_{L_0}(\ga_{(\frac{\pi}{T}-\varepsilon)I_{2n}})=0.\ee
So by Theorem 2.1 we
have
 \bea &&i_{\sqrt{-1}}^{L_0}((\frac{\pi}{T}-\varepsilon)I_{2n})=i_{L_0}(\ga_{(\frac{\pi}{T}-\varepsilon)I_{2n}}^2)
 -i_{L_0}(\ga_{(\frac{\pi}{T}-\varepsilon)I_{2n}})=0.
 \eea
 Then by (\ref{y20}) and Lemma 3.1 and
Corollary 3.1 we have
 \be 0\le
i_{-1}(B_0,\frac{T}{2})+\nu_{-1}(B_0,\frac{T}{2})\le
i_{-1}((\frac{\pi}{T}-\varepsilon)I_{2n},\frac{T}{2}) =0.\ee So we
have
 \be i_{-1}(B_0,\frac{T}{2})+\nu_{-1}(B_0,\frac{T}{2})=0.\ee

Hence by Theorem 1.1 or Corollary 1.1, the conclusion of Corollary
1.2 holds.\hfill\hb

Also a natural question is that can we prove the minimal period is
$T$ in this way? We have the following remark.

\noindent{\bf Remark 5.1.} Only use the Maslov-type index theory to
estimate the iteration time of the symmetric $T$-periodic brake
solution $x_T$ obtained in the proof of Theorem 1.5 with $B_0=0$, we
can not hope to prove $T$ is the minimal period of $x_T$. Even
$H''(z)>0$ for all $z\in \R^{2n}\setminus\{0\}$. For $n=1$ and
$T=6\pi$, we can not exclude the following case:
   \bea
    &&x_T(t)=\left(\begin{array}{c}\sin t\\ \cos
    t\end{array}\right),\nn\\
     && H'(x_T(t))=x_T(t),\nn\\
    &&H''(x_T(t))\equiv I_{2n}.\eea
It is easy to check that $\ga_{x_T}(t)=R(t)$ for $t\in[0,3\pi]$.
Hence by Theorem 2.1 and Lemma 5.1 of \cite{Liu4} or the proof of
Lemma 3.1 of \cite{LZZ} we have
  \be
  i_{\sqrt{-1}}^{L_0}(\ga_{x_T})=\sum_{3\pi/4\le s<3\pi}\nu_{L_0}(\ga_{x_T})(s)=1.\ee
 In this case the minimal period of $x_T$ is
$\frac{T}{3}$. Similarly for $n>1$ we can construct examples to
support this remark.

\setcounter{equation}{0}
\section{Appendix on Maslov-type indices $(i_\om,\nu_\om)$  }

 We first recall briefly the Maslov-type index theory of
$(i_\om,\nu_\om)$. All the details can be found in \cite{Long5}.

 For any
$\om\in \U$, the following codimension 1 hypersuface in $\Sp(2n)$ is
defined by:
    $$ \Sp(2n)_\om^0=\{M\in \Sp(2n)|\det(M-\om I_{2n})=0\}.$$
For any two continuous path $\xi$ and $\eta$: $[0,\tau]\to\Sp(2n)$
with $\xi(\tau)=\eta(0)$, their joint path is defined by
    \be  \eta * \xi(t)=\left\{\begin{array}{lr} \xi(2t)\qquad &{\rm
     if}\,
     0\le t\le \frac{\tau}{2},\\ \eta(2t-\tau)\quad &{\rm if}\,
     \frac{\tau}{2}\le t\le \tau .\end{array}\right.\lb{zzz}\ee
Given any two $(2m_k\times 2m_k)$- matrices of square block form
$M_k=\left(\begin{array}{cc}A_k&B_k\\C_k&D_k\end{array}\right)$ for
$k=1,2$, as in \cite{Long5}, the $\diamond$-product of $M_1$ and
$M_2$ is defined by the following $(2(m_1+m_2)\times
2(m_1+m_2))$-matrix $M_1\diamond M_2$:
     $$M_1\diamond M_2=\left(\begin{array}{cccc}A_1&0&B_1&0\\0&A_2&0&B_2\\
        C_1&0&D_1&0\\0&C_2&0&D_2\end{array}\right).$$
A special path $\xi_n$ is defined by
        $$\xi_n(t)=\left(\begin{array}{cc}2-\frac{t}{\tau}&0\\0&(2-\frac{t}{\tau})^{-1}\end{array}\right)^{\diamond
        n}
        , \qquad \forall t\in[0,\tau].$$
{\bf Definition 6.1.} For any $\om\in\U$ and $M\in\Sp(2n)$, define
     \be \nu_\om(M)=\dim_\C\ker(M-\om I_{2n}).\ee
For any $\ga\in \mathcal{P}_\tau(2n)$, define
  \be \nu_\om(\ga)=\nu_\om(\ga(\tau)).\ee
If $\ga(\tau)\notin\Sp(2n)_\om^0$, we define
     \be i_\om(\ga)=[\Sp(2n)_\om^0\,:\,\ga*\xi_n],\lb{m1}\ee
where the right-hand side of (\ref{m1}) is the usual homotopy
intersection number and the orientation of $\ga*\xi_n$ is its
positive time direction under homotopy with fixed endpoints.

If $\ga(\tau)\in\Sp(2n)_\om^0$, we let $\mathcal{F}(\ga)$ be the set
of all open neighborhoods of $\ga$ in $\mathcal{P}_\tau(2n)$, and
define
 \be i_\om(\ga)=\sup_{U\in\mathcal{F}(\ga)}\inf\{i_\om(\beta)|\,\beta(\tau)\in U \,{\rm and}\,
 \beta(\tau)\notin\Sp(2n)_\om^0\}.\ee
Then $(i_\om(\ga),\nu_\om(\ga))\in \Z\times \{0,1,...,2n\}$, is
called the index function of $\ga$ at $\om$.

\noindent{\bf Lemma 6.1.} (Lemma 5.3.1 of \cite{Long5}) {\it For any
$\ga\in \mathcal{P}_\tau(2n)$ and $\om\in \U$, there hold
 \bea i_\om(\ga)=i_{\bar{\om}}(\ga),\qquad
 \nu_\om(\ga)=\nu_{\bar{\om}}(\ga).\eea}
\quad As in \cite{Long2}, for any $M\in\Sp(2n)$ we define
 \bea \Om(M)=\{P\in \Sp(2n)&|&\sg(P)\cap\U=\sg(M)\cap \U\nn\\
                            && {\rm and}\,
                            \nu_\lm(P)=\nu_\lm(M),\;\; \forall
                            \lm\in\sg(M)\cap\U\}.\eea
We denote by $\Om^0(M)$ the path connected component of $\Om(M)$
containing $M$, and call it the {\it homotopy component} of $M$ in
$\Sp(2n)$.

The following symplectic matrices were introduced as {\it basic
normal forms} in \cite{Long5}:
      \bea D(\lm)=
      \left(\begin{array}{cc}\lm&0\\0&\lm^{-1}\end{array}\right),\qquad
     && \lm=\pm 2,\\
      N_1(\lm,b)=\left(\begin{array}{cc}\lm&b\\0&\lm\end{array}\right),\qquad
     && \lm=\pm1,\,b=\pm1,\,0,\\
      R(\theta)=\left(\begin{array}{cc}\cos(\theta)&-\sin(\theta)\\
      \sin(\theta)&\cos(\theta)\end{array} \right), \qquad
      &&\theta\in(0,\pi)\cup(\pi,2\pi),\\
      N_2(\om,b)=\left(\begin{array}{cc}R(\theta)&b\\0&R(\theta)\end{array}\right),
        \qquad
      &&\theta\in(0,\pi)\cup(\pi,2\pi),\eea
where $b= \left(\begin{array}{cc}b_1&b_2\\b_3&b_4\end{array}\right)$
with $b_i\in\R$ and $b_2\neq b_3$.

For any $M\in \Sp(2n)$ and $\om\in\U$, {\it splitting number} of $M$
at $\om$ is defined by
   \be S_M^{\pm}=\lim_{\epsilon\to 0^+}
   i_{\om\exp(\pm\sqrt{-1}\epsilon)}(\ga)-i_\om(\ga)\ee
for any path $\ga\in\mathcal{P}_\tau(2n)$ satisfying $\ga(\tau)=M$.

Splitting numbers possesses the following properties.

\noindent{\bf Lemma 6.2.} (cf. \cite{Long4}, Lemma 9.1.5 and List
9.1.12 of \cite{Long5}) {\it Splitting number $S_M^{\pm}(\om)$ are
well defined; that is they are independent of the choice of the path
$\ga\in\mathcal{P}_\tau(2n)$ satisfying $\ga(\tau)=M$. For
$\om\in\U$ and $M\in\Sp(2n)$, $S_N^{\pm}(\om)$ are constant for all
$N\in\Om^0(M)$. Moreover we have

(1) $(S_M^+(\pm1),S_M^-(\pm1))=(1,1)$ for $M=\pm N_1(1, b)$ with
$b=1$ or $0$;

(2) $(S_M^+(\pm1),S_M^-(\pm1))=(0,0)$ for $M=\pm N_1(1, b)$ with
$b=-1$;

 (3) $(S_M^+(e^{\sqrt{-1}\theta}),S_M^-(e^{\sqrt{-1}\theta}))=(0,1)$ for
   $M=R(\theta)$ with $\theta\in(0,\pi)\cup (\pi,2\pi)$;

(4) $(S_M^+(\om),S_M^-(\om))=(0,0)$
 for $\om\in \U\setminus \R$ and $M=N_2(\om, b)$ is {\bf trivial} i.e., for

  sufficiently small
 $\alpha>0$, $MR((t-1)\alpha)^{\diamond n}$ possesses no eigenvalues
 on $\U$ for $t\in [0,1)$.

(5) $(S_M^+(\om),S_M^-(\om)=(1,1)$
 for $\om\in \U\setminus \R$ and $M=N_2(\om, b)$ is {\bf non-trivial}.

 (6) $(S_M^+(\om),S_M^-(\om)=(0,0)$
for any $\om\in\U$ and $M\in \Sp(2n)$ with $\sg(M)\cap
\U=\emptyset$.

 (7) $S_{M_1\diamond M_2}^\pm(\om)=S_{M_1}^\pm(\om)+S_{M_2}^\pm(\om)$, for any $M_j\in \Sp(2n_j)$ with
    $j=1,2$ and $\om\in\U$.}

By the definition of splitting numbers and Lemma 6.2, for $0\le
\theta_1<\theta_2<2\pi$ and $\ga\in\mathcal{P}_\tau(2n)$ with
$\tau>0$, we have
 \bea
 i_{\exp(\sqrt{-1}\theta_2)}(\ga)&=&i_{\exp(\sqrt{-1}\theta_1)}+S^+_{\ga(\tau)}(e^{\sqrt{-1}\theta_1})\nn\\
&+&
\sum_{\theta\in(\theta_1,\theta_2)}\left(S^+_{\ga(\tau)}(e^{\sqrt{-1}\theta})-S^-_{\ga(\tau)}(e^{\sqrt{-1}\theta})\right)
 -S^-_{\ga(\tau)}(e^{\sqrt{-1}\theta_2}).\lb{ddd}\eea

For any symplectic path $\ga\in\mathcal{P}_\tau(2n)$ and $m\in \N$,
we define its $m$th iteration in the periodic boundary sense
$\ga(m): [0,m\tau]\to \Sp(2n)$ by
   \be \ga(m)(t)=\ga(t-j\tau)\ga(\tau)^j\qquad {\rm for}\,
   j\tau\le t\le (j+1)\tau,\;j=0,1,...,m-1.\ee

\noindent{\bf Definition 6.2.}(cf.\cite{Long4}, \cite{Long5}) For
any $\ga\in\mathcal{P}_\tau(2n)$ and $\om\in\U$, we define
   \be (i_\om(\ga,m),\nu_\om(\ga,m))=(i_\om(\ga(m)),\nu_\om(\ga(m))),\qquad\forall m\in\N.\ee
We have the following Bott-type iteration formula.

\noindent{\bf Theorem 6.1.} (cf. \cite{Long4}, Theorem 9.2.1 of
\cite{Long5}) {\it For any $\tau>0$, $\ga\in\mathcal{P}_\tau(2n)$,
$z\in \U$, and $m\in \N$,
    \bea i_z(\ga, m)=\sum_{\om^k=z}i_\om(\ga),\qquad
    \nu_z(\ga,m)=\sum_{\om^m=z}\nu_\om(\ga).\eea}
By Theorem 8.1.4 of \cite{Long5}, we have the following Lemma.

{\noindent{\bf Lemma 6.3.} {\it For $\ga\in \mathcal{P}_\tau(2)$
with $\tau>0$, the following results hold.

1. If $N_1(1,1)\in \Om^0(\ga(\tau))$, then
  \bea && i_1(\ga,m)=m(i_1(\ga)+1)-1,\qquad \nu_1(\ga, m)=1,\quad
  \forall m\in\N, \\
    && i_1(\ga)\in 2\Z+1.\eea

2. If $N_1(1,1)\in \Om^0(\ga(\tau))$, then
 \bea &&i_1(\ga,m)=m(i_1(\ga)+1)-1,\qquad \nu_1(\ga, m)=2,\quad
  \forall m\in\N, \\
    &&i_1(\ga)\in 2\Z+1.\eea

3. If $N_1(1,-1)\in \Om^0(\ga(\tau))$, then
  \bea &&i_1(\ga,m)=m(i_1(\ga),\qquad \nu_1(\ga, m)=1,\quad
  \forall m\in\N, \\
   && i_1(\ga)\in 2\Z.\eea}

Denote by $\U^+=\{\om\in\U|\,  Im\, \om\ge 0\}$ and
$\U^-=\{\om\in\U|\,  Im\, \om\le 0\}$.
 The following theorem was proved by Liu and Long in
\cite{LiLo1,LiLo2}, which plays a important role in the proof of our
main results in Sections 4-5.

\noindent{\bf Theorem 6.2.} (Theorem 10.1.1 of \cite{Long5}) {\it

1. For any $\ga\in \mathcal{P}_\tau(2n)$ and $\om\in
\U\setminus\{1\}$, it always holds that
 \be i_1(\ga)+\nu_1(\ga)-n\le i_\om(\ga)\le
 i_1(\ga)+n-\nu_\om(\ga).\lb{mm}\ee

 2. The left equality in (\ref{mm}) holds for some $\om\in
\U^+\setminus\{1\}$ (or $\U^-\setminus\{1\}$) if and only if
$I_{2p}\diamond N_1(1,-1)^{\diamond q}\diamond K\in
\Om^0(\ga(\tau))$ for some non-negative integers $p$ and $q$
satisfying $0\le p+q\le n$ and $K\in \Sp(2(n-p-q))$ with $\sg(K)\in
\U\setminus\{1\}$ satisfying the condition that all eigenvalues of
$K$ located with the arc between $1$ and $\om$ including
$\U^+\setminus\{1\}$ (or $\U^-\setminus\{1\}$)possess total
multiplicity $n-p-q$. If $\om \neq -1$, all eigenvalues of $K$ are
in $\U\setminus\R$ and those in $\U^+\setminus\R$ (or
$\U^-\setminus\R$) are all Krein-negative (or Krein-positive)
definite. If $\om=-1$, it holds that $(-I_{2s})\diamond
N_1(-1,1)^{\diamond t}\diamond H\in \Om^0(\ga(\tau))$ for some
non-negative integers $s$ and $t$ satisfying $0\le s+t\le n-p-q$,
and some $H\in \Sp(2(n-p-q-s-t))$ satisfying $\sg(H)\subset
\U\setminus \R$ and that all elements in $\sg(H)\cap\U^+$ (or
$\sg(H)\cap\U^-$) are all Krein-negative (or Krein-positive)
definite.

3. The left equality of (\ref{mm}) holds for all $\om\in
\U\setminus\{1\}$ if and only if $I_{2p}\diamond N_1(1,-1)^{\diamond
(n-p)}\in\Om^0(\ga(\tau))$ for some integer $p\in [0,n]$. Especially
in this case, all the eigenvalues of $\ga(\tau)$ are equal to $1$
and $\nu_\ga=n+p\ge n$.

4. The right equality in (\ref{mm}) holds for some $\om\in
\U^+\setminus\{1\}$ (or $\U^-\setminus\{1\}$) if and only if
$I_{2p}\diamond N_1(1,1)^{\diamond r}\diamond K\in \Om^0(\ga(\tau))$
for some non-negative integers $p$ and $r$ satisfying $0\le p+r\le
n$ and $K\in \Sp(2(n-p-r))$ with $\sg(K)\in \U\setminus\{1\}$
satisfying the condition that all eigenvalues of $K$ located with
the arc between $1$ and $\om$ including $\U^+\setminus\{1\}$ (or
$\U^-\setminus\{1\}$)possess total multiplicity $n-p-r$. If $\om
\neq -1$, all eigenvalues of $K$ are in $\U\setminus\R$ and those in
$\U^+\setminus\R$ (or $\U^-\setminus\R$) are all Krein-positive (or
Krein-negative) definite. If $\om=-1$, it holds that
$(-I_{2s})\diamond N_1(-1,1)^{\diamond t}\diamond H\in
\Om^0(\ga(\tau))$ for some non-negative integers $s$ and $t$
satisfying $0\le s+t\le n-p-r$, and some $H\in \Sp(2(n-p-q-r-t))$
satisfying $\sg(H)\subset \U\setminus \R$ and that all elements in
$\sg(H)\cap\U^+$ (or $\sg(H)\cap\U^-$) are all Krein-positive (or
Krein-negative) definite.

5. The right equality of (\ref{mm}) holds for all $\om\in
\U\setminus\{1\}$ if and only if $I_{2p}\diamond N_1(1,1)^{\diamond
(n-p)}\in\Om^0(\ga(\tau))$ for some integer $p\in [0,n]$. Especially
in this case, all the eigenvalues of $\ga(\tau)$ are equal to $1$
and $\nu_\ga=n+p\ge n$.

6. Both equalities of (\ref{mm}) holds for all $\om\in
\U\setminus\{1\}$ if and only if $\ga(\tau)=I_{2n}$.
 }

\noindent{\bf Acknowledgements} Part of the work was finished during
the author's visit at University of Michgan, he sincerely thanks
Professor Yongbin Ruan for his invitation and the Department of
Mathematics of University of Michigan for its hospitality.

\bibliographystyle{abbrv}

\end{document}